\documentclass[opre,nonblindrev]{informs3arXiv} 
\pdfoutput=1	
\usepackage{pdflscape}
\usepackage{amsmath}
\usepackage{amssymb}
\usepackage{multirow}
\usepackage{setspace}
\usepackage{algorithm}
\usepackage[flushleft]{threeparttable}
\usepackage{algorithmicx}
\usepackage{rotating}
\usepackage{subfigure}
\usepackage{amsfonts}
\usepackage{xcolor}
\usepackage{color}
\usepackage{pdflscape}
\usepackage{multirow}
\usepackage{graphicx}
\usepackage{calc}
\usepackage{verbatim}
\usepackage{mathtools}
\usepackage{tikz}
\usepackage{epstopdf}
\usepackage{tablefootnote}
\usepackage{algorithm}
\usepackage[noend]{algpseudocode}

\makeatletter
\def\BState{\State\hskip-\ALG@thistlm}
\makeatother

\definecolor{presentBlue}{RGB}{51,51,178}
\definecolor{baseBlue}{RGB}{0,77,77}
\definecolor{varBlue1}{RGB}{31,118,118}
\definecolor{varBlue2}{RGB}{13,97,97}
\definecolor{varBlue3}{RGB}{0,54,54}
\definecolor{varBlue4}{RGB}{0,30,30}

\definecolor{presentGreen}{RGB}{0,102,0}
\definecolor{varGreen1}{RGB}{41,158,41}
\definecolor{varGreen2}{RGB}{18,129,18}
\definecolor{varGreen3}{RGB}{0,72,0}
\definecolor{varGreen4}{RGB}{0,102,0}

\definecolor{baseRed}{RGB}{128,0,0}
\definecolor{varRed1}{RGB}{197,51,51}
\definecolor{varRed2}{RGB}{161,22,22}
\definecolor{varRed3}{RGB}{90,0,0}
\definecolor{varRed4}{RGB}{50,0,0}

\definecolor{presentGray}{RGB}{63,63,63}

\newcommand{\bA}{ \mathbf{A} }

\newcommand{\bc}{ \mathbf{c} }

\newcommand{\bd}{ \mathbf{d} }
\newcommand{\beee}{ \mathbf{e} }

\newcommand{\bff}{\mathbf{f} }

\newcommand{\bI}{ \mathbf{I} }

\newcommand{\bp}{ \mathbf{p} }

\newcommand{\bx}{ \mathbf{x} }

\newcommand{\by}{ \mathbf{y} }

\newcommand{\bzero}{ \mathbf{0} }

\newcommand{\balpha}{ \boldsymbol{\alpha} }

\newcommand{\bOmega}{ \boldsymbol{\Omega} }

\OneAndAHalfSpacedXI 

\usepackage{endnotes}
\let\footnote=\endnote
\let\enotesize=\normalsize
\def\notesname{Endnotes}%
\def\makeenmark{$^{\theenmark}$}
\def\enoteformat{\rightskip0pt\leftskip0pt\parindent=1.75em
  \leavevmode\llap{\theenmark.\enskip}}

\usepackage{natbib}
 \bibpunct[, ]{(}{)}{,}{a}{}{,}%
 \def\bibfont{\small}%
 \def\bibsep{\smallskipamount}%
 \def\bibhang{24pt}%
\TheoremsNumberedThrough     
\ECRepeatTheorems

\EquationsNumberedThrough    

\begin{document}


\RUNAUTHOR{Boutilier and Chan}

\RUNTITLE{Ambulance Emergency Response Optimization}

\TITLE{Ambulance Emergency Response Optimization in Developing Countries}

\ARTICLEAUTHORS{%
\AUTHOR{Justin J. Boutilier, Timothy C.Y. Chan}
\AFF{Department of Mechanical and Industrial Engineering, University of Toronto, Toronto, Ontario M5S 3G8, Canada \\ \EMAIL{j.boutilier@mail.utoronto.ca} \EMAIL{ tcychan@mie.utoronto.ca}}}


\ABSTRACT{The lack of emergency medical transportation is viewed as the main barrier to the access and availability of emergency medical care in low and middle-income countries (LMICs). In this paper, we present a robust optimization approach to optimize both the location and routing of emergency response vehicles, accounting for uncertainty in travel times and spatial demand characteristic of LMICs. We traveled to Dhaka, Bangladesh, the sixth largest and third most densely populated city in the world, to conduct field research resulting in the collection of two unique datasets that inform our approach. This data is leveraged to estimate demand for emergency medical services in a LMIC setting and to predict the travel time between any two locations in the road network for different times of day and days of the week. We combine our prediction-optimization framework with a simulation model and real data to provide an in-depth investigation into four policy-related questions. First, we demonstrate that outpost locations optimized for weekday rush hour lead to good performance for all times of day and days of the week. Second, we find that the performance of the current system could be replicated using one-third of the current outpost locations and one-half of the current number of ambulances. Lastly, we show that a fleet of small ambulances have the potential to significantly outperform traditional ambulance vans. In particular, they are able to capture approximately three times more demand while reducing the median average response time by roughly 10-18\% over the entire week and 24-35\% during rush hour due to increased routing flexibility offered by more nimble vehicles on a larger road network. Our results provide practical insights for emergency response optimization that can be leveraged by hospital-based and private ambulance providers in Dhaka and other urban centers in developing countries.
}


\maketitle

\section{Introduction}\label{Intro}

Time-sensitive medical emergencies are a major health concern in low and middle income countries (LMICs), comprising one third of all deaths \citep{Razzak2002}. Examples of such emergencies include cardiac arrest, motor vehicle accidents, and maternal health issues such as childbirth. Over the last decade, researchers and international organizations have stressed the need for increased focus on emergency medical care in LMICs \citep{UN2010,WHO2013}. In particular, the 66th World Health Assembly passed a resolution (60.22) that \emph{``recognizes the necessity of evidence-based approaches to development of emergency care and asks WHO to promote emergency medicine research''} \citep{WHA6022,Anderson}. However, despite widespread evidence that emergency medical care in LMICs save lives \citep{Sodemann1997, Schmid2001}, poor access and availability continues to be a major problem \citep{Kobusingye2005, Levine2007} with the lack of emergency medical transportation noted as being the main barrier \citep{Lungu2001, Macintyre1999}.

Optimizing the transport of emergency patients in urban centers in LMICs comes with unique challenges that are not present in high-income countries. First and foremost, traffic can be extremely unpredictable, and route disruptions caused by political demonstrations or extreme congestion occur regularly \citep{Jain2012, Pojani2015}. Second, it is not the norm, and often not possible due to congestion, for motorists to yield for emergency vehicles. As a result, route optimization (and vehicle outpost location, by extension) becomes a critical component for improving emergency vehicle response times. Third, LMICs generally do not have historical emergency call data that can be used to forecast future emergency demand. In fact, most LMICs do not have a centralized emergency response system, so the prospect of collecting a large, high-quality dataset is itself a major challenge. Together, these challenges lead to a high degree of uncertainty in both travel times and spatial demand. The nature of these uncertainties directly impacts any modeling approach, which must be compatible with ``small data'' environments characteristic of LMICs.

In this paper, we develop a robust optimization approach to optimize both the location and routing of emergency response vehicles, accounting for uncertainty in travel times and spatial demand characteristic of LMICs. We traveled to Dhaka, Bangladesh, the sixth largest and third most densely populated city in the world, to conduct field research resulting in the collection of two unique datasets that inform our approach. First, we obtained a field dataset that includes patient travel data associated with several thousand hospital arrivals. This data, acting as a proxy for historical call data available in all modern, high-income countries, is leveraged to develop a framework for estimating emergency medical services incidents in a LMIC setting. Second, we equipped five vehicles with custom-built GPS devices that recorded their time and location over a period of $30$ days, allowing us to understand traffic and road network characteristics in Dhaka. We then developed a machine learning model that uses the GPS data, along with census data, to predict the travel time between any two locations in the road network for different times of day and days of the week. For both demand and travel times, our predictions are leveraged to create data-driven uncertainty sets that are input into our robust location-routing model. Overall, our paper highlights the opportunity to creatively combine optimization with machine learning to solve a challenging emergency response problem in a resource-limited setting.

Like many urban centers in developing countries, Dhaka does not have a fleet of ambulances that form a centralized emergency response system. Instead, patients use a variety of transportation modes to reach hospitals in emergencies, including rickshaws, auto-rickshaws (i.e., three-wheeled motorcycles), private cars, and private or hospital-based ambulance services. Our modeling framework is well-suited to handle different transportation modes, which are accounted for via differences in road network connectivity according to vehicle type. Smaller and more nimble vehicles can traverse roads that larger vehicles cannot access. Therefore, the consideration of transportation mode affects the ultimate computational tractability of our models. In this paper, we focus on traditional \emph{van ambulances} and the locally inspired \emph{small ambulances}, which are based on three-wheeled motorcycles that have platforms that can be used for patient transport. In Dhaka (and many developing countries), most traditional van ambulances lack advanced medical equipment and are not typically staffed by paramedics, meaning that small ambulances are essentially equivalent from a medical equipment standpoint. Small ambulances have been recently proposed in Bangladesh \citep{Wadud2017}, but are not yet implemented and their potential impact on response times and patient outcomes has not been studied in the scientific literature.

Ambulance services in Dhaka are currently decentralized, meaning there are both private ambulance service providers, which are for-profit businesses, and ambulance fleets that belong to hospitals. Both types of organizations are incentivized to increase the number of patient transports they make, but lack appropriate decision support tools to optimize their operations. For example, hospitals do not currently strategically preposition their ambulances in the city, but rather position their entire fleet at the hospital. \cite{Savas1969} demonstrated the potential improvements over a similar hospital-based strategy in New York City. Therefore, private or hospital-based ambulance services are natural knowledge users of our research. Until recently, contact information for these services was also decentralized and unique to each provider, providing significant access challenges for patients. However, in December 2017, Bangladesh introduced the first centralized emergency services number ``999'' \citep{DhakaT}. The insights derived from our results can inform government policy on how to build a centralized emergency response system and aid non-government organizations to determine how to best position emergency response vehicle outposts. In particular, we use our real data and a simulation framework to answer four policy-related questions and derive practical insights for emergency response optimization in Dhaka and other LMICs: 

\begin{enumerate}
\item \textbf{Should different outpost locations be used for different times of day? (Section~\ref{TimeofDay})} In some high-income countries, ambulance locations are adjusted spatiotemporally throughout the day, but does that value persist in LMICs?
\item \textbf{What performance improvements are possible by optimizing outpost locations? (Section~\ref{AmbLocExps})} How different would a centralized, optimized system be from the current situation where ambulances are parked at hospitals? How does repositioning outpost locations compare to adding new locations?
\item \textbf{How much can the system be improved by using small ambulances? (Section~\ref{CNGvalue})} Can small ambulances capture additional demand that is currently unserved (or under-served) by existing van ambulances? What is the potential value of increased routing flexibility offered by small ambulances given their ability to traverse smaller roads in the network that are inaccessible to vans?
\item \textbf{How important is it to consider uncertainty when designing an emergency response network? (Section~\ref{RobustValue})} What is the performance improvement of our robust optimization model compared to a deterministic model? How does our robust approach compare to the perfect information case?
\end{enumerate}

The problem of optimizing emergency vehicle response has historically been cast as a facility location problem \citep{Toregas1971}. Although the facility location literature is rich, there is no unified framework for optimizing emergency vehicle response under both edge-based travel time uncertainty and demand uncertainty \citep{Ahmadi2017}. A key distinction between this paper and previous work is how we model travel time uncertainty. Our model provides a general edge-based framework for travel time uncertainty, whereas previous research has focused on modeling travel time uncertainty using a path-based approach \citep{Snyder2006}. Edge-length uncertainty is critical for our model because many of the underlying causes of travel time uncertainty in Dhaka (e.g., intersections without signal control, floods, strikes, etc.) impact small subsets of edges as opposed to the entire path. Uncertainty on individual edges can affect multiple routes and must be accounted for during optimization. Our routing problem is effectively a robust shortest path problem and, depending on how we model edge-length uncertainty, is equivalent to a regularized shortest path problem. The equivalence between robustness and regularization has been noted in domains such as regression \citep{Bertsimas2017}, but has not been previously demonstrated for the shortest path problem.

Overall, we use the aforementioned challenges faced by LMICs and gaps in the facility location literature to motivate the development of a novel location-routing model that is tailored for emergency response optimization in developing urban centers. We make the following contributions:

\begin{itemize}
\item We develop a novel edge-based reformulation of the classical path-based $p$-median problem. The $p$-median problem seeks to locate $P$ facilities relative to a set of demand nodes such that the total demand weighted distance to all demand nodes is minimized \citep{Hakimi1964, Hakimi1965}. This reformulation forms the foundation of a two-stage robust optimization model that considers both uncertain edge lengths (travel time) and node weights (demand). Our approach generalizes previous emergency facility location models based on the $p$-median architecture and provides a unified framework for emergency response optimization under travel time and demand uncertainty that is suitable for LMICs. We develop several approaches to solve our model. First, we develop an equivalent single-stage mixed-integer linear optimization problem. Second, we develop an exact scenario (i.e., row and column) generation algorithm that can improve the solution time by several orders of magnitude. For application to large-scale problems representative of the real road network in Dhaka, we develop a novel heuristic algorithm by extending a state-of-the-art $p$-median heuristic to work with edge-length uncertainty (Section~\ref{Opt}). All theorem proofs can be found in the Electronic Companion.
\item We develop a methodology to predict emergency demand spatially for urban centers without historical demand data by decomposing demand into components that can be estimated using census data and a regularized logistic regression model. This approach represents the first attempt to predict emergency demand in a developing urban center. Our complete dataset, including census, survey, and hospital location data is unique because, to the best of our knowledge, hospital arrival surveys and patient travel data have never been collected together previously in any LMIC (Sections~\ref{DemandEst} and~\ref{EmgDemand}).
\item We develop and compare several machine learning models to predict travel time on the Dhaka road network by time of day and day of week, using a dataset of vehicle trips collected by our custom-made GPS devices. We find that a random forest model performs the best, with a $43.3-64.2\%$ improvement in prediction accuracy over several baseline approaches. This paper is the first to use real travel time data from a LMIC for optimization (Sections~\ref{TT} and~\ref{Apd:TT}).
\item Using a simulation framework and real data from Dhaka, we provide an in-depth investigation into the four policy-related questions posed above (Section~\ref{DhakaExp}):
\begin{enumerate}
\item In contrast to developing countries where researchers have estimated performance improvements from repositioning ambulances according to the time and day, there is little to gain in Dhaka by optimizing outpost locations spatiotemporally. Instead, using outpost locations optimized for weekday rush hour leads to good performance for all times of day and days of the week. 
\item A centralized network designed from a clean slate can replicate the performance of the current system using roughly one-half of the ambulances and one-third of the outpost locations currently in use.
\item A fleet of small ambulances has the potential to significantly outperform traditional van ambulances. In particular, they can capture over three times the demand as van ambulances while reducing the median average response time by roughly 10-18\% over the entire week and 24-35\% during rush hour. This gain requires emergency response providers to tailor outpost locations specifically for small ambulances, instead of locating them at outposts optimized for traditional van ambulances.
\item Our robust solutions can reduce the median and worst-case response times by up to 33\% and 45\%, respectively, compared to a deterministic solution that does not take uncertainty into account. Furthermore, the performance of the robust solution is comparable to a solution that has access to perfect information on the uncertainty. 
\end{enumerate}
\end{itemize}

\section{Literature review}\label{LitRev}

Our work is related to three major streams of literature: 1) demand prediction in the context of emergency response optimization, 2) vehicle travel time prediction, and 3) facility location.

\subsection{Demand prediction}\label{Lit:DemPred}

While most papers use historical emergency call data as a direct estimate for future demand, a growing and more relevant body of literature uses that data to develop machine learning models that can predict future demand. Early approaches considered only spatial demand, using multiple linear regression to relate the magnitude of demand for ambulances with population and other socio-economic factors (e.g., \citealt{Schuman1977,Kamenetzky1982}). Key covariates can be summarized into three main groups: measures of population (e.g., household size), measures of economic status (e.g., employment rate, poverty level), and measures of social status (e.g., literacy rate, marriage rate). Temporal-only approaches were developed to forecast emergency calls at various time scales, including daily \citep{Baker1986}, multi-hour blocks \citep{Trudeau1989}, and hourly \citep{Channouf2007,Matteson2011}. Finally, there exist methods to predict future emergency demand at fine spatiotemporal resolutions \citep{Setzler2009,Zhou2015,Zhou2016}.

The aforementioned approaches rely on granular historical call data to train prediction models. High-income countries tend to be data-rich, so research efforts have focused on advanced demand prediction techniques using this abundant and granular data. However, in most LMICs, historical call data is not available \citep{Bradley2017}. In this paper, we develop a new approach that does not use historical call data and instead makes use of the limited spatiotemporal data available in many LMICs. 

\subsection{Travel time prediction}\label{Lit:TTPred}

Research on predicting edge-based travel times for ambulances has focused on developing non-linear relationships between travel time and distance \citep{Kolesar1975, Budge2010,Hofleitner2012a, Hofleitner2012b,Westgate2016}. However, almost all prior research depends directly on the availability of historical emergency transport data collected by a centralized system, which typically does not exist in LMICs.

In recent years, machine learning approaches that leverage decentralized travel time data have gained popularity and demonstrated superior prediction accuracy for regular vehicle travel time estimation \citep{Vlahogianni2014}. In contrast to ambulances, regular vehicle travel times are highly dependent on the time of day and the day of the week \citep{Kok2012,Woodard2017}. Travel times for emergency vehicles and regular vehicles are similar in LMICs because road users do not yield for ambulances. As a consequence, we employ a general travel time prediction approach similar to that of \cite{Zhang2015}, who use a random forest model that accounts for distance, time of day, and day of week. We extend their model by incorporating demographic and geographic characteristics for the origin and destination nodes, which encodes spatial information about the trip.

\subsection{Facility location}\label{Lit:FL}

Facility location is a very well-studied field and we provide only a brief review of the relevant literature. For a general review of facility location, please see \cite{Owen1998} or \cite{Melo2009}, and for a comprehensive review of facility location in the context of emergency medical services, please see \cite{Li2011}, \cite{Basar2012}, or \cite{Ahmadi2017}.

\subsubsection{Emergency response.}\label{Lit:EmergResp}

Facility location models have been applied extensively to emergency medical services location problems with the majority of previous research focusing on ambulances. There have been many papers that investigate ambulance response optimization in \emph{urban areas in high-income countries} (e.g., \citealt{Brandeau1986, Ingolfsson2008}), in \emph{rural areas in high-income countries} (e.g., \citealt{Adenso1997, Chanta2014}), and in \emph{rural areas in LMICs} (e.g., \citealt{Bennett1982, Eaton1986}). However, there have been only a few papers that consider \emph{urban areas in LMICs} \citep{Fujiwara1987, Basar2011, Salman2015, Zhang2015}, and they differ from our work in several important aspects. First, these papers focus on upper-middle-income countries (China, Thailand, and Turkey), whereas we focus on a low-income country (Bangladesh). Second, previous urban ambulance response optimization research, including the papers listed above, has focused exclusively on regions that already have a centralized ambulance system. In contrast, our paper is the first to focus on a developing urban center without an existing ambulance system, which leads to new policy questions not considered in areas with an existing system.

\subsubsection{Demand and travel time uncertainty.}

Demand uncertainty has received significant attention in general location-allocation problems (e.g., \citealt{Shen2003, Atamturk2007, Baron2011}) as well as in the specific context of ambulance response optimization \citep{Beraldi2004, Beraldi2009, Noyan2010}. The ambulance-specific papers all use chance constraints to model uncertain demand, whereas we employ a scenario-based approach that integrates a prediction model trained with our field data.

Travel time uncertainty in the context of ambulance response optimization has been focused on path-length uncertainty \citep{Ingolfsson2008, Berman2013, Ghani2017}. For networks with edge-length uncertainty, previous research has focused on the $1$-median problem \citep{Carson1990, Averbakh2003} and networks with special structure \citep{Mirchandani1979, Mirchandani1980}. We are the first to investigate edge-length uncertainty for the general $p$-median problem applied to ambulance response optimization.

Nearly all previous literature on combining both edge-length and node-weight uncertainty has focused on the special case of the $1$-median problem \citep{Chen1998,Vairaktarakis1999}, whereas we develop a methodology for the general $p$-median problem under uncertainty. The study by \cite{Serra1998}, which considers the $p$-median problem with both uncertain \emph{path lengths} and node weights, is the closest to our work. In contrast, we consider uncertain \emph{edge lengths} and node weights, which can be interpreted as a generalization of their model.

\subsubsection{Ambulance repositioning.}\label{Lit:AmbRepo}

Ambulance repositioning has received significant attention in the emergency response literature \citep{Brotcorne2003,Saydam2013,Nasro2018}. Repositioning strategies are often motivated by temporal changes in spatial demand and coverage gaps caused by busy vehicles. Real-time repositioning, which seeks to preposition ambulances in real time to better respond to future calls, leverages projected demand patterns and GPS-based ambulance location data. Repositioning strategies combined with dispatching decisions can also be used to mitigate system uncertainty \citep{Ena2018}. However, in many LMICs, real-time repositioning strategies are unrealistic because there is no centralized emergency response system to manage the real-time repositioning decisions.

Static ambulance repositioning is a simplified version of real-time repositioning that focuses on allocating ambulances to pre-selected outposts according to shift schedules, times-of-day, or the number of available ambulances \citep{Alanis2013,VanB2016, SUD2016, VanB2017}. Although static approaches are typically less effective than real-time strategies \citep{Maxwell2010}, they are easy to implement and manage. For example, compliance tables can be used to inform ambulance providers which outpost locations should be used for specific times of day or when there are only a certain number of ambulances available. We investigate the value of static repositioning in Section~\ref{TimeofDay}, which is motivated by changing demand patterns and the impact of changing traffic patterns on travel times \citep{Schmid2010}. While traffic is less of a concern in high-income countries, emergency vehicles typically face the same traffic conditions as regular road users in LMICs since other vehicles do not (or cannot) yield to ambulances.


\section{Optimization approach}\label{Opt}

We develop a two-stage robust optimization model to determine emergency response vehicle outpost locations. The outpost locations are determined based on how vehicles will be routed from the outpost to demand points (second stage), considering uncertainty in both demand and travel times.

We begin by introducing a novel edge-based location model that we prove to be equivalent to the classical $p$-median model. The advantage of our model is that it can handle edge-length uncertainty. Next, we introduce our models of uncertainty for emergency demand and travel times. Finally, we develop and compare several solution approaches.

\subsection{Network flow formulation}\label{Opt:NFF}

Let the road network be represented as the directed graph $G=(\mathcal{N}, \mathcal{E})$. Let $|\mathcal{N}|=n$, $|\mathcal{E}|=m$, and $\bA$ denote the $n\times m$ node-arc incidence matrix. Let $\bc$ denote the vector of edge lengths (i.e., travel times) and $\bd$ denote the vector of node weights (i.e., demand in terms of average annual emergency transports required). Let $\balpha$ denote the supply available at each potential facility (i.e., number of trips that can be made from each outpost per year) and $\bOmega $ represent the $n\times n$ diagonal matrix whose entries are the $n$ elements in $\balpha$. We use $P$ to represent the number of outposts to be located and $\beee$ to denote the vector of all ones. The decision variable representing the vector of flows along each edge is denoted by $\bff$ (i.e., how many trips occur on each edge annually). The outpost location variable is given by $\by\in\{0,1\}$ where $1$ indicates an outpost is located at node $i\in\mathcal{N}$. Note that all defined vectors are column vectors. In vector form (see~\ref{Apd:ProofThm1} for the non-vectorized version), our deterministic network flow formulation (\textbf{NFF}) is:
\allowdisplaybreaks
\begin{equation} \label{NFF1}
\begin{aligned}
\textbf{NFF:}\quad\quad\underset{\by,\bff}{\mathrm{minimize}} \hspace*{1em} & \bc'\bff \\
\mbox{subject to}\hspace*{1em} & \beee'\by = P, \\
&\bA\bff \leq \bOmega\by - \bd,\\
& \bff\geq \bzero, \\
& \by\in\{0,1\}^n.
\end{aligned}
\end{equation}

The second constraint accounts for supply nodes, ensures that all demand is met, and allows for transshipment flow. In scalar form, the constraint can be written as:
$$\sum_{j\in O(i)} f_{ij} - \sum_{j\in I(i)} f_{ji} \leq \alpha_i y_i - d_i, \forall i \in N,$$
where $I(i) =\{ j\in N | (j,i)\in \mathcal{E}\}$ and $O(i)=\{ j\in N | (i,j)\in \mathcal{E}\}$. If $y_i=1$, then node $i$ becomes a source node that produces up to $\alpha_i-d_i$ trips per year. If $y_i=0$, then node $i$ becomes a demand node and the constraint reduces to $\sum_{j\in I(i)} f_{ji} - \sum_{j\in O(i)} f_{ij} \geq d_i$. This ensures that at least $d_i$ trips flow into node $i$ (thereby satisfying demand), but also allows for trips to flow into and out of node $i$ en route to another location.

To ensure that \eqref{NFF1} is feasible for any value of $P$, we require the following assumption.

\begin{assumption}{$\alpha_i \geq\sum\limits_{i=j}^n d_j, \; \forall i \in \mathcal{N}$.}\end{assumption}

This assumption states that each outpost has enough capacity to service the entire system (i.e., all demand nodes) and to ensure feasibility, we set $\alpha_i = \sum\limits_{i=1}^n d_i, \; \forall i \in \mathcal{N}$. We do not consider queuing in our model because our primary focus is to determine where to strategically locate emergency response outposts, rather than determining the total number of emergency response vehicles. However, we later evaluate the tactical performance of our solutions with respect to queuing and system congestion using a simulation model. Lemma~\ref{lem1} follows immediately from this assumption (proof omitted). 
\begin{lemma}\label{lem1}
There exists an optimal solution to \textbf{NFF} such that each demand node is assigned to exactly one outpost.
\end{lemma}

This result generally holds true for uncapacitated facility location models such as the $p$-median. Finally, using Lemma~\ref{lem1}, we can show the equivalence between \textbf{NFF} and the $p$-median problem.

\begin{theorem}\label{thm:NFFpmedequiv}
A solution is optimal for \textbf{NFF} if and only if it is optimal for the $p$-median problem.
\end{theorem}


The proof of Theorem~\ref{thm:NFFpmedequiv} provides a constructive approach to obtain an optimal solution of \textbf{NFF} given an optimal solution of the $p$-median problem, and vice versa. Mathematically, this approach provides a polynomial-time many-one reduction between the \textbf{NFF} and the $p$-median problem in both directions \citep{Post1944, Karp1972}.

\subsection{Robust optimization model}\label{Opt:RNFF}
In this section, we present our two-stage robust optimization model, considering both the travel times $\bc$ and demands $\bd$ as uncertain with $\mathcal{C}$ and $\mathcal{D}$ representing the corresponding uncertainty sets, respectively. Our general two-stage robust network flow formulation is:
\allowdisplaybreaks
\begin{equation} \label{NFFRED}
\begin{aligned}
\textbf{R-NFF:}\quad\quad\mathrm{\min_{\by}\;\max_{\bc\in\mathcal{C}, \bd\in\mathcal{D}}\;\min_{\bff}} \hspace*{1em} & \bc'\bff \\
\mbox{subject to}\hspace*{1em} & \beee'\by = P, \\
&\bA\bff \leq \balpha\bI\by - \bd,\\
& \bff\geq\bzero, \\
& \by\in\{0,1\}^n.
\end{aligned}
\end{equation}
The two-stage nature of our formulation is well-suited to the problem of emergency outpost location and vehicle routing. In the first stage, \textbf{R-NFF} determines the optimal outpost locations considering both $\bc$ and $\bd$ as uncertain. Intuitively, determining these locations is a high-level strategic decision that must be made under uncertainty, before demand or traffic are realized. Then, given the realized demand and travel time conditions, the second stage determines the optimal routes from the outposts to reach each demand point (i.e., patient location). Routing is a secondary decision that is used to inform the first stage location decision because the suitability of an outpost location is influenced by the route options emanating from that outpost.

\subsubsection{Demand uncertainty set ($\mathcal{D}$).}\label{Opt:DUnc}

To model uncertainty in emergency transport demand, we use a scenario-based uncertainty set. We use this approach to preserve tractability while still capitalizing on the richness of our demand predictions. For $N$ scenarios, the resulting uncertainty set is defined as $\mathcal{D}=\{\bd^1,\bd^2,...,\bd^{N}\}$, where the dimension of $\bd$ is equal to the number of nodes in the network. To generate the scenarios that form the uncertainty set, we employ a form of bootstrapping and simulate possible realizations of demand vectors using our framework from Sections~\ref{DemandEst} and~\ref{EmgDemand}. 

\subsubsection{Travel time uncertainty set ($\mathcal{C}$).}\label{Opt:TTUnc}

Uncertainty in travel time is modeled using an interdiction-based uncertainty set with an overall budget constraint \citep{Wood1993}. Intuitively, this set models an adversary who is adding traffic (i.e., increasing travel time) to the baseline traffic on each edge. The budget constraint restricts the total amount of travel time that can be added across the network. The mathematical formulation of this uncertainty set is
$\mathcal{C} =\left\{c_{ij}, (i, j)\in\mathcal{E} \;\bigg|\; c_{ij} = \hat{c}_{ij} + w_{ij}, \sum_{(i,j)\in\mathcal{E}}w_{ij}\leq B,  w_{ij}\geq 0, \forall (i,j)\in \mathcal{E} \right\}.$
We estimate the baseline travel time $\hat{c}_{ij}$ for each edge using the final random forest model from Section~\ref{TT:Preds}. In our numerical experiments, we perform a detailed sensitivity analysis on the budget $B$.

\subsection{Solution Algorithms}\label{Opt:SolAlg}

In this section, we present several methods to solve \textbf{R-NFF}. First, we show that there is an equivalent single-stage mixed-integer optimization model for \textbf{R-NFF}. Then, we present an exact row-and-column generation algorithm to solve this equivalent problem. Finally, for the integer master problem, we devise an efficient heuristic that is needed for large-scale instances. See \ref{CompExp} for a detailed numerical comparison of the solution times and optimality gaps between the mixed-integer model, exact solution algorithm, and heuristic solution algorithm.

\subsubsection{Equivalent mixed-integer optimization model.}\label{Opt:MILP}

First, we replicate $\bff$ for each of the scenarios in the demand uncertainty set ($\mathcal{D}$). Formally, we define $\bff^k$ as the flow decision variable for scenario $k=1,\dots,N$ and $\zeta^k$ to be the dual variable corresponding to scenario $k$ for the travel time uncertainty set constraint $\sum_{(i,j)\in\mathcal{E}}w^k_{ij}\leq B$ in $\mathcal{C}$. The flow variable $\bff^k$ corresponding to the limiting scenario for the first set of constraints in \eqref{NFFsolve} is an optimal flow vector for \eqref{NFFRED}.

\begin{theorem}\textbf{R-NFF} is equivalent to the following mixed-integer linear optimization problem:
\begin{equation}\label{NFFsolve}
\begin{aligned}
\mathrm{\underset{\by,t, \bff^k, \zeta^k}{minimize}}\quad&t \\
\mathrm{subject\;to}\quad & t\geq\hat{\bc}'\bff^k + \zeta^kB,\quad k = 1,...,N,\\
& \bA\bff^k\leq\balpha\bI\by-\bd^k, \quad k = 1,...,N,\\
& \bff^k \leq \zeta^k\beee, \quad k = 1,...,N,\\
& \bff^k \geq \bzero,\quad k = 1,...,N,\\
& \zeta^k  \geq 0, \quad k = 1,...,N,\\
& \beee'\by = P, \\
& \by\in\{0,1\}^n.
\end{aligned}
\end{equation}
\end{theorem}


Formulation~\eqref{NFFsolve} quickly becomes intractable as the number of scenarios increases and the size of the graph grows. We address these two challenges in the next two subsections. First, we develop a scenario generation algorithm that scales efficiently with the number of scenarios. Similar decomposition algorithms have been developed by \cite{Atamturk2007, Zeng2013, Gabrel2014} and \cite{Chan2017} for related two-stage problems. Second, we develop a heuristic to efficiently solve the master problem associated with the scenario generation approach.

\subsubsection{Scenario Generation.}\label{Opt:ScenGen}

Consider a subset of the demand scenarios $\mathcal{D}_{|S|} = \{\bd^1,\bd^2,...,\bd^{|S|}\} \subset \mathcal{D}$, where $S$ is an index set for the vectors in $\mathcal{D}_{|S|}$, and the corresponding relaxation of formulation \eqref{NFFsolve} with $\mathcal{D}_{|S|}$ in place of $\mathcal{D}$:
\begin{equation}\label{NFFsolveMP}\
\begin{aligned}
\textbf{R-NFF-MP:}\quad\quad\mathrm{\underset{\by,t, \bff^s, \zeta^s}{minimize}}\quad&t \\
\mbox{subject to}\quad & t\geq\hat{\bc}'\bff^s + \zeta^sB,\quad \forall s\in S,\\
& \bA\bff^s\leq\balpha\bI\by-\bd^s, \quad \forall s\in S,\\
& \bff^s \leq \zeta^s\beee, \quad \forall s\in S,\\
& \bff^s \geq \bzero,\quad \forall s\in S,\\
& \zeta^s  \geq 0, \quad \forall s\in S,\\
& \beee'\by = P, \\
& \by\in\{0,1\}^n.
\end{aligned}
\end{equation}

The relaxed master problem, \eqref{NFFsolveMP}, produces a lower bound on the optimal value of \eqref{NFFRED} that can be tightened by adding additional scenarios to the set $\mathcal{D}_{|S|}$. Given an optimal solution $\bar{\by}$ to \eqref{NFFsolveMP}, we solve the following sub-problem, which is a linear optimization problem, for every $\bd^k\in\mathcal{D}$:

\begin{equation}\label{NFFsolveSP}
\begin{aligned}
\textbf{R-NFF-SP-k:}\quad\quad Z^k_{SP}\;=\;\mathrm{\underset{\bff^k, \zeta^k}{minimize}}\quad&\hat{\bc}'\bff^k + \zeta^kB \\
\mbox{subject to}\quad &\bA\bff^k\leq\balpha\bI\bar{\by}-\bd^k,\\
& \bff^k \leq \zeta^k\beee, \\
& \bff^k \geq \bzero,\\
& \zeta^k  \geq 0.
\end{aligned}
\end{equation}

We choose the scenario $k^* \in \underset{k=1,...,N}{\mathrm{arg\, max}}\{Z_{SP}^k\}$ and add the decision variables $\bff^{k^*}$ and $\zeta^{k^*}$, plus their corresponding constraints, to \eqref{NFFsolveMP}. Hence, this approach generates both rows and columns. The scenario generation algorithm terminates when the optimal value of $\eqref{NFFsolveMP}$ is equal to $Z_{SP}^{k^*}$.

Finally, we comment on the structure of the subproblem \eqref{NFFsolveSP} and connect it to a stream of research that draws an equivalence between robust optimization and regularization. Since $\zeta^k$ is being minimized in \eqref{NFFsolveSP}, the constraint $\bff^k\leq\zeta^k\beee$ identifies the maximum value of $f^k_{ij}$ over all $(i,j)\in\mathcal{E}$. Thus, we can rewrite \eqref{NFFsolveSP} as (we drop the index $k$ for simplicity):

\begin{equation}\label{RSP}
\begin{aligned}
\mathrm{\underset{\bff}{minimize}}\quad&\hat{\bc}'\bff + B\|\bff\|_{\infty} \\
\mbox{subject to}\quad &\bA\bff\leq\balpha\bI\bar{\by}-\bd,\\
& \bff \geq \bzero.
\end{aligned}
\end{equation}

Formulation~\eqref{RSP} is a ``regularized'' shortest path problem. Without the term $B\|\bff\|_\infty$ in the objective, \eqref{RSP} is exactly a shortest path problem. The extra term balances finding the shortest path with minimizing the maximum flow along any edge, which is weighted by the budget $B$. In our application, larger values of $B$ correspond to higher levels of traffic uncertainty. Thus, for large $B$, an optimal solution to~\eqref{RSP} would prefer to spread out the flows (smaller maximum $f_{ij}$), forcing nature to expend more budget to ``lengthen'' multiple edges. Equivalently, if flows are concentrated on a few arcs, then nature has easy targets for adding traffic to cause maximal disruption. Our reformulation elucidates a clear connection between a robust shortest path problem and a regularized shortest path problem, similar to the way equivalences have been derived in regression \citep{Xu2010,Bertsimas2017}. For example, if we replace the constraint $\sum_{(i,j)\in\mathcal{E}}w_{ij}\leq B$ in $\mathcal{C}$ with $w_{ij}\leq B, \forall(i,j)\in\mathcal{E}$, then our subproblem is equivalent to a L1-regularized (lasso) problem.

\subsubsection{Master problem heuristic.}\label{Opt:Heur}

To solve the large-scale, real-world instances considered in our Dhaka experiments, we require a heuristic for the master problem, which is in essence a $p$-median problem. Although there are many heuristics that have been developed for the $p$-median problem, we cannot apply these algorithms directly because they are unable to handle edge-length uncertainty. Instead, we adapt the heuristic developed by \cite{Densham1992} for the classical $p$-median problem. This heuristic, designed for large-scale problems, leverages both the interchange heuristic proposed by \cite{Teitz1968} and the alternate heuristic proposed by \cite{Maranzana1964}. A key benefit of this type of algorithm is that it scales well with both the size of the graph and the number of facilities ($P$). In fact, our heuristic represents the first tractable approach to solving large-scale instances of location problems with edge-length uncertainty. Our approach involves three main phases. 

\textbf{Initialization phase.}
We initialize our algorithm by randomly selecting $P$ nodes to serve as initial outpost locations, encoded by $\bar{\by}$. We solve \eqref{NFFsolveSP} with this $\bar{\by}$ for every $\bd^k\in\mathcal{D}_{\mathcal{S}}$, and identify $k^* \in \arg\max_{k \in \mathcal{S}}\{Z_{SP}^k\}$, $\bff^{k^*}$, and $\zeta^{k^*}$. The corresponding cost of this solution is $\hat\bc'\bff^{k^*} + \zeta^{k^*}B$. An advantage of a random initialization phase is that our algorithm can be embedded in a meta-heuristic or a simple approach that considers multiple random starts. We investigate the impact of the number of random starts in our numerical experiments.

\textbf{Interchange phase.}
In the interchange phase, we randomly swap a current outpost location node with a candidate node that is not currently in the solution. The new objective value is calculated as before after solving \eqref{NFFsolveSP} for every $\bd^k\in\mathcal{D}_{\mathcal{S}}$. Swaps that reduce the objective value are accepted. We consider $\ell$ random interchanges per outpost location, where $\ell$ is a user-chosen parameter. 

\textbf{Alternate phase.} In the alternate phase, we use the incumbent solution from the interchange phase to partition the network into $P$ connected subgraphs that are disjoint from each other. Each subgraph contains exactly one outpost location and all demand nodes served by that outpost. We solve \eqref{NFFsolveMP} for $P=1$ (i.e., the robust $1$-median problem) on each subgraph to determine the optimal outpost location. We then re-combine all subgraphs and the new optimal outpost locations to obtain an updated set of outpost locations, $\bar{y}$, in the full network. We compute the cost of this solution as before, by solving \eqref{NFFsolveSP} for every $\bd^k\in\mathcal{D}_{\mathcal{S}}$. The alternate phase continues to partition and re-combine outpost locations until it has reached a local optimum. The algorithm then proceeds back to the interchange phase.

\textbf{Termination.}
The algorithm iterates between the interchange and alternate phases until a solution from the alternate phase is found that does not result in any swaps during the interchange phase. The algorithm terminates with a solution to a single instance of the master problem~\eqref{NFFsolveMP}.

\textbf{Integration with scenario generation algorithm.}
The returned solution from the heuristic either terminates the scenario generation algorithm (when the the optimal value of $\eqref{NFFsolveMP}$ is equal to $Z_{SP}^{k^*}$) or is used as input to the sub-problem~\eqref{NFFsolveSP}. 

\section{Application to Dhaka}\label{DhakaApp}

In this section, we outline the application of our methodology to Dhaka, Bangladesh. Section~\ref{TT:Roads} describes the road networks, Section~\ref{DemandEst} details our approach for estimating spatiotemporal demand for emergency transportation, Section~\ref{TT} outlines our travel time predictions, Section~\ref{Sim} presents our tactical simulation model, and Section~\ref{SolApproach} describes our experimental setup.

\subsection{Road networks}\label{TT:Roads}

\begin{figure}[t]
\centering
\subfigure[\ Ambulance \label{AN}]{
\includegraphics[width=.2\textwidth]{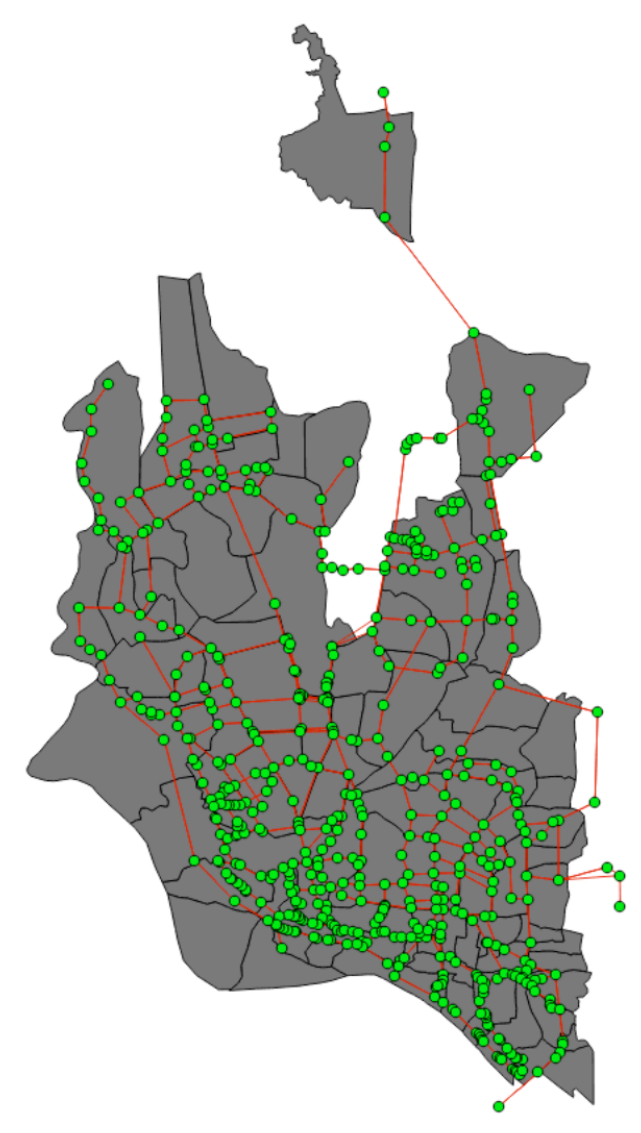}}
\hspace*{5em}
\subfigure[\ Complete \label{CN}]{
\includegraphics[width=.25\textwidth]{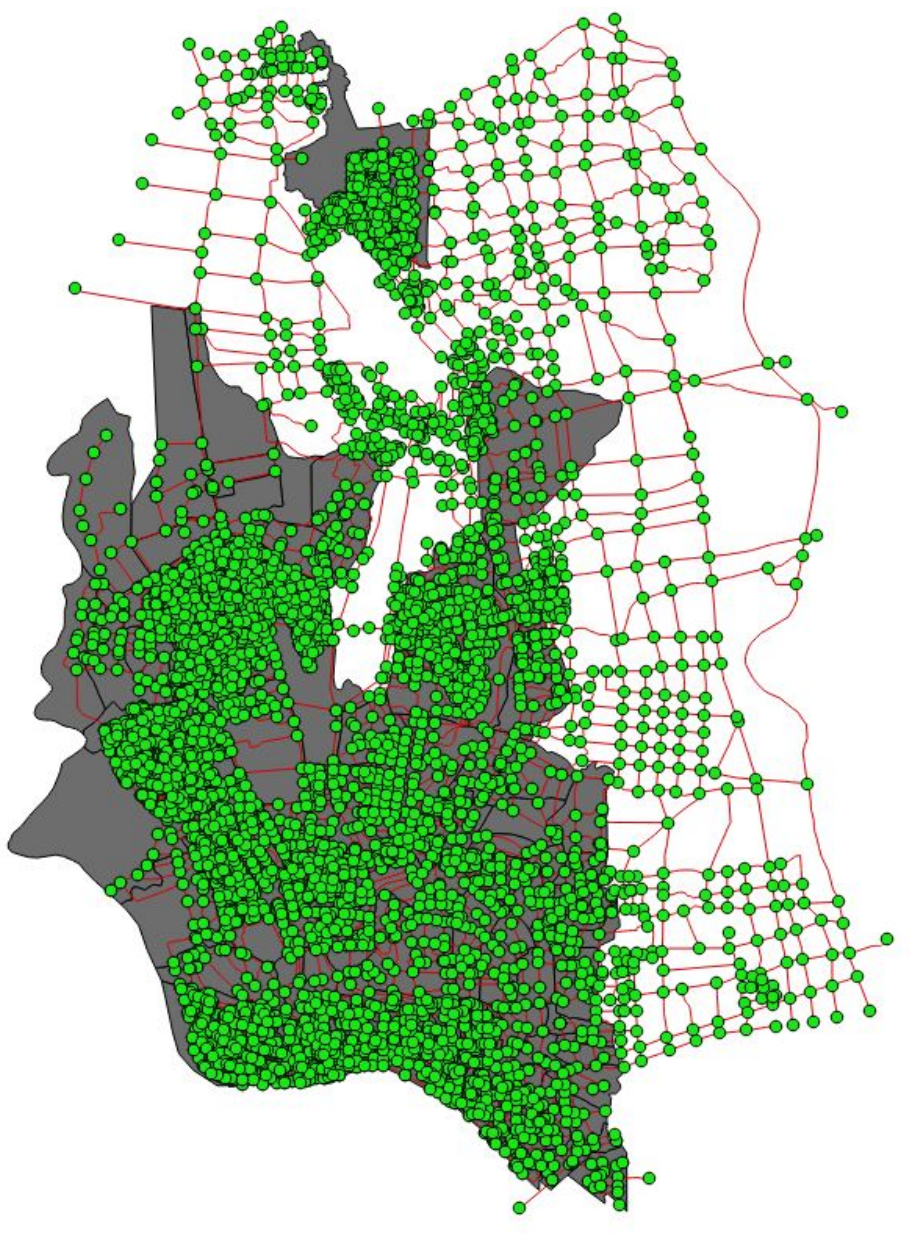}}
\caption{Two road networks overlaid on a ward map of Dhaka.\label{Roads}}
\end{figure}

We consider two different road networks in Dhaka. The first road network that we consider is the \emph{ambulance network}. In consultation with a transportation engineer in Dhaka and using a detailed map of the entire city, we determined exactly which roads are feasible for ambulance travel (many roads are too narrow for a van ambulance). The ambulance network has $530$ nodes and $1,280$ edges. A node is defined as the intersection of edges (i.e., roads). The second road network we consider is the \emph{complete network}. This network -- a superset of the ambulance network -- includes all roads ranging from large arterial roads to small alleyways that can only be traversed by small vehicles like rickshaws, motorcycles, and auto-rickshaws. The complete road network has $5,358$ nodes and $16,538$ edges. Figure \ref{Roads} displays both networks overlaid on Dhaka's $92$ wards. 

\subsection{Demand for emergency transportation}\label{DemandEst}

In this section, we outline our framework for estimating spatiotemporal demand for emergency transportation. We do not have data on the total number of emergency transports as we would in North America because Dhaka does not have a centralized emergency medical system. Instead, we propose a two step process that leverages the limited data at our disposal (see~\ref{ED:Descriptive} for a detailed description of our data). First, we provide a novel decomposition of a standard metric for emergency demand: the annual number of emergency trips from ward $w$ at time $\tau$ via mode $m$ (Section~\ref{annualtrips}). Second, we develop a simulation framework to estimate the precise time and location for each emergency transport (Section~\ref{tripsim}).

\subsubsection{Estimating the annual number of emergency trips.}\label{annualtrips}

We decompose the estimated annual number of emergency trips for each ward $w$, time of day $\tau$, and mode $m$, denoted by $d_{w,\tau,m}$, into three components: 
\begin{equation}\label{Demand0}
\begin{split}
d_{w,\tau,m} &= \xi \, n_{w,\tau} \, \delta_{w,m},
\end{split}
\end{equation}
where $\xi$ represents the average annual number of emergency trips per person, $n_{w,\tau}$ represents the population in ward $w$ at time $\tau$, and $\delta_{w,m}$ represents the proportion of emergency trips from ward $w$ that arrived via mode $m$.

Equation~\eqref{Demand0} suggests an approach to estimating $d_{w,\tau,m}$ by estimating its constituent terms $\xi$, $n_{w,\tau}$, and $\delta_{w,m}$. To do this, we consider two time periods (daytime ($D$) and nighttime ($N$)) and two modes of transport (van ambulance ($V$) and small ambulance ($S$)). We consider two time periods because emergency demand is known to follow a circadian rhythm \citep{Circ2013, McCormack2015}, meaning that demand is much higher during the day than at night. We consider two modes of transport because of the multi-modal nature of decentralized ambulance services in LMICs.

\begin{figure}[t]
\centering
\hspace*{1em}
\subfigure[\ Van ambulance - daytime \label{AmbPred1}]{
\includegraphics[width=.31\textwidth]{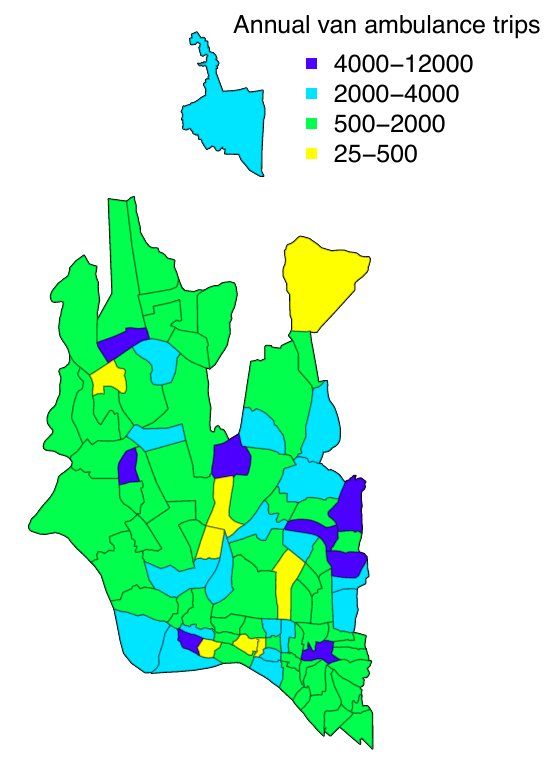}}
\hspace*{1em}
\subfigure[\ Van ambulance - nighttime \label{AmbPred2}]{
\includegraphics[width=.3\textwidth]{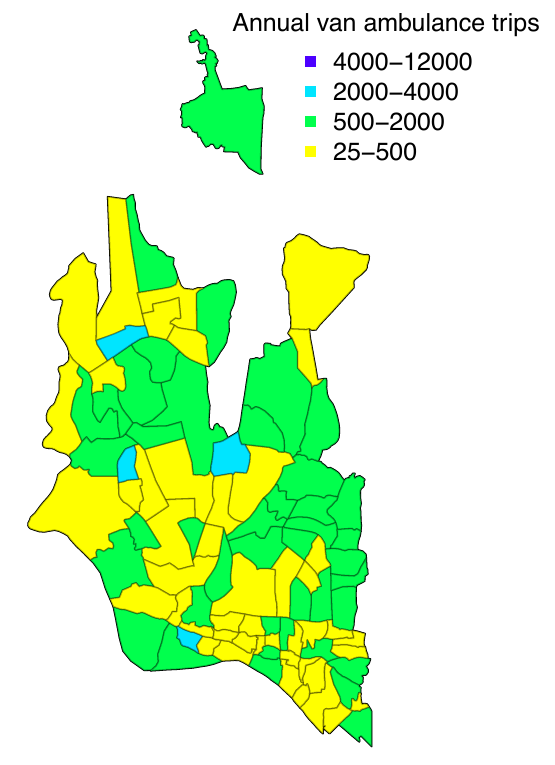}}
\hspace*{1em}

\subfigure[\  Small ambulance - daytime \label{CNGPred1}]{
\includegraphics[width=.285\textwidth]{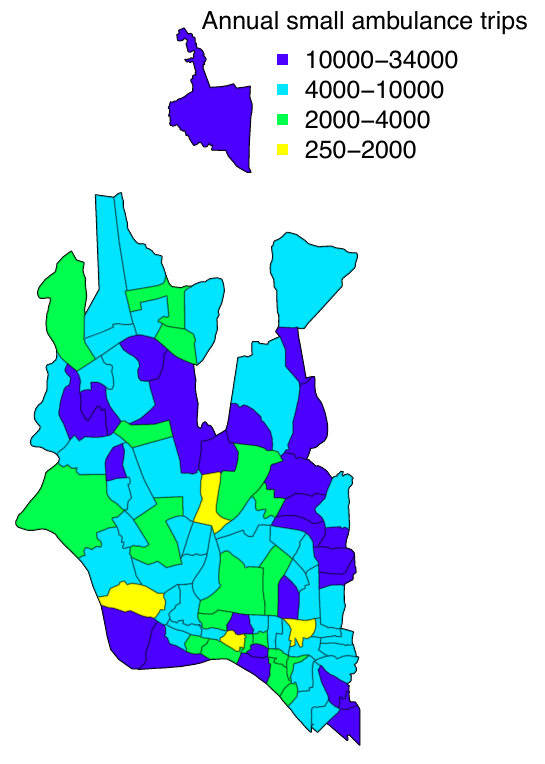}}
\hspace*{1em}
\subfigure[\ Small ambulance - nighttime  \label{CNGPred2}]{
\includegraphics[width=.29\textwidth]{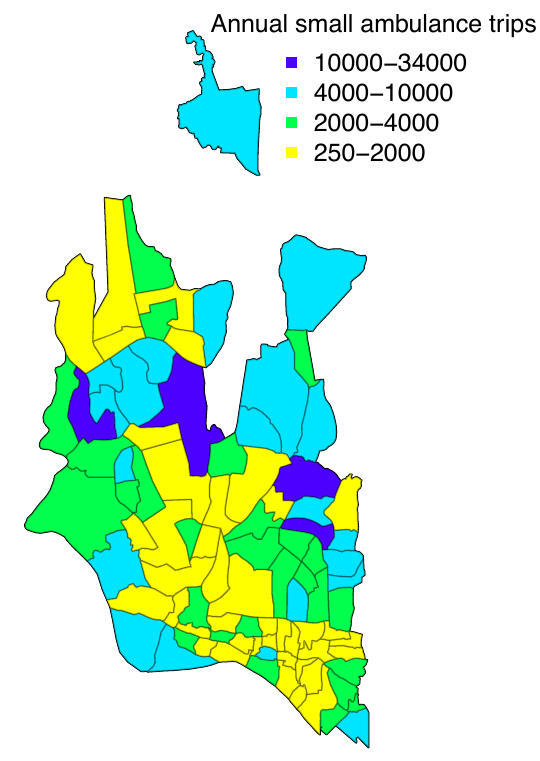}}
\caption{Expected annual number of van and small ambulance trips arising in each ward.\label{Finaldemand}}
\end{figure}

In total, there are 369 quantities to estimate: two per ward for population ($n_{w,\tau}$), two per ward for mode ($\delta_{w,m}$), and a single value for the average annual number of ED visits per person across the entire city ($\xi$). The estimation of these three sets of parameters are described in \ref{Apd:DemandEst1}, \ref{Apd:DemandEst2}, and \ref{Apd:DemandEst3}, respectively. Figure~\ref{Finaldemand} shows the final estimation for the expected annual number of daytime and nighttime trips arising from each ward, for both van and small ambulances.

We use our estimations to simulate scenarios for the uncertainty set described in Section~\ref{Opt:DUnc}. To do this, we assume that the population in each ward follows a triangle distribution on the interval between the estimated daytime and nighttime population, with a peak at the midpoint; we assume that $\xi$ follows a triangle distribution on the estimated interval $[0.23-0.46]$ (see \ref{Apd:DemandEst2} for details), where the peak occurs at $0.40$ for conservatism; and we assume that $\delta_{w,m}$ follows a truncated normal distribution with a mean equal to the predicted ward value (Figures EC.\ref{Pred1} and EC.\ref{Pred2}) and a standard deviation equal to the median error (see Figures EC.\ref{Err} and EC.\ref{Err2}).

The simulated demand vectors need to be adjusted so that the demand in each ward is spread proportionally to the road network nodes in that ward. In other words, we need to map the predicted demand based on the 92 wards to the $\sim$500 or $\sim$5,000 nodes in the ambulance and complete networks, respectively. To make this adjustment, we generate a fine grid of nodes spaced $25$m apart across all $92$ wards, resulting in over 200,000 grid nodes. We distribute the simulated demand in each ward uniformly among the grid nodes in that ward. Then, we assign each grid node and its corresponding demand to the closest road network node using Euclidean distance.

\subsubsection{Simulating spatiotemporal emergency trips.}\label{tripsim}

As written, \eqref{Demand0} provides a high-level approach to estimate the spatial and temporal heterogeneity in emergency demand. These annual estimations are useful for our optimization model because locating ambulance outposts is a high-level strategic decision that may be fixed for long time period. However, our simulation model, which evaluates the tactical performance of the optimization results, requires the exact time and node location for each emergency trip. To do this, we develop a novel procedure that maps annual ward-based demand to a fine spatiotemporal resolution.

\begin{figure}[t]
\centering
\includegraphics[width=.5\textwidth]{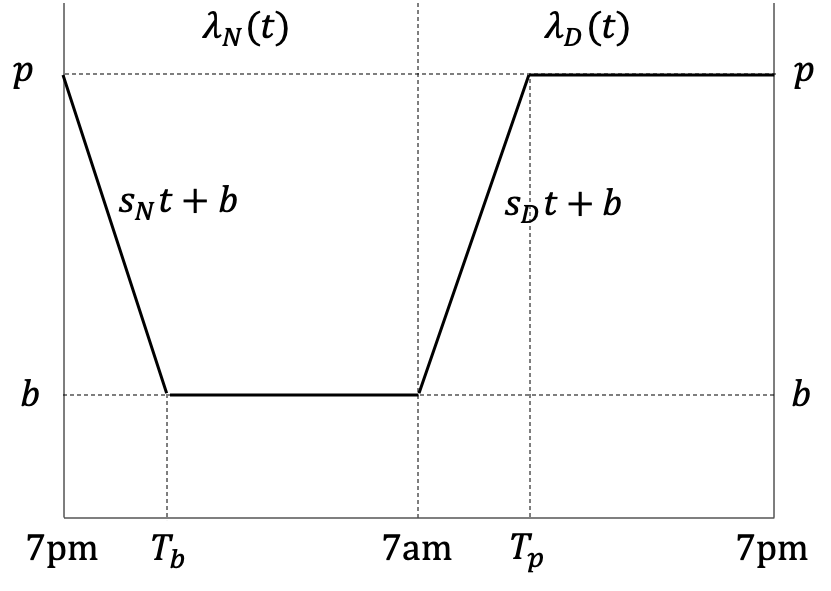}
\caption{A visualization of $\lambda_N(t)$ and $\lambda_D(t)$ (we drop the ward and mode indices). \label{PDF}}
\end{figure}

We approximate the time-dependent arrival rate for emergency demand with the piecewise linear function shown in Figure~\ref{PDF}. For each ward and mode of transportation, we partition the daily arrival rate function into separate daytime ($\lambda_{w,D,m}(t)$) and nighttime ($\lambda_{w,N,m}(t)$) components. We drop the $w$ and $m$ indices for the remainder of this section. We translate our annual estimates for emergency trips ($d_{\tau}$) into daily estimates ($\hat{d}_{\tau}$) by assuming that each day has an equal number of expected emergency trips (i.e., $\hat{d}_{\tau} = \frac{d_{\tau}}{365}$) \citep{Circ2013, McCormack2015}. We also assume that $T_p = 9$am and $T_b=11$pm, based on estimates from the literature \citep{Circ2013}. The boundary conditions of $\lambda_{D}(t)$ and $\lambda_{N}(t)$ are set to ensure continuity of the overall function (i.e., $b = b_{N} = b_{D}$ and $p = p_{N} = p_{D}$). Hence, for each ward and mode, there are four unknowns: $s_{D}, s_{N}, b,$ and $p$. 


\begin{figure}[t]
\centering
\includegraphics[width=.5\textwidth]{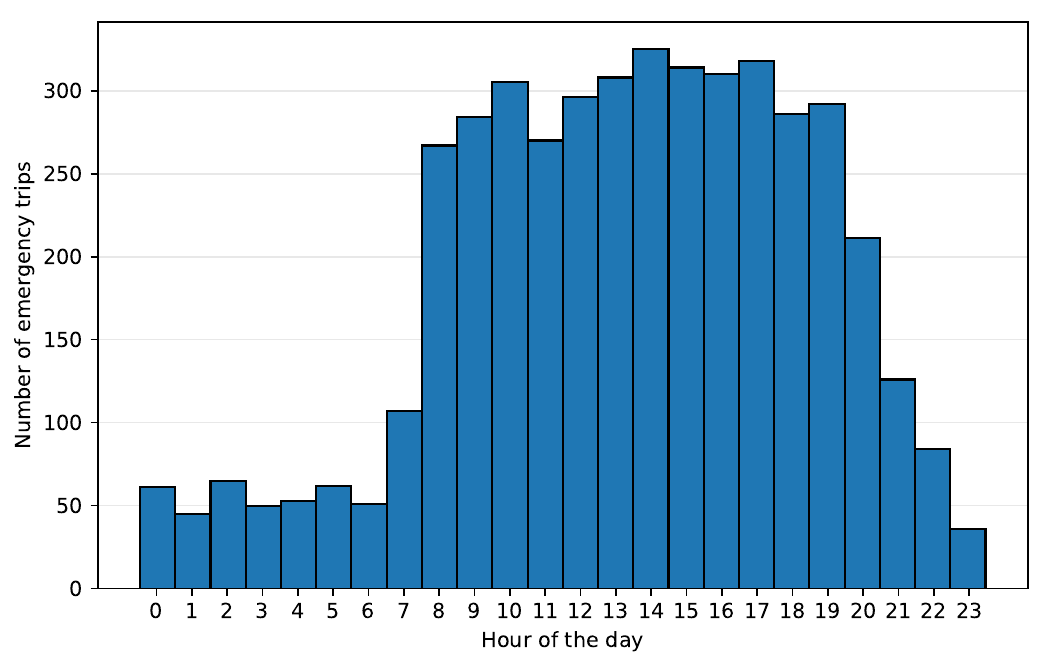}
\caption{A histogram displaying one week of van ambulance emergency calls summed across all wards and binned according to the hour of the day. \label{CallHist}}
\end{figure}

We obtain closed form solutions for the four parameters by leveraging the fact that $\hat{d}_{D} = \int_0^{12}\lambda_{D}(t)dt$ and $s_{D} = \frac{p-b}{T_p}$ (similar equations hold for nighttime). Once we determine the parametric form for the arrival rate function, we use the order statistic sampling method to generate the exact time for each emergency trip \citep{Cox1966, NHPP}:
\begin{enumerate}
\item Generate the number of emergency trips: $n\sim$ Poisson$(\hat{d}_{\tau})$
\item Independently generate $n$ random variates $t^{\prime}_1, t^{\prime}_2, \dots, t^{\prime}_n$ from the cumulative distribution function given by $F_{\tau}(t) = \frac{1}{\hat{d}_{\tau}}\int_0^{t}\lambda_{\tau}(t)dt$
\item Order $t^{\prime}_1, t^{\prime}_2, \dots, t^{\prime}_n$ and return the ordered times.
\end{enumerate}
We repeat these steps for each ward ($w$), time of day ($\tau$), and mode of transport ($m$). For each simulated 24-hour period, this procedure produces a series of times for each ward and mode that follow the distribution shown in Figure~\ref{PDF}. Figure~\ref{CallHist} displays one week of van ambulance emergency calls summed across all wards and binned according to the hour of the day. 

We map the times to nodes on the road network using a modified version of the procedure outlined in Section~\ref{annualtrips}, which maps the total ward demand to nodes in the road network using a finely spaced grid. In other words, each road network node captures the demand of $g$ grid nodes. For a given ward, we randomly assign each trip to a road network node with a probability corresponding to the proportion of grid nodes captured by that road network node. In summary, our demand simulation procedure estimates the exact time and road network node location for each van and small ambulance trip. We use the estimates as input to the tactical simulation model described in Section~\ref{Sim}. 

\subsection{Travel time prediction}\label{TT}

We compare four machine learning models and several baseline approaches for predicting travel time on the Dhaka road network according to time of day and day of week, using a dataset of vehicle trips collected by our custom-made GPS devices. We find that a random forest model performs the best, with a $43.3-64.2\%$ improvement in prediction accuracy (measured with root mean squared error) over the baseline approaches. See Section~\ref{TT:Preds} for details. We use the final random forest model trained with all available data to estimate the baseline travel time $\hat{c}_{ij}$ for each edge, which is used as part of the uncertainty set described in Section~\ref{Opt:TTUnc}. To the best of our knowledge, this paper is the first to use real travel time data from a LMIC for optimization.

\subsection{Tactical simulation model}\label{Sim}

The main focus of the simulation model is to capture the effects of congestion (i.e., waiting time) on overall response times. Our approach is similar to the model developed by \cite{McCormack2015} with three key differences (due to the lack of historical data):
\begin{enumerate}
\item We simulate the time and location of emergency trips using the procedure outlined in Section~\ref{tripsim}.
\item We simulate the travel time between the ambulance base and the patient (and between the patient and hospital) by solving the robust shortest path problem with edge lengths predicted according to the hour of the day and day of the week. 
\item We simulate the scene time using an exponential distribution with an average scene time of 15 minutes because there is no data on scene time from Dhaka \citep{Brown2016, Nagata2016}.
\end{enumerate}
\ref{Apd:Sim} provides a detailed description of our simulation framework. The output of the simulation model is the waiting time, drive time, scene time, transport to hospital time, and the return to home base travel time (if applicable) for each emergency trip. We use response time to denote the summation of waiting time and drive time.

\subsection{Experimental setup}\label{SolApproach}

For all experiments in Section~\ref{DhakaExp}, we solve formulation~\eqref{NFFsolve} using the heuristic scenario generation (HSGen) algorithm with $10$ random starts and $10$ random interchanges. The optimization model solutions are then input to the simulation model described in Section~\ref{Sim} to evaluate the tactical system performance over a seven day period. We use a three day warm up period to reach the steady state system.

Unless otherwise indicated, we use the uncertainty sets described in Sections~\ref{Opt:DUnc} and \ref{Opt:TTUnc} with 100 ambulance demand scenarios and a travel time budget ($B$) of 1000 seconds. Through a detailed sensitivity analysis on the travel time budget (see \ref{Apd:BudgetSens}), we find that optimizing outpost locations using a budget of 1000 seconds generates solutions that perform comparably to solutions optimized for other budgets. Sections~\ref{TimeofDay}, \ref{AmbLocExps}, and \ref{RobustValue} use the ambulance road network, while Section~\ref{CNGvalue} uses both the ambulance and complete road networks. All optimization experiments were programmed using \texttt{MATLAB 2016a} and linear programming sub-problems were solved using \texttt{Gurobi 7.0}. All simulations were programming using \texttt{Python 3.5}. The HSGen algorithm was able to solve each large-scale problem instance in under one hour, and most were solved within 10 minutes. These real-world instances are comparable with the largest problems solved in the facility location literature and papers that focus on problems this large exclusively use heuristic methods \citep{Fischetti2017}.

\section{Policy experiments}\label{DhakaExp}

In this section, we demonstrate the application of our models using data from Dhaka. Each of the following subsections addresses a policy question relevant to the design of an emergency response system: 1) Should different outposts be used for different times of day? (Section~\ref{TimeofDay}) 2) What performance improvements are possible by optimizing outpost locations? (Section~\ref{AmbLocExps}) 3) How much can the system be improved by using small ambulances? (Section~\ref{CNGvalue}) and 4) How important is it to consider uncertainty when designing an emergency response network? (Section~\ref{RobustValue}). 
~\ref{Apd:NumAper},  \ref{Apd:BudgetSens}, and~\ref{Apd:OptVSim} quantify the impact of the number of ambulances per outpost, the impact of the robust travel time budget, and the differences between the optimization and simulation results, respectively.

\subsection{Should different outposts be used for different times of day?}\label{TimeofDay}

In this section, we quantify the benefit of using different outpost locations for different times of day and days of the week, which we refer to as temporal snapshots. In many developed countries, demand is estimated at a fine spatiotemporal resolution \citep{Zhou2015}, allowing ambulances to be repositioned and response to be optimized for different snapshots \citep{VanB2017,Nasro2018}. However, there is a second key motivation for intra-day changes in ambulance locations in LMICs, which is the impact of changing traffic patterns on travel times. We observed first-hand on several occasions during our field work the dramatic increase in travel times in different parts of the city during the evening rush hour. While traffic is less of a concern in high-income countries, emergency vehicles typically face the same traffic conditions as regular road users in LMICs since other vehicles do not (or cannot) yield to ambulances. Thus, our experiments in this section compare the performance of a system that changes outpost locations according to time of day versus a configuration that keeps the ambulance outposts static at all times.

We use baseline travel times for three different temporal snapshots: weekday rush hour (Monday between 6pm and 7pm), weekday overnight (Monday between 2am and 3am), and weekend midday (Saturday between 12pm and 1pm). We use daytime population scenarios for rush hour and we use nighttime population scenarios for weekday overnight and weekend midday. For all three snapshots, we solve \eqref{NFFsolve} with $P=20$. We simulate the performance of each set of outpost locations on all three temporal snapshots using seven ambulances per outpost, chosen based on our investigation of the impact of the number of ambulances per outpost (see~\ref{Apd:NumAper} for details).

\begin{figure}[t]
\centering
\includegraphics[width=.45\textwidth]{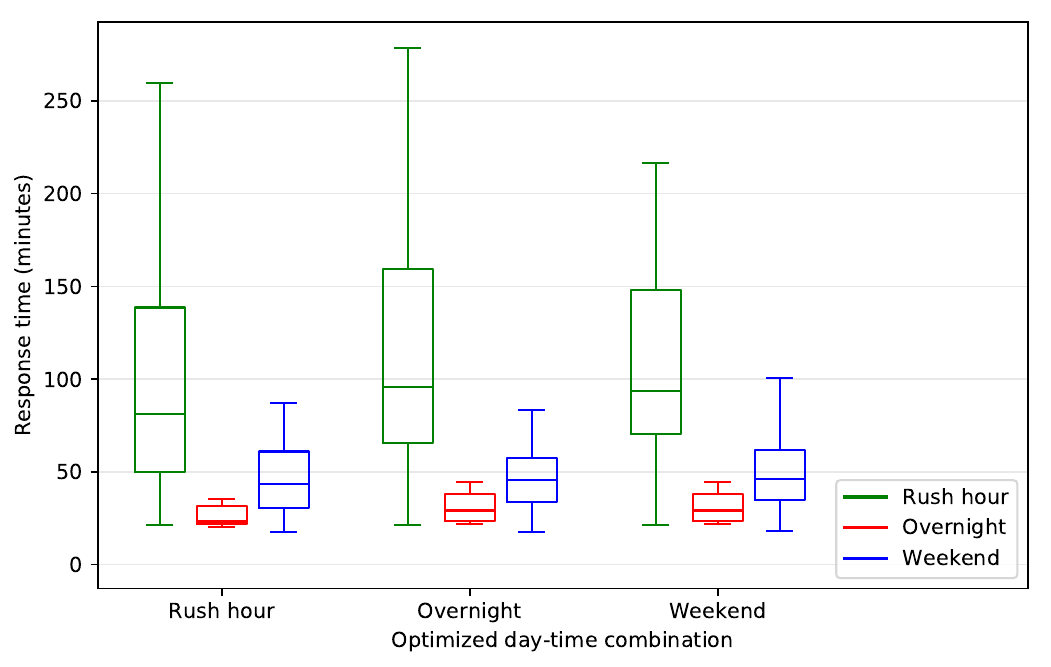}
\caption{The response time performance of outpost locations optimized for one specific snapshot and applied to other snapshots. \label{DayTime}}
\end{figure}

Figure~\ref{DayTime} displays the distribution of response times from the simulation model corresponding to outpost locations optimized for each of the three temporal snapshots. We find that the median ambulance response time is 58.0 and 38.1 minutes longer during rush hour as compared to overnight and weekend, respectively. During rush hour, the rush hour-optimized locations have a median response time that is 14.4 min (15.0\%) and 12.5 min (13.4\%) and faster than the median response time of  the overnight- and weekend-optimized locations, respectively. During overnight and weekend, the rush hour-optimized locations are 5.9 min (20\%) and 0.2 min (0.5\%) better than the best outpost locations, respectively.


\subsubsection{Discussion and policy implications.}

Our results suggest that ambulance providers in Dhaka do not need to optimize outpost locations by time of day or day of week. Instead, providers can use static outpost locations optimized for daytime rush hour. The rush hour-optimized locations produce significant gains in response time during rush hour, while maintaining similar performance to specialized outpost locations at other times of the day. This finding is important because it supports reduced system complexity by removing the need to reposition emergency vehicles. In LMICs, it has been shown that complex solutions are far less likely to succeed compared to simple ones \citep{Bradley2017}. Thus, our rush hour-optimized solution is recommended since it is optimal for the busiest time of day, close to optimal otherwise, and more likely to be implemented than a solution that involves regular repositioning. 

\subsection{What performance improvements are possible by optimizing outpost locations?}\label{AmbLocExps}

Given the results in Section~\ref{TimeofDay}, we turn our attention to designing a static ambulance emergency response network for daytime rush hour and quantifying the gain from shifting away from the current practice of having hospital-based ambulances. There are currently $87$ hospitals with emergency departments. Many of these hospitals have their own ambulance services, while others rely on private services. In both cases, ambulance providers typically position their fleets at the hospitals. We estimate the total number of ambulances in Dhaka by assuming that each of the 19 government hospitals has a fleet of seven ambulances, while each of the 68 private hospitals has a fleet of two ambulances. We obtain these estimates based on the volume differences between government and private hospitals, and based on our field experience. In total, we estimate that Dhaka has approximately 269 ambulances. To calculate response times, we assign each hospital and its fleet of ambulances to the closest node on the ambulance road network, resulting in $67$ unique locations. Using these locations, Section~\ref{PolicyQ1} determines the baseline performance of the current hospital-based outpost locations in Dhaka.

We then consider three policy experiments for improving baseline ambulance response times that may inform the decision making of existing ambulance providers interested in improving or expanding their operations, as well as possibly new entrants or the government looking to design a system from scratch. In particular, Section~\ref{RepoOutpost} quantifies the value of repositioning current outposts. Section~\ref{NewOutpost} quantifies the value of adding additional outpost locations to the current network. Finally, Section~\ref{GreenF} quantifies the performance of an ambulance network that is designed from scratch, without consideration of current outpost locations. We measure performance of the ambulance networks over an entire week using our simulation model. 

\subsubsection{What is the baseline performance of current hospital-based outpost locations?}\label{PolicyQ1}

The median response time of the current outpost locations is $47.2$, $24.1$, and $33.3$ minutes, during rush hour, overnight, and weekend, respectively. The variability in average response time is much larger during rush hour, with a $168.6$ minute difference between the best and worst response times, compared to $91.0$ and $94.5$ minute differences between the best and worst response times during overnight and weekend, respectively.

\begin{figure}[t]
\centering
\subfigure[\ Repositioned outposts \label{BL2}]{
\includegraphics[width=.3\textwidth]{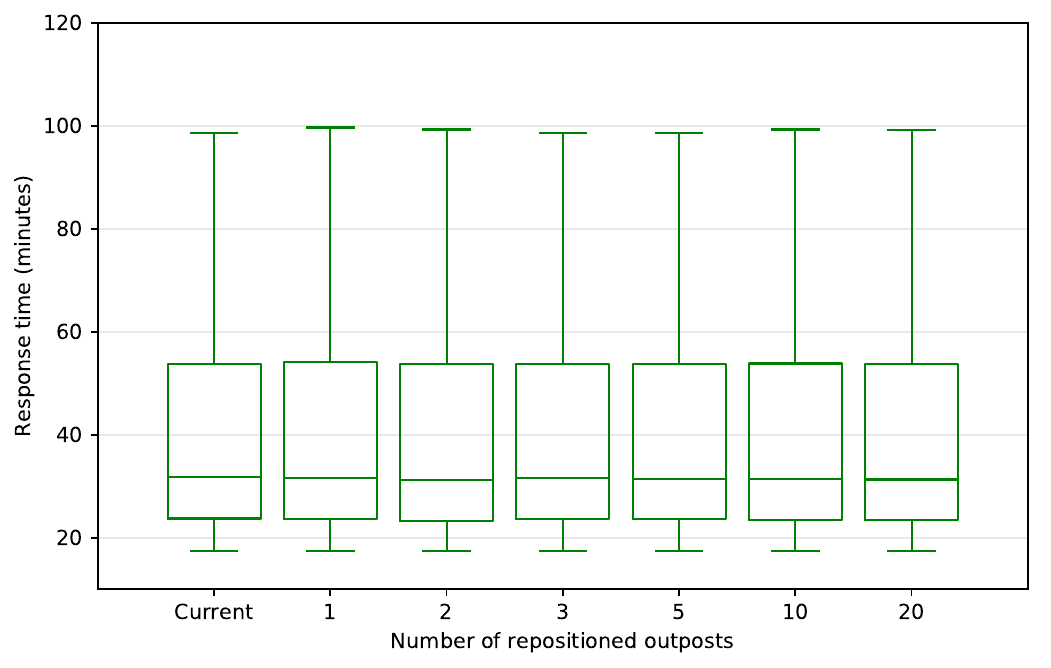}}
\subfigure[\ Additional outposts \label{BL3}]{
\includegraphics[width=.3\textwidth]{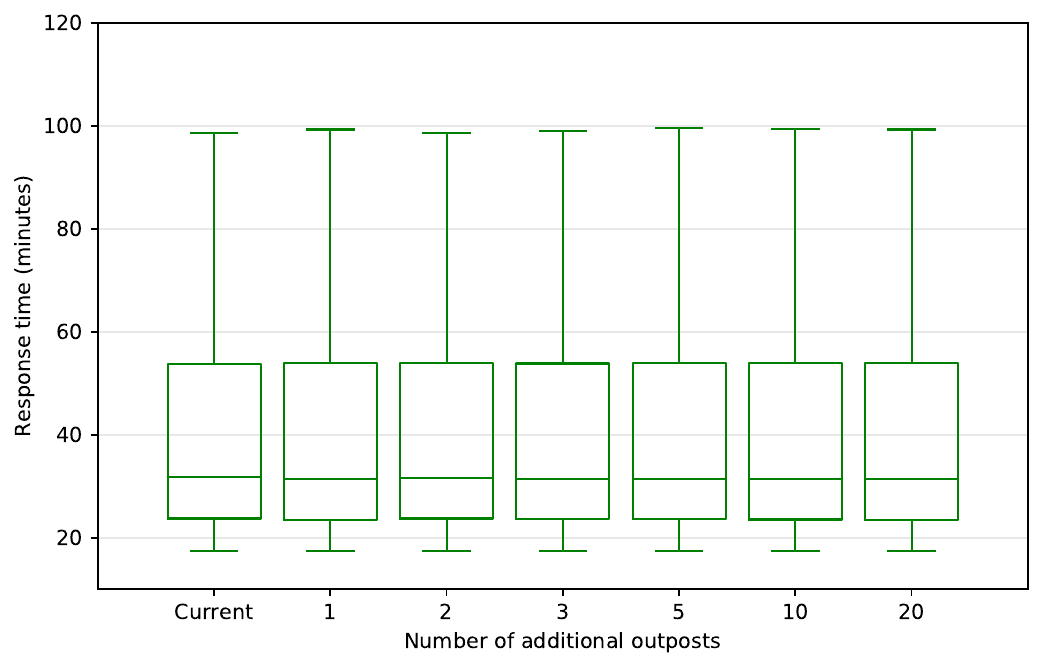}}
\subfigure[\ New outposts \label{BL4}]{
\includegraphics[width=.3\textwidth]{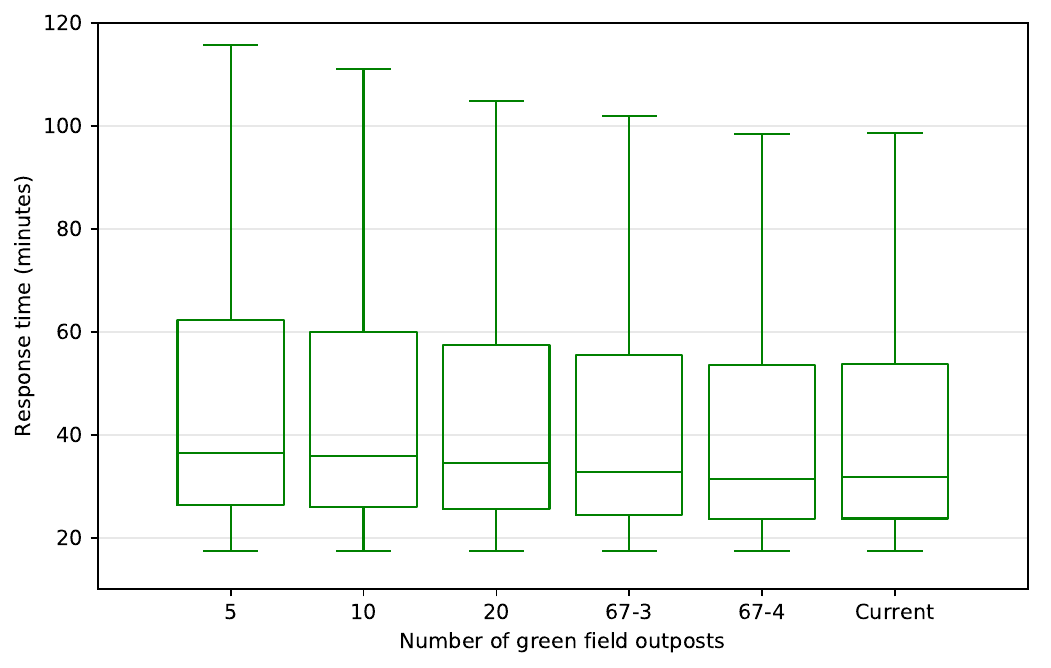}}
\caption{Response time performance for different emergency response network configurations. The  67-3 and 67-4 labels in (c) correspond to 67 outposts with 3 and 4 ambulances per outpost, respectively. \label{Baseline}}
\end{figure}

\begin{figure}[t]
\centering
\subfigure[\ Single repositioned outpost during rush hour \label{BLmaps1}]{
\includegraphics[width=.23\textwidth]{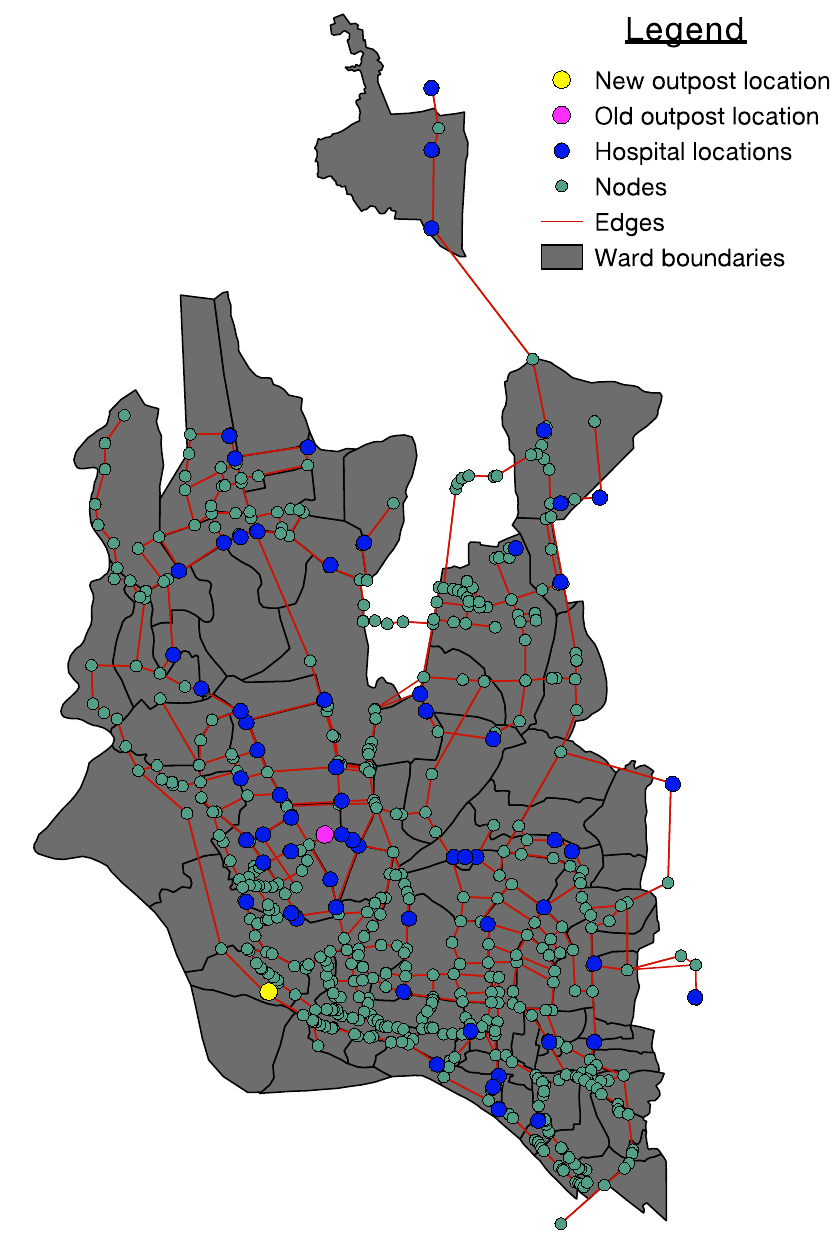}}
\hspace*{3em}
\subfigure[\ Single additional rush hour, overnight, and weekend outpost \label{BLmaps2}]{
\includegraphics[width=.25\textwidth]{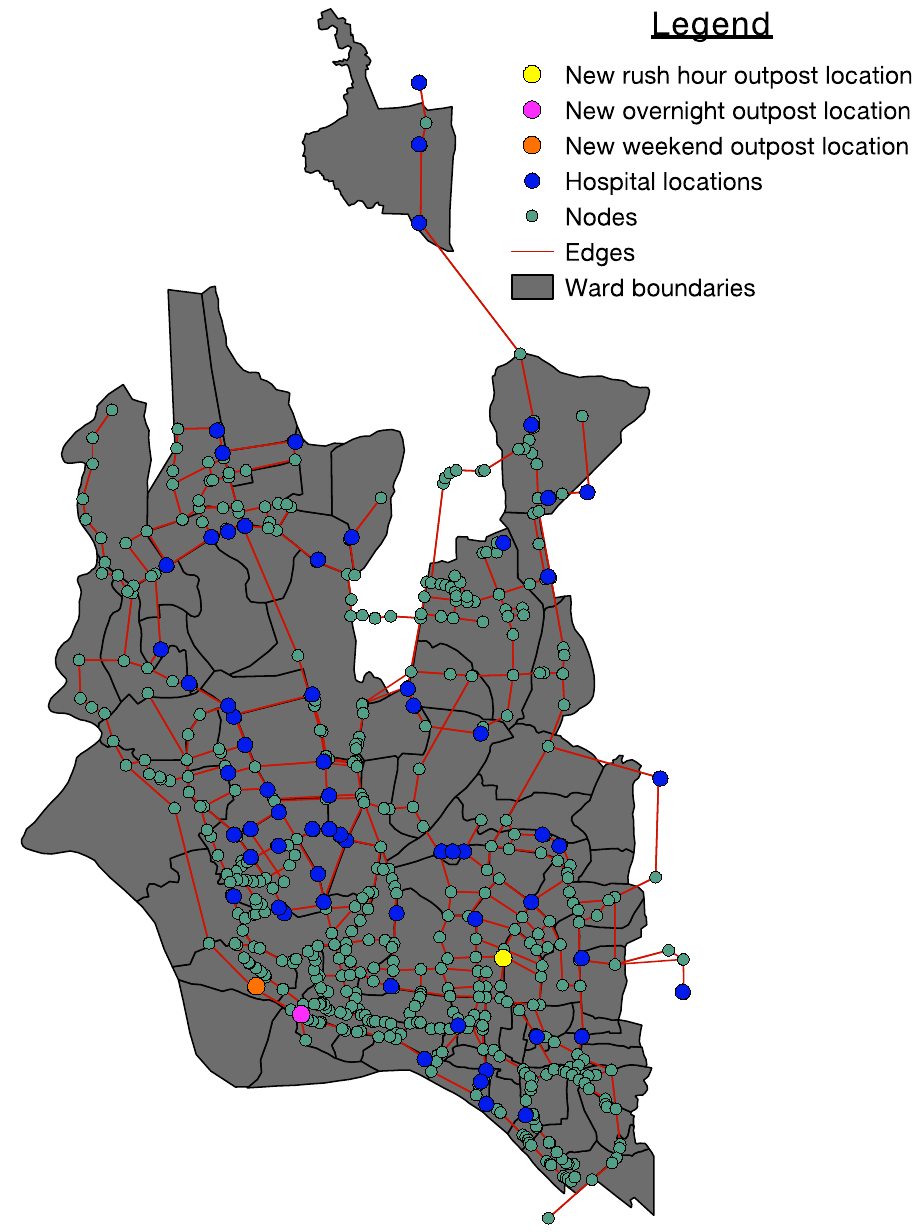}}
\hspace*{3em}
\subfigure[\ 20 new outposts during rush hour \label{BLmaps4}]{
\includegraphics[width=.3\textwidth]{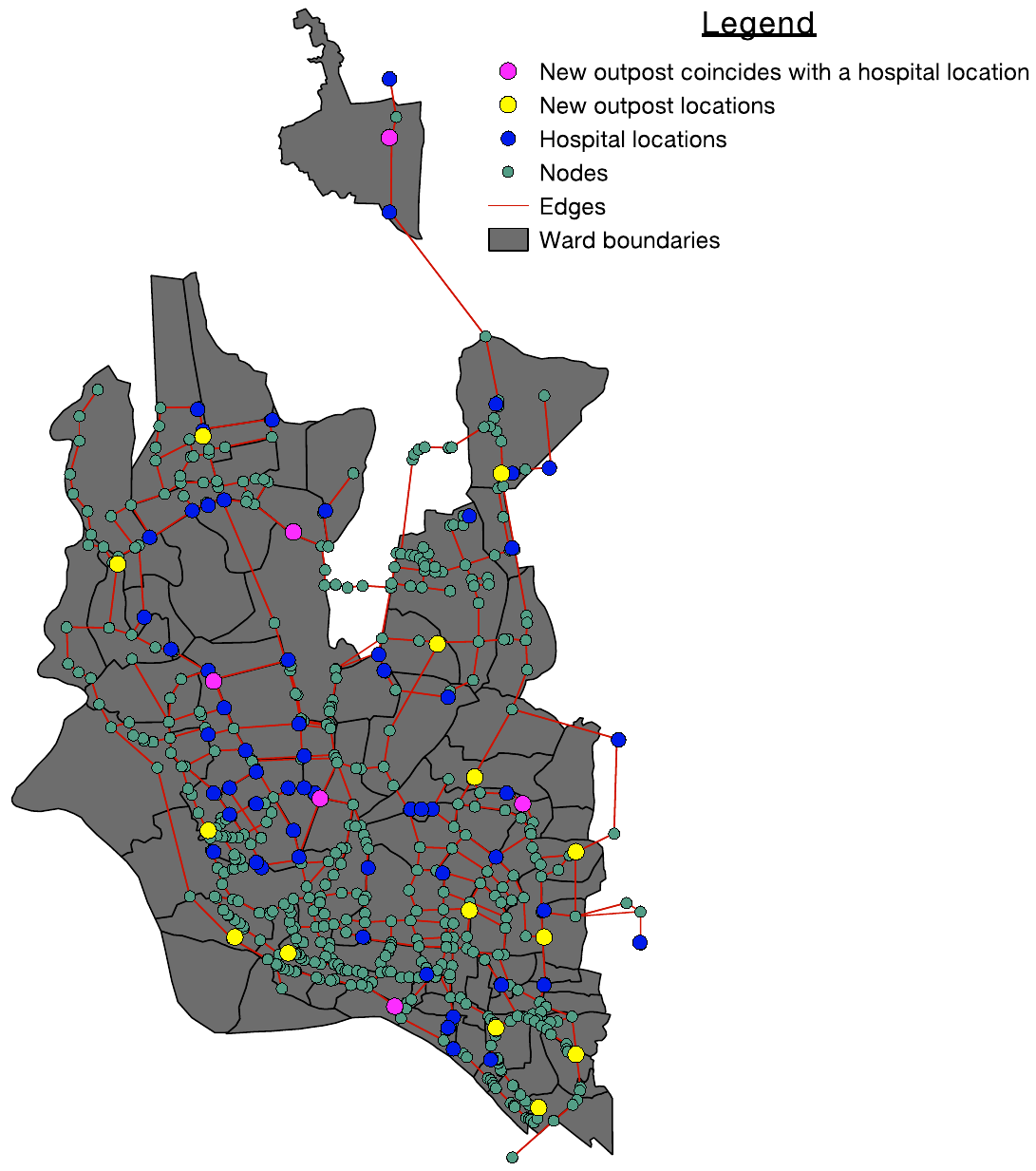}}
\caption{A visualization of the outpost locations for different improvment policies.\label{BaselineLocations}}
\end{figure}

\subsubsection{What is the value of repositioning current outpost locations?}\label{RepoOutpost}

We use a modified version of HSGen for these experiments. For each random start, we randomly choose the required number of outposts to reposition from the current locations and fix all the other outposts for the remainder of the algorithm. As a result, the problem is reduced to determining the location of a specified number of outposts given a set of incumbent outpost locations. When an outpost is repositioned, all ambulances at that outpost are also repositioned.

Figure~\ref{BL2} displays the distribution of response times for each number of repositioned outposts. Repositioning outposts provides only marginal improvements in response time. For example, re-locating one outpost provides no improvement in response time, while repositioning $20$ outposts provides only a 0.5 min (1.2\%) response time improvement.

Figure~\ref{BLmaps1} shows a representative solution from the one-outpost repositioning problem.  The current outpost locations are blue, the old outpost location is pink, and the repositioned outpost location is yellow. Although the Euclidean distance between the old and new outpost locations is only 3.3 km, the time to travel between them is 403.8 minutes during rush hour, meaning that the new location can provide quicker service to an area that would otherwise see significant delays during rush hour. 

\subsubsection{What is the value of adding additional outpost locations to the current network?}\label{NewOutpost}

Figure~\ref{BL3} displays the distribution of response times as a function of outposts added. Note that additional outposts are selected from a candidate set that includes all nodes without a facility and additional outposts are staffed with two ambulances. The addition of new outposts provides nearly the same value as repositioning outposts, suggesting that some current outposts provide minimal value. Figure~\ref{BLmaps2} displays the location of a single additional rush hour, overnight, and weekend outpost. The additional weekend outpost is the same as the repositioned outpost shown in Figure~\ref{BLmaps1}. Although the additional rush hour location is quite different, it is located in an area with many business, government offices, and universities; during rush hour, this area is particularly busy with people commuting home. 

\subsubsection{What is the value of designing a new emergency response network?}\label{GreenF}

Figure~\ref{BL4} displays the response time distribution for newly optimized networks. Ambulances are distributed uniformly over the the new network outposts so that the entire system has a total of 140 ambulances (129 fewer than the current system), unless otherwise stated. We observe steady response time improvements for all new greenfield solutions. The response time performance of 20 new outposts is only 2.9 minutes (6.3\%) worse than the current 67 outposts, suggesting that similar performance can be achieved with only one-third of the current locations and roughly half as many ambulances. Figure~\ref{BLmaps4} displays the location of 20 new outposts in relation to the current outposts. The new outposts are more strategically spread out compared to the current hospital locations that are concentrated in central Dhaka. For example, new outposts are added to the southwest and east of the city, which include low-income areas there were previously under-served by hospital-based outposts.

\subsubsection{Discussion and policy implications.}

Our first two experiments (Sections~\ref{RepoOutpost} and \ref{NewOutpost}) measure gains from local changes to the current network. The results suggest that policies focused on repositioning current outposts or adding additional outposts provide little value. Furthermore, the improvements from repositioning current outposts are nearly identical to adding new outposts. Practically, this result suggests that some of the current outpost locations are contributing very little to the overall response time calculation (i.e., they are rarely the fastest responding outpost to any given demand point).

If we consider a move towards centralization and a complete redesign of the current system, our third experiment shows that we can achieve roughly the performance of the current system with one-third of the outpost locations and half as many ambulances. The non-governmental organization behind the newly implemented 999-number or a formal government agency seeking to implement a centralized emergency response system may consider a complete re-design. Examining the 20 optimized locations from this experiment we find that nine of them coincide with hospital locations, while the other 11 are located off-site. Another way to view these results is that over 40 of the current hospital-based ambulance outposts can be removed without much impact on city-wide response times, or put to better use by concentrating the ambulances at fewer, more strategically located outposts.

Our experiments recommend putting outposts in the southwest, southeast, and northeast wards, suggesting that these areas are generally under-served. The southwest seems particularly under-served since both repositioned and newly added outposts are located there. Knowing the demographics of the city, this result is not particularly surprising: the southwest wards form part of Old Dhaka and encompasses very dense low-income areas (see Figure~\ref{Maps}) that have poor access to emergency transportation.

Overall, the key takeaway is that the current ambulance network in Dhaka is a dominated solution: response times in Dhaka can be significantly reduced without adding new resources, or equivalently, many fewer resources can be employed to match the current level of performance. Our modeling framework can play a pivotal role in the process to help decision makers strategically position their current ambulance resources. Of course, complementary initiatives will be required to achieve these gains, such as better public education about emergency medical transport and awareness of the newly created 999 number, which became operational in December 2017.

\subsection{How much can the system be improved by using small ambulances?}\label{CNGvalue}

In this section, we consider the hypothetical situation where the city is served by a fleet of small ambulances that are able to traverse every road in the complete road network. Compared to the ambulance road network, the complete road network provides access to a larger portion of the city, including many dense low-income areas that are not accessible to van ambulances. In some areas of the city, an entire sub-network of the complete road network is reduced to a single node in the ambulance road network. However, the distance between nodes in the sub-network and the nearest ambulance network node may be quite far, and it may be unrealistic to assume patients will coordinate multiple modes of transportation for different legs of their trip. As a result, we hypothesize that much of the emergency demand that arises from these low-income areas is lost or unserved. In Section~\ref{CNGdemand}, we quantify the potential emergency demand lost as a result of lack of access via the ambulance road network. To do this, we generate 100 demand scenarios using the prediction models for both van and small ambulance demand from Section~\ref{ED:DemandEst} and map this demand to nodes on the complete road network. Nodes that belong to the ambulance network retain the sum of the van and small ambulance demand, while demand corresponding to complete network nodes that are not present in the ambulance network are assumed to be lost. 

In addition to potentially capturing more demand, the complete road network also provides more routing options for small ambulances, which in turn may enable them to better avoid congestion and deal with travel time uncertainty. In Section~\ref{CNGroutes}, we quantify the value of increased routing flexibility provided by the complete network. We start by mapping demand (van plus small ambulance demand) to nodes in the ambulance network. Then, we use our simulation model to evaluate the response time performance of the current 67 hospital-based outpost locations as well as 20 new locations on both the ambulance and complete road networks. The 20 new locations are optimized for the corresponding road network, so they represent two distinct solutions.

\subsubsection{How much potential demand is lost by van ambulances restricted to the ambulance road network?} \label{CNGdemand}

The complete network captures an average of 769,790 small ambulance trips per year, while the ambulance road network only captures 225,559 trips per year, representing a potential loss of 544,231 ambulance (70.7\%) trips. These numbers represent an upper bound on the true number of ambulance trips because we have implicitly assumed that all available ambulance trips will be captured (in reality, some may be missed). Figure~\ref{DL2} visualizes the lost ambulance demand. The 530 green nodes are those that capture demand in both the ambulance and complete networks, while the 4,828 blue nodes only capture demand in the complete network and therefore, represent lost demand for the ambulance road network. 

\subsubsection{What is the value of the increased routing flexibility offered by small ambulances?}\label{CNGroutes}

\begin{figure}[t]
\centering
\subfigure[\ Lost demand \label{DL2}]{
\includegraphics[width=.3\textwidth]{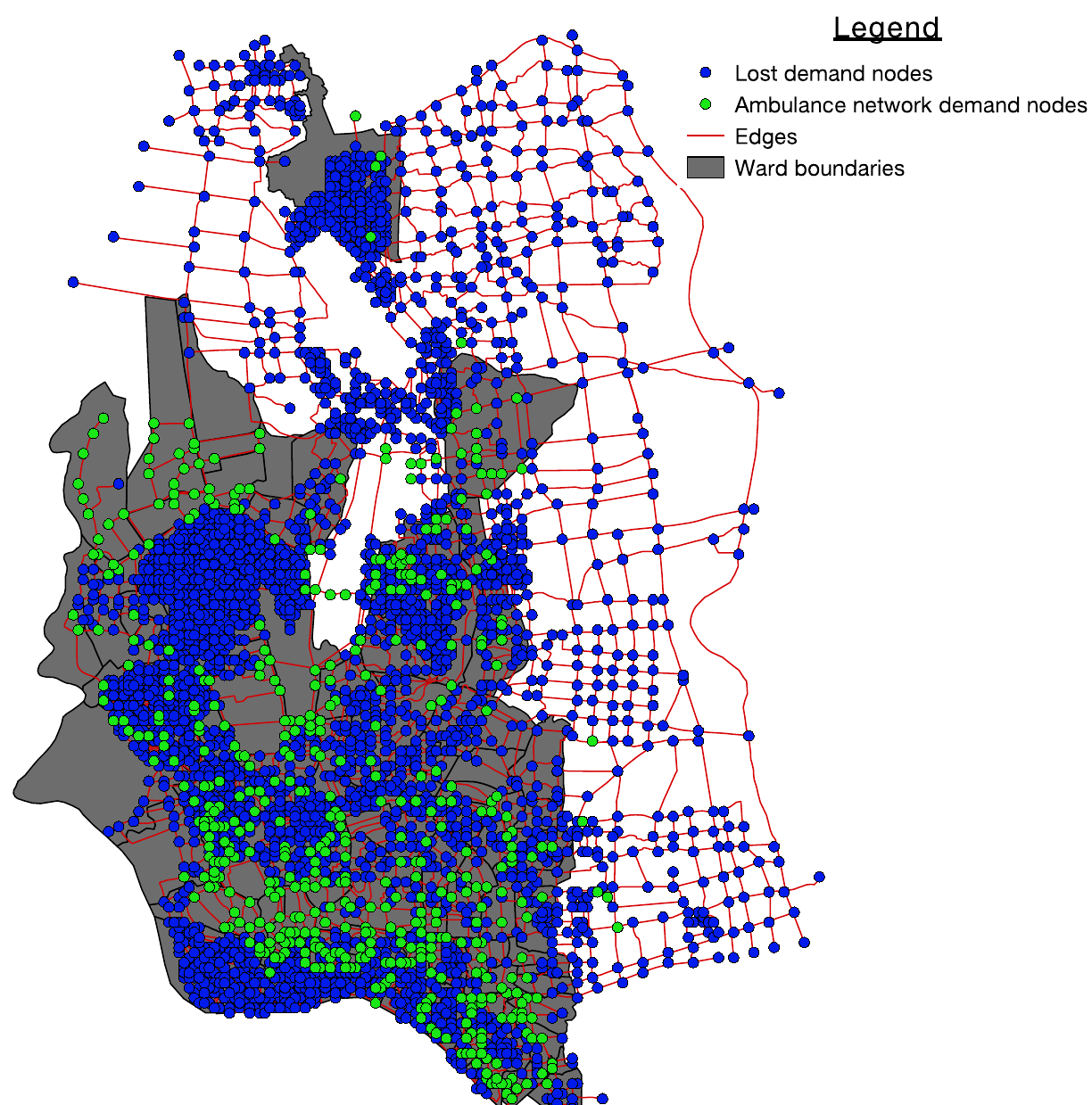}}
\subfigure[\ Response time performance for baseline and new locations\label{AC1}]{
\includegraphics[width=.33\textwidth]{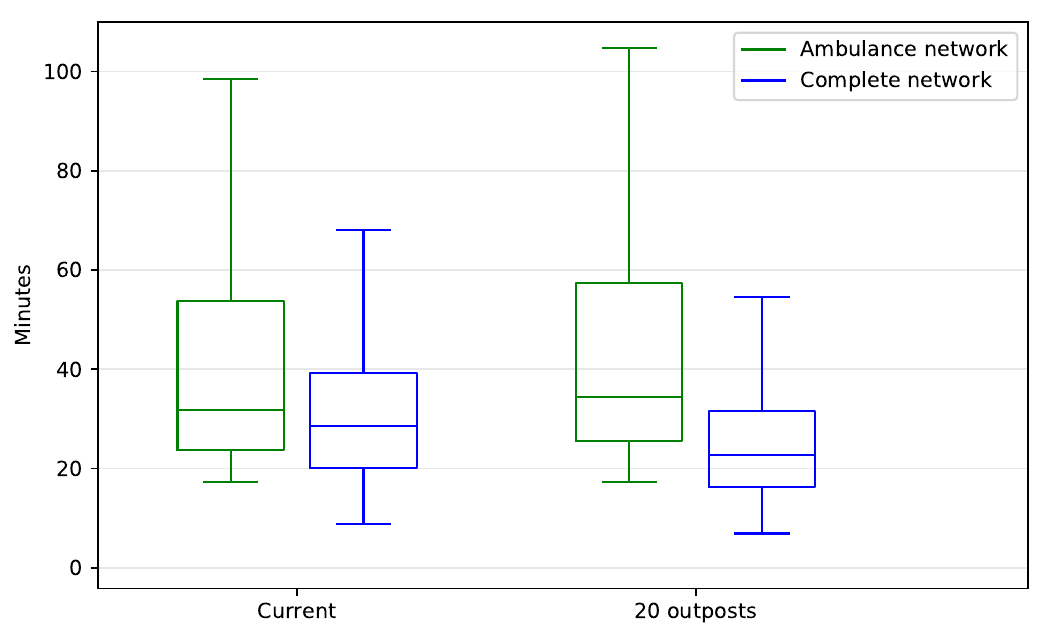}}
\hspace*{1em}
\subfigure[\ Twenty new outpost locations optimized for the two road networks\label{AC2}]{
\includegraphics[width=.3\textwidth]{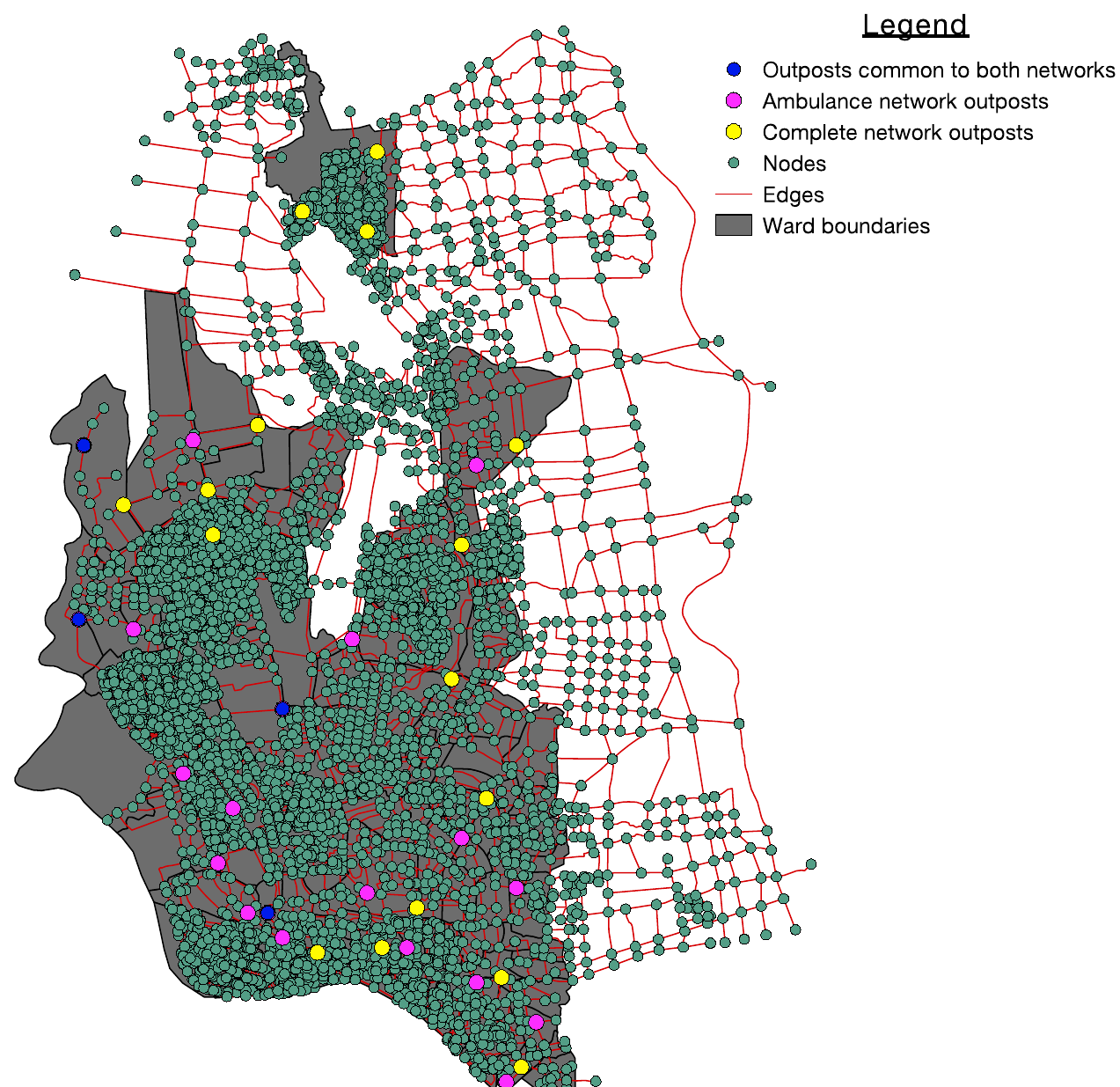}}
\caption{Comparing the performance of the ambulance and complete road networks \label{ANCN}}
\end{figure}

Figure~\ref{AC1} displays the response time performance for the current 67 baseline outpost locations and 20 new outpost locations on both the ambulance and complete road networks. The median response time of the current 67 locations on the ambulance road network is 31.8 minutes over an entire week. Small ambulances located at the same outposts are able to reduce the median response time to 28.5 minutes (a 10.1\% reduction). If we consider 20 new outpost locations, we get a larger reduction in median response time of 17.8\%, from 34.5 minutes to 28.4 minutes. During the busiest time of week, rush hour, the improvements are even larger at 23.7\% and 35.2\% for current and new outpost locations, respectively. Figure~\ref{AC2} shows the 20 new outpost locations for both the ambulance (pink nodes) and complete (yellow nodes) road networks. The four blue nodes represent outpost locations that are the same in both networks. 


\subsubsection{Discussion and policy implications.}

The results in this section represent the first attempt to provide quantitative evidence of the potential benefit of small ambulances in an LMIC.    

Note that that 23\% of survey respondents indicated that traditional ambulance vans were either not available or too slow to reach their location (see Section~\ref{ED:Survey}). Small ambulances offer a potential solution to both these issues. Our results have three policy implications:
\begin{enumerate}
\item Smaller response vehicles can potentially capture three times more emergency demand than traditional van ambulances in Dhaka. Much of the additional demand captured is generated from nodes in hard-to-reach and low-income areas, such as urban slums (southern and western clusters of nodes in Figure~\ref{DL2}). These areas are known to already suffer from poor access and availability of emergency medical care.
\item Smaller response vehicles are able to reduce the median average response time by roughly 10-18\% over the entire week and 24-35\% during rush hour. These reductions are entirely due to increased routing flexibility offered by having nimbler vehicles navigating a larger road network. These results may even be somewhat conservative because we did not incorporate the fact that small ambulance are typically able to travel faster than larger ambulances.
\item Our results demonstrate that the outpost locations chosen for small ambulances are very different from those chosen for traditional van ambulances. This result emphasizes the importance of considering small ambulances independently; we cannot assume they should be positioned alongside traditional ambulances, even if the ambulance outpost locations are themselves optimized, because they are optimized for a different road network.
\end{enumerate}

Overall, the key takeaway from these experiments is that small ambulances have the potential to not only significantly improve system efficiency through lower response times, but also simultaneously improve equity and access by capturing substantial demand in the hardest to reach areas of the city. Although both van and small ambulances have similarly limited medical capabilities and are not typically staffed by trained paramedics, further research is needed to evaluate their medical and operational impact in LMICs.

\subsection{How important is it to consider uncertainty when designing an emergency response network?}\label{RobustValue}

In this section, we quantify the value of robustness by comparing our robust optimization model to the deterministic model (\textbf{NFF}), as well as to a perfect information formulation that solves \textbf{NFF} after the uncertainty has been realized. We focus on the situation where 20 new outposts are being located. We also examine how the performance gaps vary as the travel time budget is varied. The deterministic formulation uses the average demand and baseline travel times with no uncertainty, while the perfect information formulation finds a unique solution for each demand scenario. In this section, we directly report the optimization results, rather than evaluating them via simulation.

\begin{figure}[t]
\centering
\subfigure[\ Median \label{Rob1}]{
\includegraphics[width=.4\textwidth]{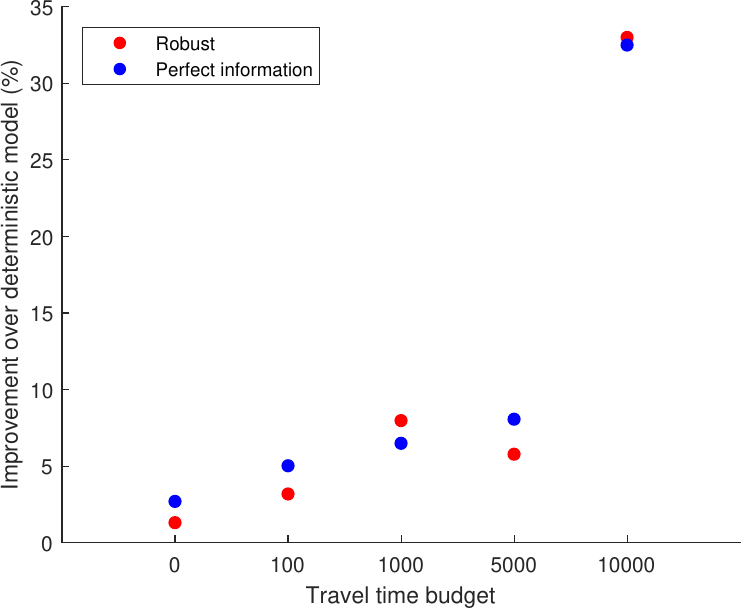}}
\subfigure[\ Worst-case \label{Rob1}]{
\includegraphics[width=.4\textwidth]{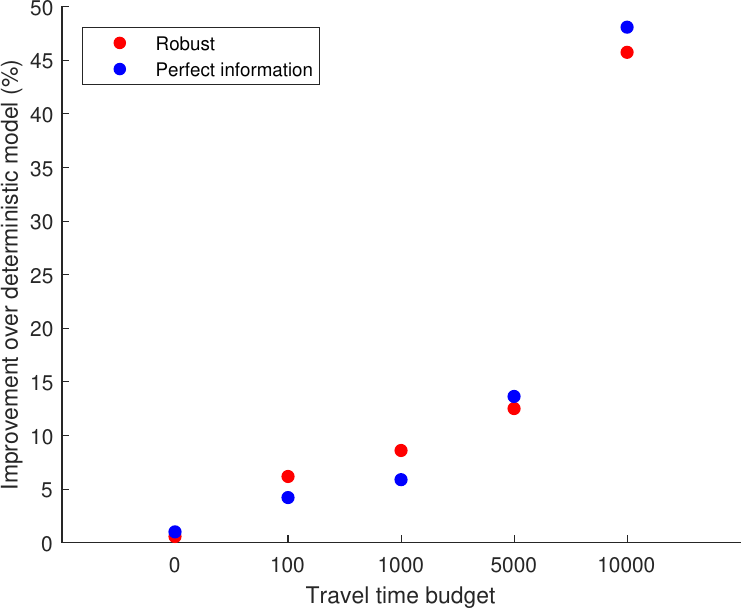}}
\caption{Response time improvement of the robust and perfect information formulations over the deterministic formulation as a function of travel time uncertainty.\label{Robust}}
\end{figure}

Figure~\ref{Robust} displays the response time improvement of the robust and perfect information solutions over the deterministic solution for different levels of travel time uncertainty. Because both models are solved using a heuristic, there are instances where the robust solution slightly outperforms the perfect information solution. For a travel time budget of 1000 seconds, the robust solution generates a 8.0\% and 8.6\% improvement over the deterministic solution in the median and worst-case average response time, respectively. Compared to the perfect information solution, the robust solution has a median average response time that is only 0.8\% worse, with a worst-case average response time that is 1.9\% better. As expected, the gains from the robust model increase as the size of the uncertainty set grows. For example, with a budget of 10,000 seconds, the robust solution improves upon the deterministic solution by 33.0\% in median and 45.8\% in worst-case average response time. At the same time, the performance of the robust solution continues to track the performance of the perfect information solution quite closely.

\subsubsection{Discussion and policy implications.}
Our results demonstrate that a robust optimization framework tailored for the uncertainties faced by LMICs is able to produce solutions that significantly outperform solutions that do not consider uncertainty. As expected, the performance gains increase with the amount of uncertainty considered. Furthermore, our robust solutions are comparable to those derived from a perfect information model.  Overall, these results further reinforce the importance of robustness for designing emergency response solutions in environments with substantial uncertainty characteristic of LMICs.

\section{Conclusion}\label{Conc}

In this paper, we developed a comprehensive framework for emergency response optimization that combines two machine learning approaches and a simulation model with a robust optimization model tailored to address the specific challenges faced by LMICs. Our optimization model generalizes previous emergency response models in both high, middle, and low-income countries and provides a unified framework for emergency response optimization under travel time and demand uncertainty. We use two unique datasets that we collected in Dhaka, Bangladesh to train our machine learning models and build our uncertainty sets.

Using our real data and modelling framework, we address four policy questions related to the design of an emergency response system in LMICs, using Dhaka, Bangladesh as a target site. First, we demonstrated that daily population migration has a minimal impact on response times and that outpost locations optimized specifically for rush hour perform well throughout the day and week. Second, we demonstrated that a centralized network designed from a clean slate can replicate the performance of the current system using roughly half of the ambulance resources and one-third of the outpost locations currently in use. Half of the new outposts would coincide with current outpost locations, while the other half should be strategically positioned in the lower-income parts of the city. 
Third, we show that small ambulances may be able to capture three times more demand than van ambulances due to their ability to access parts of the city with narrow roads such as slums. In addition, the routing flexibility offered by the larger road network available to small ambulances can reduce the median average response time by roughly 10-18\% over the entire week and 24-35\% during rush hour, based on our experiments. Our final experiment demonstrated that our robust optimization framework is able to produce networks with average response times that are up to 33\% faster than a deterministic solution, comparable to a network designed with perfect information on the uncertainty. 


\ACKNOWLEDGMENT{The authors gratefully acknowledge Dr. Moinul Hossain for his input on early stages of the project and for leading the data collection efforts. The authors thank Mehedi Hasan for his help with demand data collection, Mahfuzur Rahman Siddiquee for his help with travel time data collection, and all those who volunteered in Dhaka. We are grateful to Prof. Yu-Ling Cheng and Dr. Laurie Morrison for support and advice throughout this project. This research was supported by Grand Challenges Canada and the Natural Sciences and Engineering Research Council of Canada (NSERC).}

\bibliographystyle{informs2014}
\bibliography{AERObib}

\begin{thebibliography}{107}
\providecommand{\natexlab}[1]{#1}
\providecommand{\url}[1]{\texttt{#1}}
\providecommand{\urlprefix}{URL }

\bibitem[{Abdul~Ghani \protect\BIBand{} Ahmad(2017)}]{Ghani2017}
Abdul~Ghani N, Ahmad N (2017) Analysis of mclp, q-malp, and mq-malp with travel
  time uncertainty using monte carlo simulation. \emph{Journal of Computational
  Engineering} .

\bibitem[{Adenso-Diaz \protect\BIBand{} Rodriguez(1997)}]{Adenso1997}
Adenso-Diaz B, Rodriguez F (1997) A simple search heuristic for the mclp:
  Application to the location of ambulance bases in a rural region.
  \emph{Omega} 25:181--187.

\bibitem[{Ahmadi-Javid et~al.(2017)Ahmadi-Javid, Seyedi, \protect\BIBand{}
  Syam}]{Ahmadi2017}
Ahmadi-Javid A, Seyedi P, Syam SS (2017) A survey of healthcare facility
  location. \emph{Computers and Operations Research} 79:223--263.

\bibitem[{Ahmed et~al.(2015)Ahmed, Rahman~Siddiquee, Karim, Zaman, Monzur,
  \protect\BIBand{} Hossain}]{Ahmed2015}
Ahmed N, Rahman~Siddiquee M, Karim R, Zaman M, Monzur R, Hossain M (2015) Map
  matching on sparse gps data: A perspective of a developing city. \emph{10th
  International Conference of Eastern Asia Society For Transportation Studies.}
  10.

\bibitem[{Alanis et~al.(2013)Alanis, Ingolfsson, \protect\BIBand{}
  Kolfal}]{Alanis2013}
Alanis R, Ingolfsson A, Kolfal B (2013) A markov chain model for an ems system
  with repositioning. \emph{Production and Operations Management}
  22(1):216--231,
  \urlprefix\url{http://dx.doi.org/10.1111/j.1937-5956.2012.01362.x}.

\bibitem[{Anderson et~al.(2012)Anderson, Suter, Mulligan, Bodiwala, Razzak,
  \protect\BIBand{} Mock}]{Anderson}
Anderson PD, Suter RE, Mulligan R, Bodiwala G, Razzak JA, Mock C (2012) World
  health assembly resolution 60.22 and its importance as a health care policy
  tool for improving emergency care access and availability globally.
  \emph{Annals of Emergency Medicine} 60:35--44.

\bibitem[{Atamturk \protect\BIBand{} Zhang(2007)}]{Atamturk2007}
Atamturk A, Zhang M (2007) Two-stage robust network flow and design under
  demand uncertainty. \emph{Operations Research} 55:662--673.

\bibitem[{Averbakh(2003)}]{Averbakh2003}
Averbakh I (2003) Complexity of robust single facility location problems on
  networks with uncertain edge lengths. \emph{Discrete Applied Mathematics}
  127:505--522.

\bibitem[{Bagai et~al.(2013)Bagai, McNally, Al-Khatib, Myers, Kim, Karlsson,
  Torp-Pedersen, Wissenberg, van Diepen, Fosbol, Monk, Abella, Granger,
  \protect\BIBand{} Jollis}]{Circ2013}
Bagai A, McNally BF, Al-Khatib SM, Myers JB, Kim S, Karlsson L, Torp-Pedersen
  C, Wissenberg M, van Diepen S, Fosbol EL, Monk L, Abella BS, Granger CB,
  Jollis JG (2013) Temporal differences in out-of-hospital cardiac arrest
  incidence and survival. \emph{Circulation} 128(24):2595--2602.

\bibitem[{Baker \protect\BIBand{} Fitzpatrick(1986)}]{Baker1986}
Baker JR, Fitzpatrick KE (1986) Determination of an optimal forecast model for
  ambulance demand using goal programming. \emph{The Journal of the Operational
  Research} 37:1047--1059.

\bibitem[{Baron et~al.(2011)Baron, Milner, \protect\BIBand{}
  Naseraldin}]{Baron2011}
Baron O, Milner J, Naseraldin H (2011) Facility location: A robust optimization
  approach. \emph{Production and Operations Management} 20:772--785.

\bibitem[{Basar et~al.(2011)Basar, Catay, \protect\BIBand{}
  Unluyurt}]{Basar2011}
Basar A, Catay B, Unluyurt T (2011) A multi-period double coverage approach for
  locating the emergency medical service stations in istanbul. \emph{Journal of
  the Operational Research Society} 62:627--637.

\bibitem[{Basar et~al.(2012)Basar, Catay, \protect\BIBand{}
  Unluyurt}]{Basar2012}
Basar A, Catay B, Unluyurt T (2012) A taxonomy for emergency service station
  location problem. \emph{Optimization letters} 6:1147--1160.

\bibitem[{Bennett et~al.(1982)Bennett, Eaton, \protect\BIBand{}
  Church}]{Bennett1982}
Bennett VL, Eaton DJ, Church RL (1982) Selecting sites for rural health
  workers. \emph{Social Science and Medicine} 16:63--72.

\bibitem[{Beraldi \protect\BIBand{} Bruni(2009)}]{Beraldi2009}
Beraldi P, Bruni ME (2009) A probabilistic model applied to emergency service
  vehicle location. \emph{European Journal of Operational Research}
  196:323--331.

\bibitem[{Beraldi et~al.(2004)Beraldi, Bruni, \protect\BIBand{}
  Conforti}]{Beraldi2004}
Beraldi P, Bruni ME, Conforti D (2004) Designing robust emergency medical
  service via stochastic programming. \emph{European Journal of Operational
  Research} 158:183--193.

\bibitem[{Berchet(2015)}]{Berchet2015}
Berchet C (2015) Emergency care services: trends, drivers, and interventions to
  manage the demand. \emph{Organization for Economic Co-operation and
  Development} .

\bibitem[{Berman et~al.(2013)Berman, Hajizadeh, \protect\BIBand{}
  Krass}]{Berman2013}
Berman O, Hajizadeh I, Krass D (2013) The maximum covering problem with travel
  time uncertainty. \emph{IIE Transactions} 45:81--96.

\bibitem[{Bertsimas \protect\BIBand{} Copenhaver(2017)}]{Bertsimas2017}
Bertsimas D, Copenhaver MS (2017) Characterization of the equivalence of
  robustification and regularization in linear and matrix regression.
  \emph{European Journal of Operational Research} In Press.

\bibitem[{Bradley et~al.(2017)Bradley, Jung, Tandon-Verma, Khoury, Chan,
  \protect\BIBand{} Cheng}]{Bradley2017}
Bradley BD, Jung T, Tandon-Verma A, Khoury B, Chan TCY, Cheng YL (2017)
  Operations research in global health: a scoping review with a focus on the
  themes of health equity and impact. \emph{Health Research Policy and Systems}
  15:32.

\bibitem[{Brandeau \protect\BIBand{} Larson(1986)}]{Brandeau1986}
Brandeau ML, Larson RC (1986) Extending and applying the hypercube queueing
  model to deploy ambulances in boston. \emph{National Emergency Training
  Center} .

\bibitem[{Brotcorne et~al.(2003)Brotcorne, Laporte, \protect\BIBand{}
  Semet}]{Brotcorne2003}
Brotcorne L, Laporte G, Semet F (2003) Ambulance location and relocation
  models. \emph{European Journal of Operational Research} 147(3):451 -- 463,
  ISSN 0377-2217,
  \urlprefix\url{http://dx.doi.org/https://doi.org/10.1016/S0377-2217(02)00364-8}.

\bibitem[{Brown et~al.(2016)Brown, Rosengart, Forsythe, Reynolds, Gestring,
  Hallinan, Peitzman, Billiar, \protect\BIBand{} Sperry}]{Brown2016}
Brown JB, Rosengart MR, Forsythe RM, Reynolds BR, Gestring ML, Hallinan WM,
  Peitzman AB, Billiar TR, Sperry JL (2016) Not all prehospital time is equal:
  Influence of scene time on mortality. \emph{J Trauma Acute Care Surg}
  81(1):93--100, ISSN 2163-0763 (Electronic); 2163-0755 (Print); 2163-0755
  (Linking), \urlprefix\url{http://dx.doi.org/10.1097/TA.0000000000000999}.

\bibitem[{Budge et~al.(2010)Budge, Ingolfsson, \protect\BIBand{}
  Zerom}]{Budge2010}
Budge S, Ingolfsson A, Zerom D (2010) Empirical analysis of ambulance travel
  times: the case of calgary emergency medical services. \emph{Management
  Science} 56:716--723.

\bibitem[{Carson \protect\BIBand{} Batta(1990)}]{Carson1990}
Carson YM, Batta R (1990) Locating an ambulance on the amherst campus of the
  state university of new york at buffalo. \emph{Interfaces} 20:43--49.

\bibitem[{Chan(2017)}]{Chan2017}
Chan TCY (2017) Rise and shock: Optimal defibrillator placement in a high-rise
  building. \emph{Prehospital Emergency Care} 21(3):309--314,
  \urlprefix\url{http://dx.doi.org/10.1080/10903127.2016.1247202}, pMID:
  27858504.

\bibitem[{Channouf et~al.(2007)Channouf, L'Ecuyer, Ingolfsson,
  \protect\BIBand{} Avramidis}]{Channouf2007}
Channouf N, L'Ecuyer P, Ingolfsson A, Avramidis A (2007) The application of
  forecasting techniques to modeling emergency medical system calls in calgary,
  alberta. \emph{Health Care Management Science} 10:25--45.

\bibitem[{Chanta et~al.(2014)Chanta, Mayorga, \protect\BIBand{}
  McLay}]{Chanta2014}
Chanta S, Mayorga ME, McLay LA (2014) Improving emergency service in rural
  areas: a bi-objective covering location model for ems systems. \emph{Annals
  of Operations Research} 221:133--159.

\bibitem[{Chen \protect\BIBand{} Lin(1998)}]{Chen1998}
Chen B, Lin CS (1998) Minmax regret robust 1-median location on a tree.
  \emph{Networks} 31:93--103.

\bibitem[{Cox \protect\BIBand{} Lewis(1966)}]{Cox1966}
Cox D, Lewis P (1966) \emph{The statistical analysis of series of events} (John
  Wiley and Sons).

\bibitem[{Densham \protect\BIBand{} Rushton(1992)}]{Densham1992}
Densham PJ, Rushton G (1992) A more efficient heuristic for solving large
  p-median problems. \emph{Papers in Regional Science} 71:307--329.

\bibitem[{Eaton et~al.(1986)Eaton, H{\'e}ctor, S{\'a}nchez, \protect\BIBand{}
  Morgan}]{Eaton1986}
Eaton DJ, H{\'e}ctor ML, S{\'a}nchez U, Morgan J (1986) Determining ambulance
  deployment in santo domingo, dominican republic. \emph{Journal of the
  Operational Research Society} 113--126.

\bibitem[{Enayati et~al.(2018)Enayati, Mayorga, Rajagopalan, \protect\BIBand{}
  Saydam}]{Ena2018}
Enayati S, Mayorga ME, Rajagopalan HK, Saydam C (2018) Real-time ambulance
  redeployment approach to improve service coverage with fair and restricted
  workload for ems providers. \emph{Omega} 79:67 -- 80, ISSN 0305-0483,
  \urlprefix\url{http://dx.doi.org/https://doi.org/10.1016/j.omega.2017.08.001}.

\bibitem[{Fischetti et~al.(2017)Fischetti, Ljubic, \protect\BIBand{}
  Sinnl}]{Fischetti2017}
Fischetti M, Ljubic I, Sinnl M (2017) Redesigning benders decomposition for
  large-scale facility location. \emph{Management Science} 63:2146--2162.

\bibitem[{Fujiwara et~al.(1987)Fujiwara, Makjamroen, \protect\BIBand{}
  Gupta}]{Fujiwara1987}
Fujiwara O, Makjamroen T, Gupta KK (1987) Ambulance deployment analysis: a case
  study of bangkok. \emph{European Journal of Operational Research} 31:9--18.

\bibitem[{Gabrel et~al.(2014)Gabrel, Lacroix, Murat, \protect\BIBand{}
  Remli}]{Gabrel2014}
Gabrel V, Lacroix M, Murat C, Remli N (2014) Robust location transportation
  problems under uncertain demands. \emph{Discrete Applied Mathematics}
  164:100--111.

\bibitem[{Goins \protect\BIBand{} Conroy(2015)}]{Goins2015}
Goins S, Conroy MB (2015) New york state all payer emergency room visits.
  Technical report, New York State Department of Health Statistical Brief.

\bibitem[{Hakimi(1964)}]{Hakimi1964}
Hakimi SL (1964) Optimum locations of switching centers and the absolute
  centers and medians of a graph. \emph{Operations research} 12:450--459.

\bibitem[{Hakimi(1965)}]{Hakimi1965}
Hakimi SL (1965) Optimum distribution of switching centers in a communication
  network and some related graph theoretic problems. \emph{Operations Research}
  13:462--475.

\bibitem[{Hofleitner et~al.(2012{\natexlab{a}})Hofleitner, Herring, Abbeel,
  \protect\BIBand{} Bayen}]{Hofleitner2012b}
Hofleitner A, Herring R, Abbeel P, Bayen A (2012{\natexlab{a}}) Learning the
  dynamics of arterial traffic from probe data using a dynamic bayesian
  network. \emph{IEEE Transactions on Intelligent Transportation Systems}
  13:1679--1693.

\bibitem[{Hofleitner et~al.(2012{\natexlab{b}})Hofleitner, Herring,
  \protect\BIBand{} Bayen}]{Hofleitner2012a}
Hofleitner A, Herring R, Bayen A (2012{\natexlab{b}}) Arterial travel time
  forecast with streaming data: A hybrid approach of flow modeling and machine
  learning. \emph{Transportation Research Part B} 46:1097--1122.

\bibitem[{Ingolfsson et~al.(2008)Ingolfsson, Budge, \protect\BIBand{}
  Erkut}]{Ingolfsson2008}
Ingolfsson A, Budge S, Erkut E (2008) Optimal ambulance location with random
  delays and travel times. \emph{Health Care management science} 11:262--274.

\bibitem[{Jain et~al.(2012)Jain, Sharma, \protect\BIBand{}
  Subramanian}]{Jain2012}
Jain V, Sharma A, Subramanian L (2012) Road traffic congestion in the
  developing world. \emph{Proceedings of the 2nd ACM Symposium on Computing for
  Development} .

\bibitem[{Kamenetzky et~al.(1982)Kamenetzky, Shuman, \protect\BIBand{}
  Wolfe}]{Kamenetzky1982}
Kamenetzky RD, Shuman LJ, Wolfe H (1982) Estimating need and demand for
  prehospital care. \emph{Operations Research} 30:1148--1167.

\bibitem[{Karim et~al.(2009)Karim, Hansen, Ahmad, \protect\BIBand{}
  Lahiry}]{Karim2009}
Karim MZ, Hansen EL, Ahmad BU, Lahiry S (2009) A retrospective study of illness
  and admission pattern of emergency patients utilizing a corporate hospital in
  dhaka, bangladesh: 2006-2008. \emph{The ORION} 32:1--7.

\bibitem[{Karp(1972)}]{Karp1972}
Karp RM (1972) \emph{Complexity of computer computations} (Springer).

\bibitem[{Kobusingye et~al.(2005)Kobusingye, Hyder, Bishai, Hicks, Mock,
  \protect\BIBand{} Joshipura}]{Kobusingye2005}
Kobusingye OC, Hyder AA, Bishai D, Hicks ER, Mock C, Joshipura M (2005)
  Emergency medical systems in low- and middle-income countries:
  recommendations for action. \emph{Bulletin of the World Health Organization}
  83:626--631.

\bibitem[{Kok et~al.(2012)Kok, Hans, \protect\BIBand{} Schutten}]{Kok2012}
Kok A, Hans E, Schutten J (2012) Vehicle routing under time-dependent travel
  times: The impact of congestion avoidance. \emph{Computers and Operations
  Research} 39(5):910 -- 918, ISSN 0305-0548,
  \urlprefix\url{http://dx.doi.org/https://doi.org/10.1016/j.cor.2011.05.027}.

\bibitem[{Kolesar et~al.(1975)Kolesar, Walker, \protect\BIBand{}
  Hausner}]{Kolesar1975}
Kolesar P, Walker W, Hausner J (1975) Determining the relation between fire
  engine travel times and travel distances in new york city. \emph{Operations
  Research} 23:614--627.

\bibitem[{Levine et~al.(2007)Levine, Gadiraju, Goel, Johar, King,
  \protect\BIBand{} Arnold}]{Levine2007}
Levine AC, Gadiraju S, Goel A, Johar S, King R, Arnold K (2007) International
  emergency medicine: a review of the literature. \emph{Acad Emerg Med},
  182--1833 (14).

\bibitem[{Li et~al.(2011)Li, Zhao, Zhu, \protect\BIBand{} Wyatt}]{Li2011}
Li X, Zhao Z, Zhu X, Wyatt T (2011) Covering models and optimization techniques
  for emergency response facility location and planning: a review.
  \emph{Mathematical Methods of Operations Research} 74(3):281--310, ISSN
  1432-5217, \urlprefix\url{http://dx.doi.org/10.1007/s00186-011-0363-4}.

\bibitem[{Lungu et~al.(2001)Lungu, Kamfose, Hussein, \protect\BIBand{}
  Ashwood-Smith}]{Lungu2001}
Lungu K, Kamfose V, Hussein J, Ashwood-Smith H (2001) Are bicycle ambulances
  and community transport plans effective in strengthening obstetric referral
  systems in southern malawi? \emph{Malawi Medical Journal} 12:16--18.

\bibitem[{Macintyre \protect\BIBand{} Hotchkiss(1999)}]{Macintyre1999}
Macintyre K, Hotchkiss DR (1999) Referral revisited: community financing
  schemes and emergency transport in rural africa. \emph{Soc Sci Med},
  1473--1487 (49).

\bibitem[{Maranzana(1964)}]{Maranzana1964}
Maranzana FE (1964) On the location of supply points to minimize transport
  costs. \emph{Operational Research Quarterly} 15:261--270.

\bibitem[{Matteson et~al.(2011)Matteson, McLean, Woodard, \protect\BIBand{}
  Henderson}]{Matteson2011}
Matteson DS, McLean MW, Woodard DB, Henderson SG (2011) Forecasting emergency
  medical service call arrival rates. \emph{Annals of Applied Statistics}
  5:1379--1406.

\bibitem[{Maxwell et~al.(2010)Maxwell, Restrepo, Henderson, \protect\BIBand{}
  Topaloglu}]{Maxwell2010}
Maxwell MS, Restrepo M, Henderson SG, Topaloglu H (2010) Approximate dynamic
  programming for ambulance redeployment. \emph{INFORMS Journal on Computing}
  22(2):266--281, \urlprefix\url{http://dx.doi.org/10.1287/ijoc.1090.0345}.

\bibitem[{McCormack \protect\BIBand{} Coates(2015)}]{McCormack2015}
McCormack R, Coates G (2015) A simulation model to enable the optimization of
  ambulance fleet allocation and base station location for increased patient
  survival. \emph{European Journal of Operational Research} 247(1):294--309,
  \urlprefix\url{http://dx.doi.org/https://doi.org/10.1016/j.ejor.2015.05.040}.

\bibitem[{Melo et~al.(2009)Melo, Nickel, \protect\BIBand{}
  Saldanha-Da-Gama}]{Melo2009}
Melo MT, Nickel S, Saldanha-Da-Gama F (2009) Facility location and supply chain
  management: A review. \emph{European journal of operational research}
  196:401--412.

\bibitem[{Mirchandani \protect\BIBand{} Odoni(1979)}]{Mirchandani1979}
Mirchandani PB, Odoni AR (1979) Locations of medians on stochastic networks.
  \emph{Transportation Science} 13:85--97.

\bibitem[{Mirchandani \protect\BIBand{} Oudjit(1980)}]{Mirchandani1980}
Mirchandani PB, Oudjit A (1980) Localizing 2-medians on probabilistic and
  deterministic tree networks. \emph{Networks} 10:329--350.

\bibitem[{Nagata et~al.(2016)Nagata, Abe, Nakata, \protect\BIBand{}
  Tamiya}]{Nagata2016}
Nagata I, Abe T, Nakata Y, Tamiya N (2016) Factors related to prolonged
  on-scene time during ambulance transportation for critical emergency patients
  in a big city in japan: a population-based observational study. \emph{BMJ
  Open} 6(1), ISSN 2044-6055,
  \urlprefix\url{http://dx.doi.org/10.1136/bmjopen-2015-009599}.

\bibitem[{Nasrollahzadeh et~al.(2018)Nasrollahzadeh, Khademi, \protect\BIBand{}
  Mayorga}]{Nasro2018}
Nasrollahzadeh AA, Khademi A, Mayorga ME (2018) Real-time ambulance dispatching
  and relocation. \emph{Manufacturing \& Service Operations Management}
  20(3):467--480, \urlprefix\url{http://dx.doi.org/10.1287/msom.2017.0649}.

\bibitem[{Nations(2010)}]{UN2010}
Nations U (2010) The millennium development goals report. Technical report,
  United Nations.

\bibitem[{Neumann(1928)}]{VonNeumann1928}
Neumann J (1928) Zur theorie der gesellschaftsspiele. \emph{Math. Annalen}
  100:295--320.

\bibitem[{Noyan(2010)}]{Noyan2010}
Noyan N (2010) Alternate risk measures for emergency medical service system
  design. \emph{Annals of Operations Research} 181:559--589.

\bibitem[{of~Statistics(2010)}]{HIES}
of~Statistics BB (2010) Report of the household income and expenditure survey.
  Technical report, Ministry of Planning.

\bibitem[{Organization(2013)}]{WHO2013}
Organization TWH (2013) The world health report 2013: Research for universal
  health coverage. Technical report, The World Health Organization.

\bibitem[{Owen \protect\BIBand{} Daskin(1998)}]{Owen1998}
Owen SH, Daskin MS (1998) Strategic facility location: A review. \emph{European
  journal of operational research} 111:423--447.

\bibitem[{Pasupathy(2011)}]{NHPP}
Pasupathy R (2011) \emph{Generating Nonhomogeneous Poisson Processes} (American
  Cancer Society), ISBN 9780470400531,
  \urlprefix\url{http://dx.doi.org/10.1002/9780470400531.eorms0356}.

\bibitem[{Pojani \protect\BIBand{} Stead(2015)}]{Pojani2015}
Pojani D, Stead D (2015) Sustainable urban transport in the developing world:
  beyond megacities. \emph{Sustainability} 7:7784--7805.

\bibitem[{PoSaw et~al.(1998)PoSaw, Aggarwal, \protect\BIBand{}
  Bernstein}]{PoSaw1998}
PoSaw LL, Aggarwal P, Bernstein SL (1998) Emergency medicine in the new delhi
  area, india. \emph{Ann Emerg Med} 32:609--615.

\bibitem[{Post(1944)}]{Post1944}
Post EL (1944) Recursively enumerable sets of positive integers and their
  decision problems. \emph{Bulletin of the American Mathematical Society}
  50:284--316.

\bibitem[{Raftery(1996)}]{Raftery1996}
Raftery KA (1996) Emergency medicine in southern pakistan. \emph{Emerg Med}
  27:79--93.

\bibitem[{Razzak \protect\BIBand{} Kellerman(2002)}]{Razzak2002}
Razzak JA, Kellerman AL (2002) Emergency medical care in developing countries:
  is it worthwhile? \emph{Bulletin of the World Health Organization}
  80:900--905.

\bibitem[{ReVelle \protect\BIBand{} Swain(1970)}]{ReVelle1970}
ReVelle CS, Swain RW (1970) Central facilities location. \emph{Geographical
  Analysis} 2:30--42.

\bibitem[{Salman \protect\BIBand{} Y{\"u}cel(2015)}]{Salman2015}
Salman FS, Y{\"u}cel E (2015) Emergency facility location under random network
  damage: Insights from the istanbul case. \emph{Computers and Operations
  Research} 62:266--281.

\bibitem[{Savas(1969)}]{Savas1969}
Savas ES (1969) Simulation and cost-effectiveness analysis of new york's
  emergency ambulance service. \emph{Management Science} 15(12):B--608--B--627,
  \urlprefix\url{http://dx.doi.org/10.1287/mnsc.15.12.B608}.

\bibitem[{Saydam et~al.(2013)Saydam, Rajagopalan, Sharer, \protect\BIBand{}
  Lawrimore-Belanger}]{Saydam2013}
Saydam C, Rajagopalan HK, Sharer E, Lawrimore-Belanger K (2013) The dynamic
  redeployment coverage location model. \emph{Health Systems} 2(2):103--119,
  \urlprefix\url{http://dx.doi.org/10.1057/hs.2012.27}.

\bibitem[{Schmid et~al.(2001)Schmid, Kanenda, Ahluwalia, \protect\BIBand{}
  Kouletio}]{Schmid2001}
Schmid T, Kanenda O, Ahluwalia I, Kouletio M (2001) Transportation for maternal
  emergencies in tanzania: empowering communities through participatory problem
  solving. \emph{American Journal of Public Health} 91:1589--1590.

\bibitem[{Schmid \protect\BIBand{} Doerner(2010)}]{Schmid2010}
Schmid V, Doerner KF (2010) Ambulance location and relocation problems with
  time-dependent travel times. \emph{European Journal of Operational Research}
  207(3):1293 -- 1303, ISSN 0377-2217,
  \urlprefix\url{http://dx.doi.org/https://doi.org/10.1016/j.ejor.2010.06.033}.

\bibitem[{Schuman et~al.(1977)Schuman, Wolfe, \protect\BIBand{}
  Sepulveda}]{Schuman1977}
Schuman LJ, Wolfe H, Sepulveda J (1977) Estimating demand for emergency
  transportation. \emph{Med Care} 15:738--749.

\bibitem[{Serra \protect\BIBand{} Marianov(1998)}]{Serra1998}
Serra D, Marianov V (1998) The p-median problem in a changing network: the case
  of barcelona. \emph{Location Science} 6:383--394.

\bibitem[{Setzler et~al.(2009)Setzler, Saydam, \protect\BIBand{}
  Park}]{Setzler2009}
Setzler H, Saydam C, Park S (2009) Ems call volume predictions: A comparative
  study. \emph{Computers and Operations Research} 36:1843--1851.

\bibitem[{Shen et~al.(2003)Shen, Coullard, \protect\BIBand{} Daskin}]{Shen2003}
Shen ZJM, Coullard C, Daskin MS (2003) A joint location-inventory model.
  \emph{Transportation science} 37:40--55.

\bibitem[{{Sixtieth World Health Assembly}(2007)}]{WHA6022}
{Sixtieth World Health Assembly} (2007) Agenda item 12.14: Health systems:
  emergency-care systems. Technical report, World Health Organization.

\bibitem[{Snyder(2006)}]{Snyder2006}
Snyder LV (2006) Facility location under uncertainty: a review. \emph{IIE
  Transactions} 38:537--554.

\bibitem[{Sodemann et~al.(1997)Sodemann, Jakobsen, Molbak, Alvarenga,
  \protect\BIBand{} Aaby}]{Sodemann1997}
Sodemann M, Jakobsen MS, Molbak K, Alvarenga IC, Aaby P (1997) High mortality
  despite good care-seeking behaviour: a community study of childhood deaths in
  guinea-bissau. \emph{Bulletin of the World Health Organization} 75:205--212.

\bibitem[{Streatfield \protect\BIBand{} Karar(2008)}]{Streatfield2008}
Streatfield PK, Karar ZA (2008) Population challenges for bangladesh in the
  coming decades. \emph{J Health Popul Nutr} 26:261--272.

\bibitem[{Sudtachat et~al.(2016)Sudtachat, Mayorga, \protect\BIBand{}
  Mclay}]{SUD2016}
Sudtachat K, Mayorga ME, Mclay LA (2016) A nested-compliance table policy for
  emergency medical service systems under relocation. \emph{Omega} 58:154 --
  168, ISSN 0305-0483,
  \urlprefix\url{http://dx.doi.org/https://doi.org/10.1016/j.omega.2015.06.001}.

\bibitem[{Teitz \protect\BIBand{} Bart(1968)}]{Teitz1968}
Teitz MB, Bart P (1968) Heuristic methods for estimating generalized vertex
  median of a weighted graph. \emph{Operations Research} 16:955--961.

\bibitem[{Toregas et~al.(1971)Toregas, Swain, ReVelle, \protect\BIBand{}
  Bergman}]{Toregas1971}
Toregas CR, Swain RW, ReVelle CS, Bergman L (1971) The location of emergency
  service facilities. \emph{Operations Research} 19:1363--1373.

\bibitem[{Tribune(2017)}]{DhakaT}
Tribune D (2017) '999' emergency services begin. Technical report, Dhaka
  Tribune,
  \urlprefix\url{https://www.dhakatribune.com/bangladesh/2017/12/12/999-emergency-services-begin/}.

\bibitem[{Trudeau et~al.(1989)Trudeau, Rousseau, Ferland, \protect\BIBand{}
  Choquette}]{Trudeau1989}
Trudeau P, Rousseau Jm, Ferland JA, Choquette J (1989) An operations research
  approach for the planning and operation of an ambulance service.
  \emph{Information Systems and Operational Research} 27:95--113.

\bibitem[{Vairaktarakis \protect\BIBand{} Kouvelis(1999)}]{Vairaktarakis1999}
Vairaktarakis GL, Kouvelis P (1999) Incorporation dynamic aspects and
  uncertainty in 1-median location problems. \emph{Naval Research Logistics}
  46:147--168.

\bibitem[{van Barneveld(2016)}]{VanB2016}
van Barneveld T (2016) The minimum expected penalty relocation problem for the
  computation of compliance tables for ambulance vehicles. \emph{INFORMS
  Journal on Computing} 28(2):370--384,
  \urlprefix\url{http://dx.doi.org/10.1287/ijoc.2015.0687}.

\bibitem[{van Barneveld et~al.(2017)van Barneveld, van~der Mei,
  \protect\BIBand{} Bhulai}]{VanB2017}
van Barneveld T, van~der Mei R, Bhulai S (2017) Compliance tables for an ems
  system with two types of medical response units. \emph{Computers and
  Operations Research} 80:68 -- 81, ISSN 0305-0548,
  \urlprefix\url{http://dx.doi.org/https://doi.org/10.1016/j.cor.2016.11.013}.

\bibitem[{Vlahogianni et~al.(2014)Vlahogianni, Karlaftis, \protect\BIBand{}
  Golias}]{Vlahogianni2014}
Vlahogianni EI, Karlaftis MG, Golias JC (2014) Short-term traffic forecasting:
  Where we are and where we are going. \emph{Transportation Research Part C}
  43:3--19.

\bibitem[{Wadud(2017)}]{Wadud2017}
Wadud M (2017) Cheap solar ambulances to speed into service in rural
  bangladesh. \emph{Thomson Reuters Foundation}
  \urlprefix\url{https://www.reuters.com/article/us-bangladesh-solar-ambulance/cheap-solar-ambulances-to-speed-into-service-in-rural-bangladesh-idUSKBN15T1AP}.

\bibitem[{Westgate et~al.(2016)Westgate, Woodard, Matteson, \protect\BIBand{}
  Henderson}]{Westgate2016}
Westgate BS, Woodard DB, Matteson DS, Henderson SG (2016) Large-network travel
  time distribution estimation for ambulances. \emph{European Journal of
  Operational Research} 252:322--333.

\bibitem[{Wood(1993)}]{Wood1993}
Wood RK (1993) Deterministic network interdiction. \emph{Mathematical and
  Computer Modelling} 17:1--18.

\bibitem[{Woodard et~al.(2017)Woodard, Nogin, Koch, Racz, Goldszmidt,
  \protect\BIBand{} Horvitz}]{Woodard2017}
Woodard D, Nogin G, Koch P, Racz D, Goldszmidt M, Horvitz E (2017) Predicting
  travel time reliability using mobile phone gps data. \emph{Transportation
  Research Part C: Emerging Technologies} 75:30 -- 44, ISSN 0968-090X,
  \urlprefix\url{http://dx.doi.org/https://doi.org/10.1016/j.trc.2016.10.011}.

\bibitem[{Xu et~al.(2010)Xu, Caramanis, \protect\BIBand{} Mannor}]{Xu2010}
Xu H, Caramanis C, Mannor S (2010) Robust regression and lasso. \emph{IEEE
  Transactions in Information Theory} 56:3561--574.

\bibitem[{Zeng \protect\BIBand{} Zhao(2013)}]{Zeng2013}
Zeng B, Zhao L (2013) Solving two-stage robust optimization problems using a
  column-and-constraint generation method. \emph{Operations Research Letters}
  41:457--461.

\bibitem[{Zhang et~al.(2014)Zhang, Chenyang, Kan, Cao, Peng, Xu,
  \protect\BIBand{} Wang}]{Zhang2014}
Zhang Y, Chenyang Y, Kan H, Cao J, Peng L, Xu J, Wang W (2014) Effect of
  ambient temperature on emergency department visits in shanghai, china: a time
  series study. \emph{Environmental Health} 13:1--8.

\bibitem[{Zhang \protect\BIBand{} Li(2015)}]{Zhang2015}
Zhang ZH, Li K (2015) A novel probabilistic formulation for locating and sizing
  emergency medical service stations. \emph{Annals of Operations Research}
  229:813--835.

\bibitem[{Zhou \protect\BIBand{} Matteson(2016)}]{Zhou2016}
Zhou Z, Matteson DS (2016) Predicting melbourne ambulance demand using kernel
  warping. \emph{Ann. Appl. Stat.} 10(4):1977--1996,
  \urlprefix\url{http://dx.doi.org/10.1214/16-AOAS961}.

\bibitem[{Zhou et~al.(2015)Zhou, Matteson, Woodard, Henderson,
  \protect\BIBand{} Micheas}]{Zhou2015}
Zhou Z, Matteson DS, Woodard DB, Henderson SG, Micheas AC (2015) A
  spatio-temporal point process model for ambulance demand. \emph{Journal of
  the American Statistical Association} 110:6--15.

\end{thebibliography}


\newpage
\ECSwitch
\ECHead{Electronic Companion}

\section{Demand for emergency transportation}\label{EmgDemand}

In this section, we provide a descriptive analysis of our census and survey data (\ref{ED:Descriptive}) and present our methodology for estimating the components of emergency transport demand described in Section~\ref{annualtrips} (\ref{ED:DemandEst}).

\subsection{Descriptive analysis}\label{ED:Descriptive}

\subsubsection{Census data.}\label{ED:Census}

\begin{figure}[t]
\subfigure[\  Literacy rate \label{MapLit}]{
\includegraphics[width=.27\textwidth]{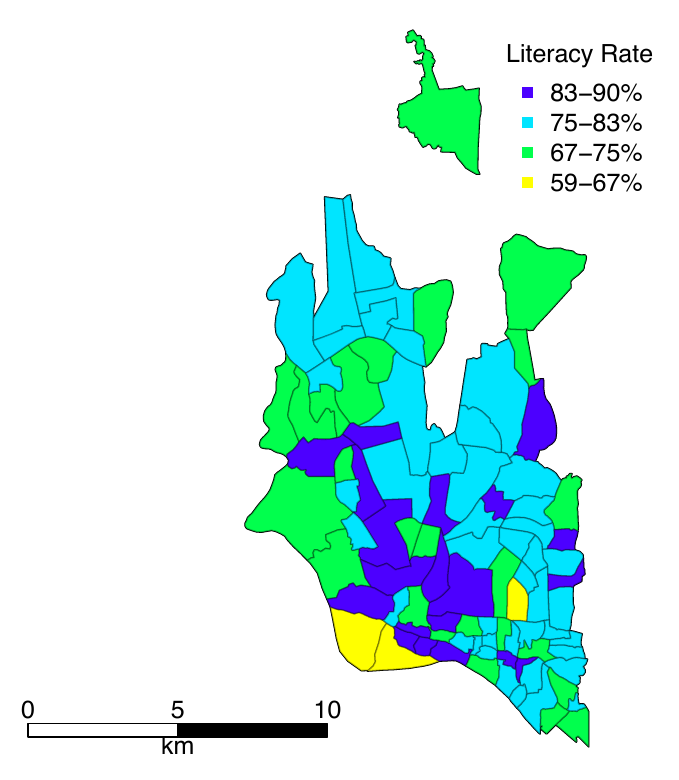}}
\hspace*{1em}
\subfigure[\  Pukka housing \label{MapHouse}]{
\includegraphics[width=.185\textwidth]{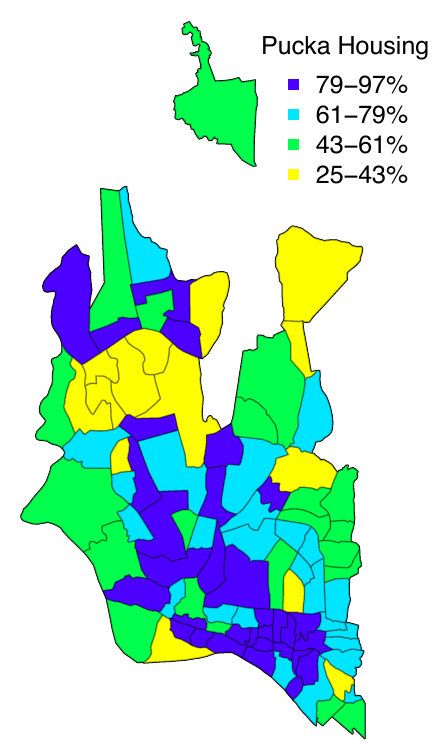}}
\hspace*{1em}
\subfigure[\  Female marriage rate \label{MapToilet}]{
\includegraphics[width=.21\textwidth]{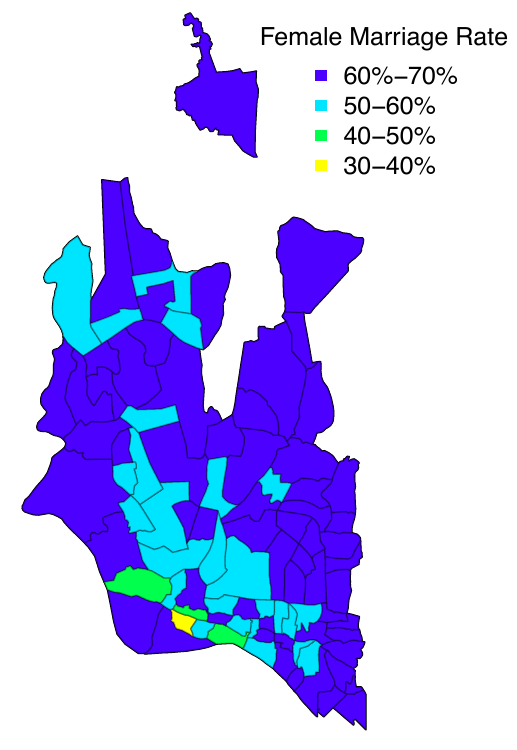}}
\hspace*{1em}
\subfigure[\ Household size\label{MapWater}]{
\includegraphics[width=.185\textwidth]{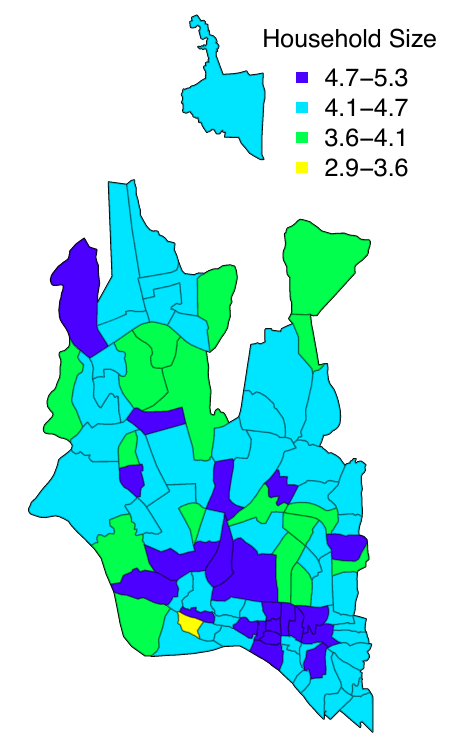}}
\caption{Ward-based rates for four demographic characteristics. A Pukka house is a solid permanent dwelling usually made from brick or stone that is reflective of higher socioeconomic status. \label{Maps}}
\end{figure}

We obtained the 2011 Dhaka census from the Bangladesh Bureau of Statistics. The census includes detailed demographic information for each of Dhaka's $92$ official wards (census tracts). Dhaka occupies a very small area of roughly $300$ km$^2$ with an estimated population of 7.35 million in 2011 (8.95 million in 2016), which is a slightly larger population than New York City in under 40\% of the area. Table~\ref{Census} provides summary statistics for key census characteristics and
Figure~\ref{Maps} illustrates the variation in four demographic characteristics across Dhaka's $92$ wards.

\setlength{\extrarowheight}{1.5pt}
\begin{table}[H]
\caption{Demographic summary statistics across Dhaka's 92 wards.\label{Census}}
\begin{threeparttable}
\centering
\begin{tabular}{l c c c c c}
\noalign{\smallskip} \hline
Characteristic & & Dhaka  & & Individual ward & \\
&& (all wards) & Minimum & Mean & Maximum \\ \hline
Population (2011 census) & & 7,349,324 & 18,170 & 79,884 & 228,870 \\
Average household size & & 4.3 & 3.0 & 4.4 & 5.3 \\
Male-female population ratio & & 1.2 & 0.9 & 1.3 & 2.5 \\
Population under 19 ($\%$) &&  36.2 & 22.9 & 35.5 & 43.5 \\
Population over 60 ($\%$) && 4.5 & 2.4 & 4.5 & 7.6 \\
Married ($\%$) && 59.0 & 29.4 &  57.4 & 66.7 \\
Literacy ($\%$) && 73.7 & 52.5 & 74.9 & 90.4  \\
Pukka$^*$ house ($\%$) && 58.4 & 24.9 & 66.4 & 96.5 \\
Jupri$^*$ house ($\%$) && 2.1 & 0 & 1.8 & 11.2 \\
Sanitary toilet ($\%$) && 58.0 & 7.4 & 60.3 & 98.1 \\
Electricity ($\%$) && 98.4 & 92.5 & 98.8 & 99.9 \\
Rent-free home ($\%$) & & 3.4 & 0.4 & 3.4 & 14.7 \\
Male-female employment ratio & & 2.0 & 0.3 & 2.8 & 12.8 \\
\hline \noalign{\smallskip}
\end{tabular}
\begin{tablenotes}
\scriptsize
\item $^*$Note that a Pukka house is a solid permanent dwelling usually made from brick or stone that is reflective of higher socioeconomic status, while a Jupri house is a temporary dwelling typically made from tin and other available supplies.
\end{tablenotes}
\end{threeparttable}
\end{table}

\subsubsection{Survey data.}\label{ED:Survey}

We obtained data from $2,808$ surveys administered by physicians to patients arriving at emergency departments (EDs) in $16$ major hospitals ($9$ private, $7$ government) across Dhaka. The survey had $14$ questions (see Table~\ref{Survey}) and was administered over $30$ days between July 7, 2014 and August 25, 2014. The survey data includes the chief complaint, date, time, and ward location of the emergency, mode and cost of transportation, and the time of arrival at the ED. Our survey data is unique because it includes patient travel data and as a result, we are able to provide insights on the current emergency medical system that have not been previously captured.
\begin{table}[h!]
\caption{Hospital survey questions and response options for all 2808 surveys.\label{Survey}}
\centering
\resizebox{30em}{!}{%
\begin{tabular}{l c l c c }
\noalign{\smallskip} \hline \noalign{\smallskip}
Question && Response Type & &   \\
\noalign{\smallskip} \hline \noalign{\smallskip}
P1. What was the approximate time of the emergency? && Free text &  & \\ \noalign{\smallskip} \hline \noalign{\smallskip}
P2. When did you decide to leave for the hospital? && Free text & & \\ \noalign{\smallskip} \hline \noalign{\smallskip}
P3. Which ward did you leave from? && Ward Number & & \\ \noalign{\smallskip} \hline \noalign{\smallskip}
P4. What time did you leave? && Free text & & \\  \noalign{\smallskip} \hline \noalign{\smallskip}
\multirow{7}{*}{P5. What was your method of transportation?} && A. Own car & & \\ && B. Rental car & & \\ && C. Rickshaw & & \\ && D. Ambulance & & \\ && E. CNG & & \\ && F. Taxicab & & \\ && G. Other & & \\ \noalign{\smallskip} \hline \noalign{\smallskip}
\multirow{4}{*}{P6. What was the cost of transportation (in BDT)} && A. $<100$ & & \\ && B. 100-500 & & \\ && C. 500-1000 & & \\ && D. 1000+ & & \\ \noalign{\smallskip} \hline \noalign{\smallskip}
\multirow{2}{*}{P7. Do you have a mobile phone?} && A. Yes & & \\ && B. No & & \\ \noalign{\smallskip} \hline \noalign{\smallskip}
\multirow{2}{*}{P8. Do you know how to contact an ambulance?} && A. Yes & & \\ && B. No & & \\ \noalign{\smallskip} \hline \noalign{\smallskip}
P9. Why did you not take an ambulance? && Free text & & \\ \noalign{\smallskip} \hline \noalign{\smallskip}
P10. Why did you come to this hospital? && Free text & & \\ \noalign{\smallskip} \hline \noalign{\smallskip}
H1. Name of hospital && Free text & & \\ \noalign{\smallskip} \hline \noalign{\smallskip}
H2. What was the arrival time of the patient? && Free text & & \\ \noalign{\smallskip} \hline \noalign{\smallskip}
H3. What is the general type of injury/complication? && Free text & & \\ \noalign{\smallskip} \hline \noalign{\smallskip}
H4. What time did the patient first receive treatment? && Free text & & \\ \noalign{\smallskip} \hline \noalign{\smallskip}
\end{tabular}}
\end{table}

\begin{figure}[t]
\begin{center}
\subfigure[\ Number of trips and cost for each mode  \label{CostMode}]{
\includegraphics[width=.44\textwidth]{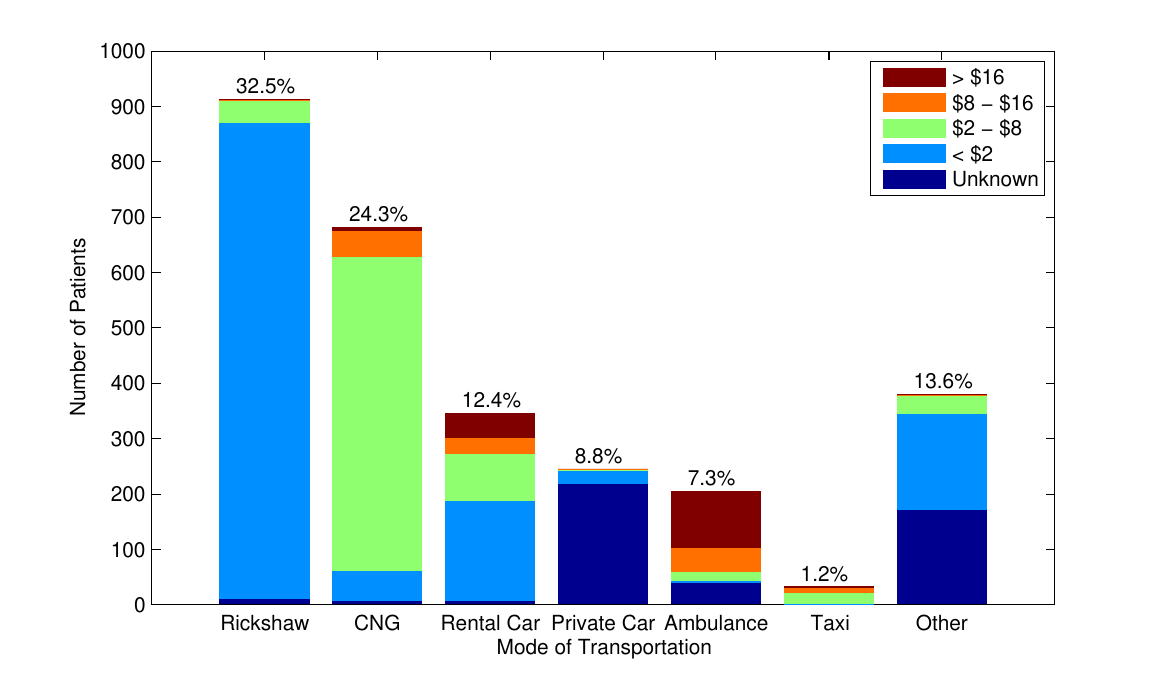}}
\subfigure[\ Mode of transportation according to severity \label{EmergencyMode}]{
\includegraphics[width=.525\textwidth]{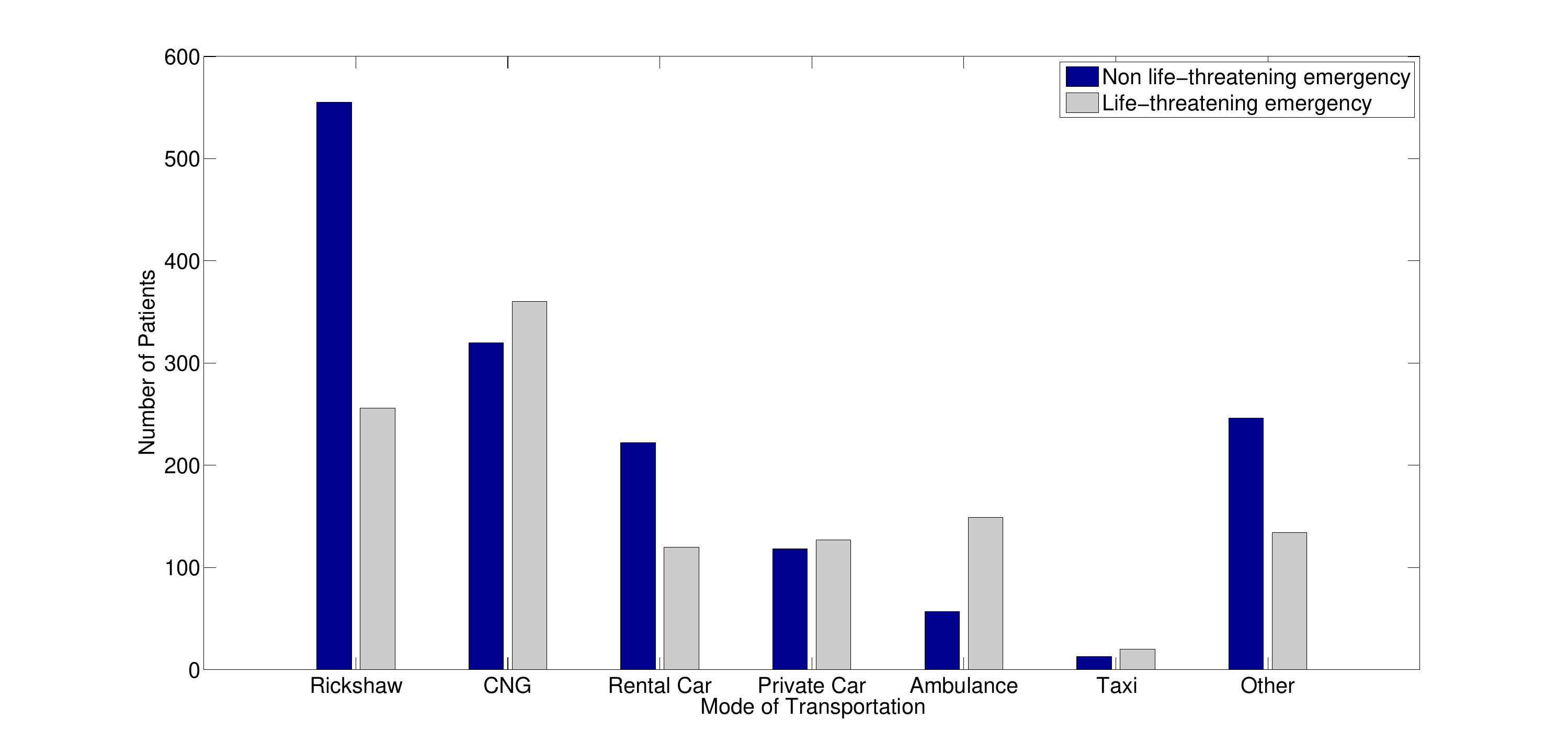}}
\caption{Histograms of trip characteristics for each mode of transportation.}
\end{center}
\end{figure}

Figure~\ref{CostMode} displays the various modes of transportation taken by patients and the costs incurred for each mode. Traditional ambulance vans were one of the least used modes of transport, comprising only $7.3\%$ of all trips. Of the survey respondents who answered the question ``Why did you not use an ambulance'', 16\% indicated that they tried but it was not available and 7\% cited slow response times, both of which are issues that can be addressed using our approach. Another major impediment was cost. Ambulances were found to be the most expensive mode of transportation with trips typically costing more than $16$ US dollars (USD). For context, the average annual income in Bangladesh is $1,260$ USD  \citep{HIES}. In contrast, rickshaws and auto-rickshaws, the two cheapest modes of transportation, comprised $34\%$ and $25\%$ of all trips, respectively. 

Aside from cost, another possible explanation for low ambulance utilization is that ambulances were one of the slowest modes of transportation, while rickshaws, private cars, and other, which includes walking, were the fastest modes of transportation. Table~\ref{ModeDist} provides a breakdown of the inter-ward distance travelled by each mode of transportation. Higher rickshaw usage, especially for short trips, is likely because rickshaws are readily available at all times and at nearly any location.
\begin{table}
\caption{A breakdown of the inter-ward distance travelled by each mode of transportation.\label{ModeDist}}
\centering
\begin{tabular}{l c c c c}
\noalign{\smallskip} \hline \noalign{\smallskip}
&&  Trips within & Median travel distance & Median travel distance  \\
Mode && ward (\%) & for all trips (m) & for out of ward trips (m)  \\ \hline \noalign{\smallskip}
Rickshaw && 30.2 & 1358 & 1544 \\
CNG && 6.8 & 4749 & 5233 \\
Rental Car && 7.8 & 6041 & 7018 \\
Private Car && 41.8 & 1367 & 2462 \\
Ambulance && 9.2 & 3379 & 3670 \\
Taxi && 5.0 & 5737 & 6262 \\
Other && 24.0 & 2017 & 4041 \\
\hline \noalign{\smallskip}
\end{tabular}
\end{table}

Although overall ambulance utilization in Dhaka is low, Figure~\ref{EmergencyMode} shows that it is the mode with the highest proportion of trips for life-threatening emergencies (classified by the attending physician). More than two thirds of all ambulance trips are for life-threatening emergencies, compared to less than one third of rickshaw trips. In life-threatening emergencies, ambulances become the third most common mode of transportation. This data suggests that patients recognize the importance of ambulances and are willing to use them for life-threatening emergencies. These findings also reinforce the importance of considering multiple vehicles types in LMICs. 

\subsection{Estimating the annual number of emergency trips}\label{ED:DemandEst}

In this section, we present our methodology for estimating the three components that comprise the annual number of emergency trips: the population in ward $w$ at time $\tau$ (\ref{Apd:DemandEst1}), the average annual number of ED visits per person (\ref{Apd:DemandEst2}), and the proportion of ED visits from ward $w$ arriving via mode $m$ (\ref{Apd:DemandEst3}).

\subsubsection{Estimating population ($n_{w,\tau}$).}\label{Apd:DemandEst1}

\begin{figure}[t]
\subfigure[\ $n_{w,D} $ \label{DayPop}]{
\includegraphics[width=.19\textwidth]{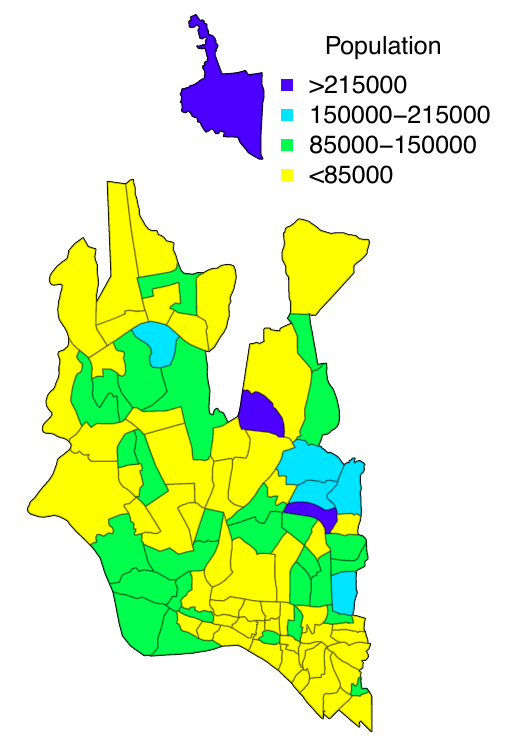}}
\hspace*{2em}
\subfigure[\ $n_{w,N}$\label{NightPop}]{
\includegraphics[width=.19\textwidth]{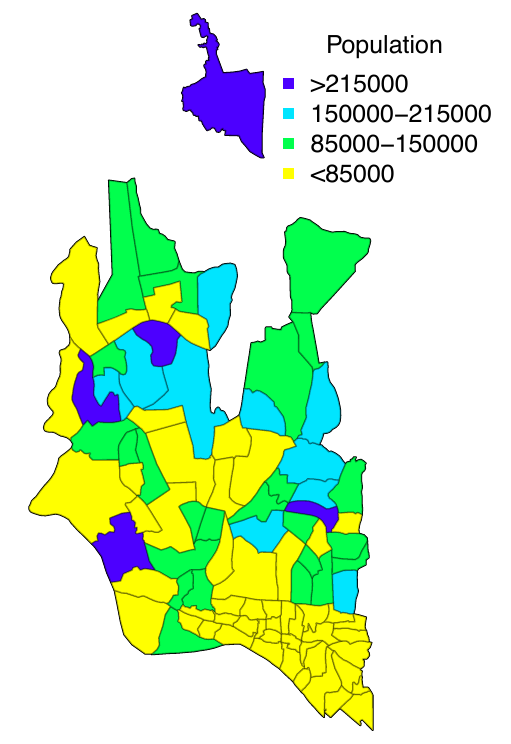}}
\hspace*{2em}
\subfigure[\ $\delta_{w,V}$ \label{Pred1}]{
\includegraphics[width=.22\textwidth]{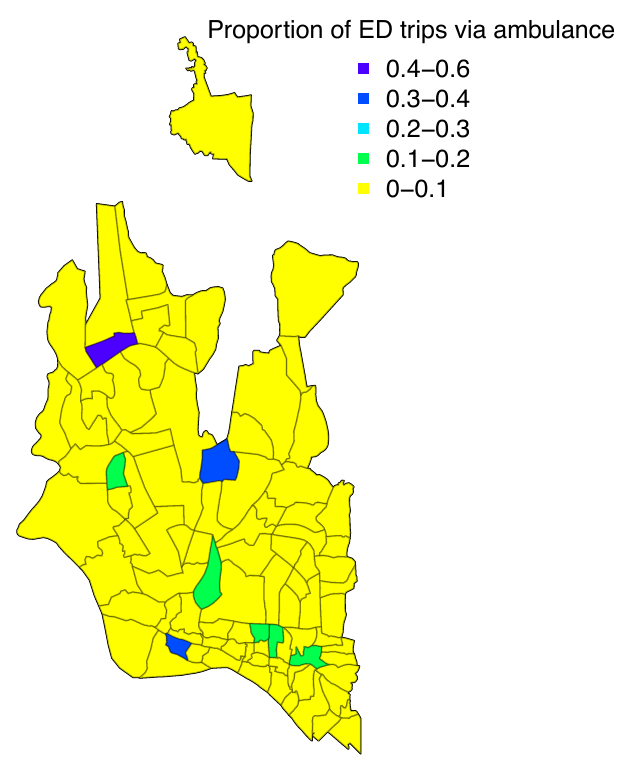}}
\hspace*{1em}
\subfigure[\  $\delta_{w,S}$\label{Pred2}]{
\includegraphics[width=.22\textwidth]{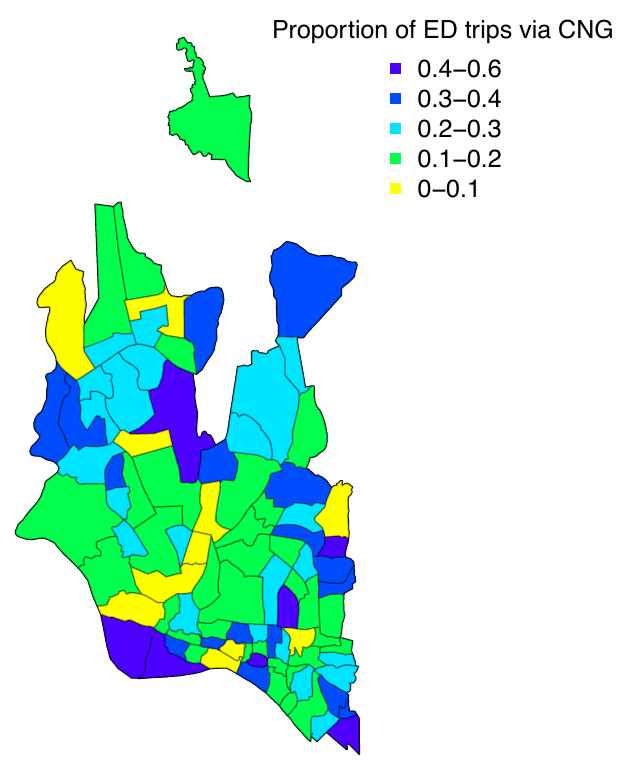}}
\caption{The geographical distribution of several estimated components of demand.}
\end{figure}

In this section, we estimate both the daytime and nighttime population for each ward (i.e., $\tau\in\{D,N\}$). Dhaka has a major difference in the spatial distribution of the daytime and nighttime population, due to daily migration. The magnitude of the daily migration out of Dhaka is estimated to be over $700,000$ as many people leave the city during the day to work in the surrounding industrial areas. We estimate the daytime population in each ward from the Earthquake Vulnerability Assessment of Dhaka, which was conducted by the Government of Bangladesh with support from the United Nations Development Programme. This assessment estimates the total population in each ward during the daytime working hours. In $2008$, the total daytime population in Dhaka was estimated to be $6.63$ million people. The nighttime population in each ward is obtained directly from the $2011$ census and was estimated to be $7.35$ million people. The population of Dhaka has been consistently growing at a rate of approximately $320,000$ people per year \citep{Streatfield2008}. Under the assumption that each ward is growing at a rate proportional to its population, we estimate the total $2016$ daytime and nighttime population in Dhaka to be $8.24$ and $8.95$ million, respectively. Figures EC.\ref{DayPop} and EC.\ref{NightPop} illustrate the estimated geographical distribution of the daytime and nighttime populations, respectively.

\subsubsection{Estimating the average annual number of ED visits per person ($\xi$).}\label{Apd:DemandEst2} 

In this section, we leverage published research from South Asian cities to estimate the average annual number of ED visits per person. Given data limitations and the coarseness of previous studies, we cannot generate ward-specific rates and instead settle on a single estimate for the population.

A recent study of a ED arrivals at a ``specialty corporate hospital" in Dhaka, found an average of $10,000$ ED visits each year \citep{Karim2009}. This result requires careful interpretation because specialty hospitals can be very expensive and serve only a limited population. To the best of our knowledge, there is no further data on ED visits in Dhaka or other cities in Bangladesh.

To supplement this lack of data, we estimate the number of ED visits using data from other similar South Asian cities. A study of three major government hospitals in Kirachi, Pakistan found an average of $70,000-100,000$ annual ED visits per hospital \citep{Raftery1996}. A similar study of two major hospitals in New Dehli, India found that a private hospital with free emergency services received $30,000$ annual ED visits, while a government funded hospital with free services received over $100,000$ annual ED visits \citep{PoSaw1998}.

From these reports, we estimate the number of annual ED visits for a government funded hospital to be between $70,000-100,000$ and we estimate the number of annual ED visits for a private hospital to be between $10,000-30,000$. Dhaka has $87$ hospitals with EDs, of which $19$ are government funded. From this information, we estimate the number of annual ED visits to be between $2.08-4.15$ million. Given that Dhaka has a population of $8.95$ million, the visit rate is between $230-460$ per $1000$ persons. Therefore, the average number of annual ED visits per capita, $\xi$, is estimated to be between $0.23-0.46$.

It is difficult to put these numbers into context because most LMIC countries do not collect data on annual ED visits. However, data is available for $19$ high income OECD countries and the average number of annual ED visits per capita across all countries is $0.31$ with a range from $0.07$ to $0.70$ \citep{Berchet2015}. These results require careful interpretation because they are from high income countries and they combine data from both rural and urban areas, which are known to have significant differences in ED visit rates. For example, in the US, urban areas have an annual rate of $0.32$ ED visits per capita as compared to $0.45$ ED visits per capita in rural areas. The only reliable data available for large urban areas is from New York City and Shanghai, which are similar to Dhaka in terms of population, but not in terms of culture or demographics. The annual per capita ED visit rates in New York City and Shanghai are $0.37$ \citep{Goins2015} and $0.33$ \citep{Zhang2014}, respectively.

\subsubsection{Estimating the proportion of ED visits from ward $w$ arriving via mode $m$ ($\delta_{w,m}$).}\label{Apd:DemandEst3}

In this section, we train a regularized (lasso) logistic regression model for predicting the proportion of ED visits from ward $w$ arriving via mode $m$ ($\delta_{w,m}$). In particular, we aim to estimate the proportion of ED visits arriving via ambulance van ($\delta_{w,V}$) and the proportion of ED visits arriving via small ambulance ($\delta_{w,S}$), for each ward in Dhaka, respectively. Although we did not explicitly assume that $\delta_{w,S} + \delta_{w,V} \leq 1$, our results naturally satisfied this inequality.

The nature of our data presents significant challenges for model training. Our survey data indicates that only $74$ of Dhaka's $92$ wards include at least one surveyed patient and as a result, we cannot directly estimate $\delta_{w,m}$ for each ward. Furthermore, only $32$ of $92$ wards include at least $20$ observations. Given that the overall ambulance usage is 7\%, roughly $20$ observations are required to ensure sufficient granularity in our estimations. One way to overcome the lack of data is to assume that $\delta_{w,m}$ is uniform across all $92$ wards and estimate a single value for each mode, $\delta_m$. Intuitively, assuming that $\delta_m$ is uniform across all wards is analogous to using population as a proxy for the annual number of emergency trips. We use this naive approach as a benchmark for our models. A second approach to overcome the lack of data is to group patients according to the ward in which the trip originated and calculate $\delta_{w,m}$ for each ward with $20$ or more patients (i.e., only $32$ wards). We employ this approach to transform our data and weight each $\delta_{w,m}$ by the number of observations (i.e., patients) from ward $w$. As a result, our final dataset includes $32$ ward observations comprising $1,843$ patients. This approach assumes that $\delta_{w,m}$ is constant across each ward and over time, which is a limitation because $\delta_{m}$ may be different for each individual patient and for different times of day. However, modeling emergency mode choice decisions at a fine spatiotemporal level requires very granular data, which does not exist in many developing countries. 

Grouping patients by ward as opposed to a patient level approach that treats each patient as a unique observation is beneficial for our application for three key reasons: 1) we are interested in estimating $\delta_{w,m}$ at the ward level, not the patient level, 2) we want our approach to be generalizable and this framework allows other regions in LMICs with only census data to apply our models, and 3) we require independent variables or features that are available for all $92$ wards. The only features available to us for all $92$ wards are from the census and we link this data to the survey data using the ward where the trip originated. In contrast, a patient level approach will cause all patients from the same ward to have identical independent variables, regardless of their mode choice.

The set of $27$ ward-level demographic features we use in each model was selected from the census data, which contains $104$ unique fields. Note that our models do not use individual-level features. To do this, we first remove all highly correlated ($R^2>0.85$) variables that appear to represent the same latent feature. In particular, we remove the minimum number of variables required to eliminate all pairwise correlations above $0.85$. Next, we combined variables to create new features that have been previously shown to correlate with ambulance demand. For example, the original data contained male population, female population, and total population, which are all highly correlated. We kept total population and created a new variable using the ratio of male to female population. After this procedure, a final set of $27$ demographic census features remained. 

Our data is well suited for logistic regression because our observations can be viewed as independent Bernoulli trials and modelled using a binomial distribution. Given the large set of features and the likelihood of overfitting, we consider a logistic regression model with L1-regularization (LASSO), where $\gamma$ denotes the regularization parameter. We optimize over $1000$ values of $\gamma$ between $0.0001$ and $0.01$. Recall, that the naive prediction approach mentioned above predicts a constant equal to the average $\delta_m$ across all wards. We also consider a weighed (by population) naive approach that predicts a constant equal to the weighted average $\delta_m$ across all wards.

We train our models using repeated $10$-fold cross validation, which partitions the data into ten sets: eight sets of three wards and two sets of four wards. Each set is used exactly once as the testing set, while the remaining $9$ are combined and used as the training set. We repeat this process $500$ times to reduce the variance in our estimations of model accuracy. We measure prediction accuracy using root mean squared error. Once the value of $\gamma$ that minimizes RMSE is determined through repeated cross validation, we train a final model using this $\gamma$ and all available data to estimate $\delta_{w,m}$ for all $92$ wards. All models were implemented using \texttt{R version 3.3.3.}

\begin{figure}[t]
\centering
\subfigure[\ $\delta_{w,v}$ \label{Err}]{
\includegraphics[width=.4\textwidth]{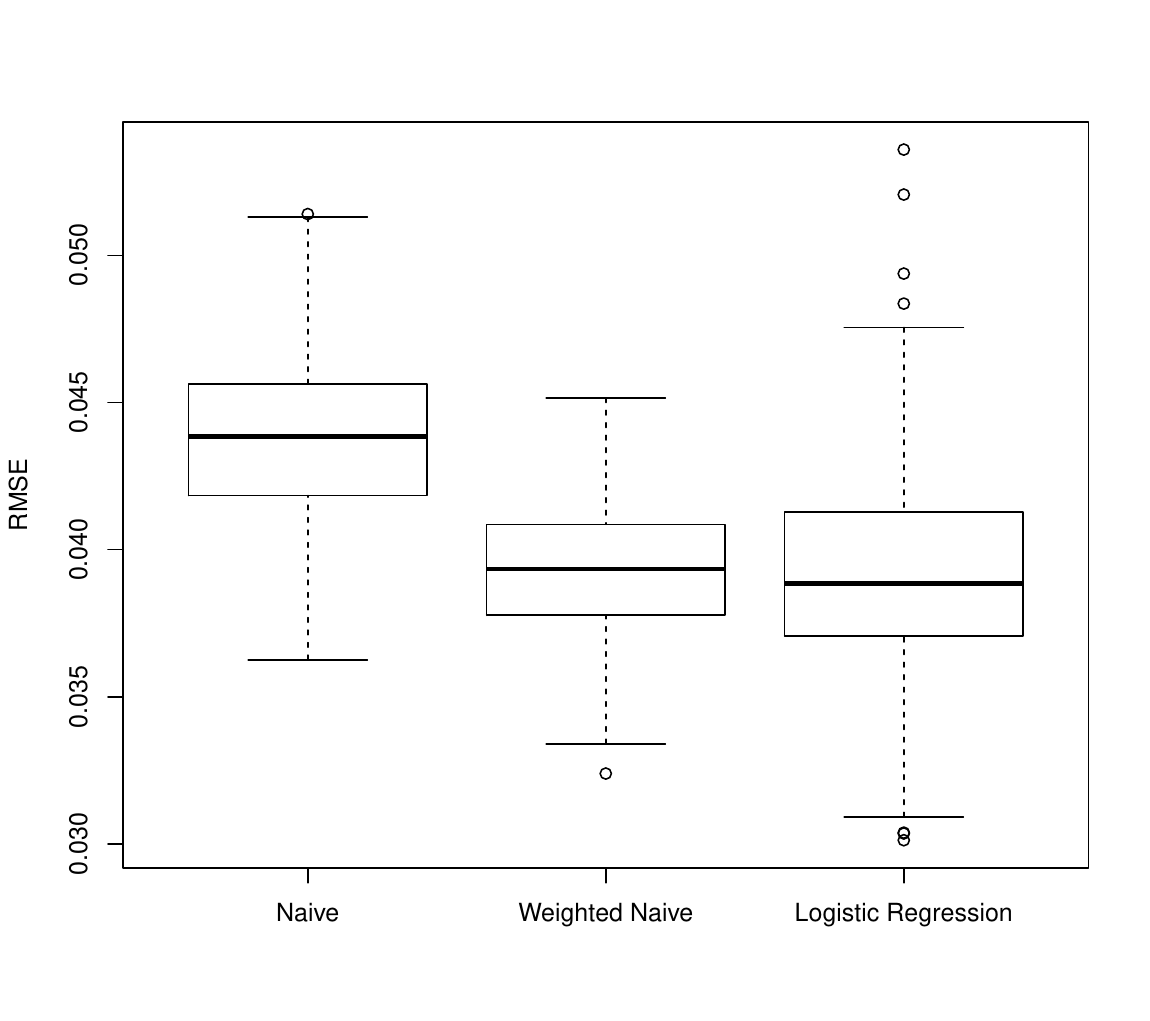}}
\subfigure[\ $\delta_{w,s}$ \label{Err2}]{
\includegraphics[width=.4\textwidth]{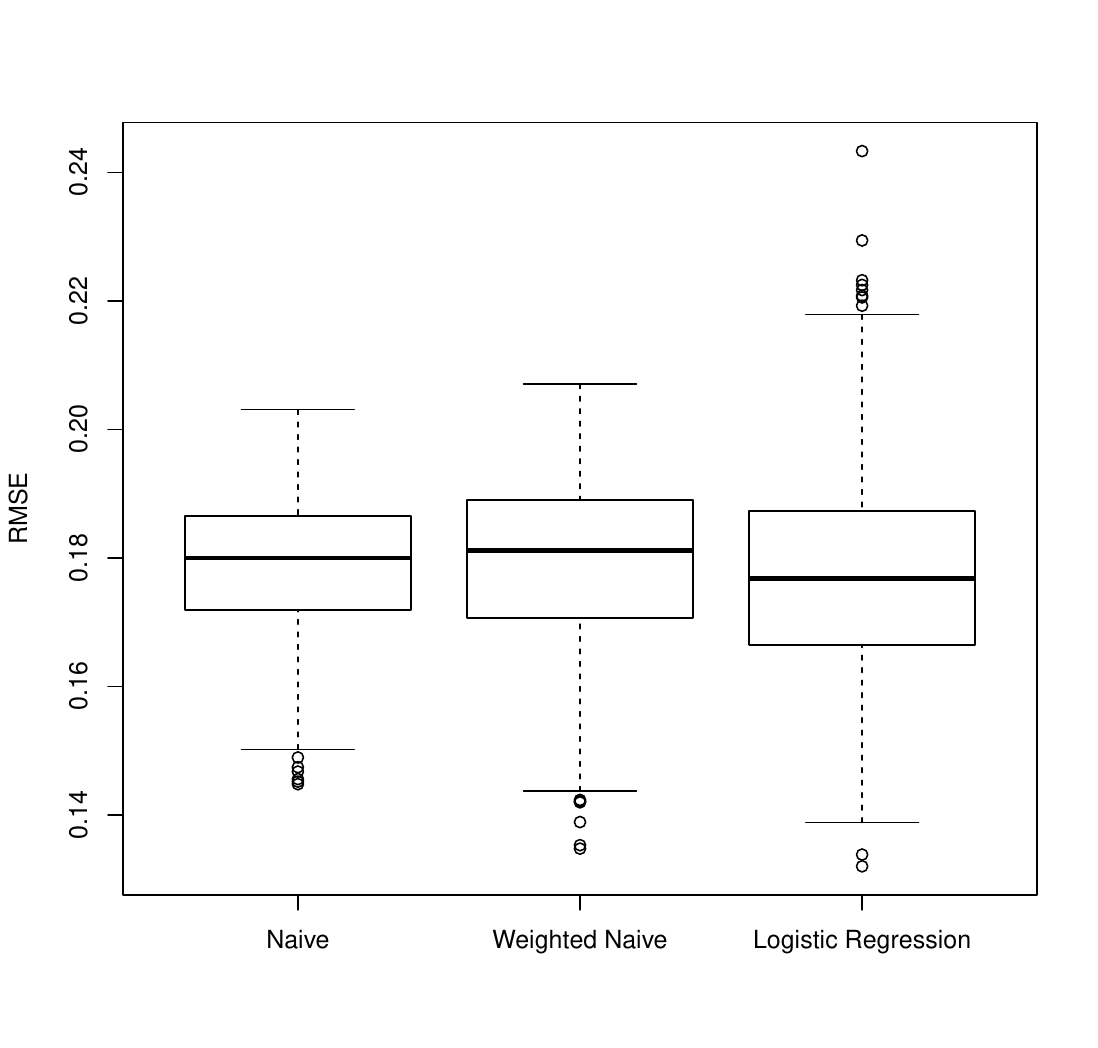}}
\caption{Comparison of RMSE between the logistic regression model and the naive approaches.}
\end{figure}

Figures EC.\ref{Err} and EC.\ref{Err2} display box plots of the RMSE distribution across the $500$ repetitions for $\delta_{w,V}$ and $\delta_{w,S}$, respectively. The solid black line indicates the median and the box indicates the interquartile range. The whiskers extend to $1.5$ times the interquartile range. For $\delta_{w,V}$, the median RMSE was $0.04385$ and $0.03935$ for the naive and weighted naive, respectively. The logistic regression model (with $\gamma=0.00183$) performed the best with a median RMSE of $0.03886$, corresponding to a $11.4\%$ improvement over the naive approach and a $1.2\%$ improvement over the weighted naive approach. For $\delta_{w,S}$, the median RMSE was $0.180$ and $0.181$ for the naive and weighted naive, respectively. The logistic regression model (with $\gamma=0.00183$) preformed the best with a median RMSE of $0.176$, corresponding to a $2.2$\% improvement over the naive approach and a $2.8\%$ improvement over the weighted naive approach. Both improvements were found to be statistically significant at $\alpha=0.05$ using the Wilcoxon signed-rank test. For both $\delta_{w,V}$ and $\delta_{w,S}$, the RMSE improvements from logistic regression were found to be statistically significant at $\alpha=0.05$ using the Wilcoxon signed-rank test. 

Although the logistic regression model only marginally improves upon the weighted naive approach, the model is able to provide insight into the demographic features that contribute to ambulance usage. Our final features (see Tables~\ref{FeatureResults} and~\ref{FeatureResults2}), consistent with prior literature, include measures of population (e.g., average household size), measures of social status (e.g., female marriage rate), and measures of economic status (e.g., access to electricity). Figures EC.\ref{Pred1} and EC.\ref{Pred2} display the model-predicted values of ED visits arriving via ambulance van ($\delta_{w,S}$) and small ambulance ($\delta_{w,S}$), respectively. We find that areas of higher socioeconomic status are more likely to use ambulance vans as compared to small ambulances. For example, the wards with the largest values of $\delta_{w,V}$ include areas with a high density of foreigners, government officials, and an area with major government offices, hospitals, and universities. In contrast, small ambulance use is highest in many of the city's outer wards, which include slums.

\begin{table}
\caption{Non-zero regression coefficients as determined by LASSO for ambulance vans ($\delta_{w,V}$).\label{FeatureResults}}
\centering
\begin{tabular}{l c r  }
\noalign{\smallskip} \hline \noalign{\smallskip}
Feature && Coefficient \\ \hline \noalign{\smallskip}
Intercept && $-4.065$ \\
Average household size (number of persons)& & $-0.085$\\
Ratio of male to female population &&  $1.358$ \\
Female marriage rate (\%) && $-0.014$\\
Access to electricity && $0.005$\\
\hline \noalign{\smallskip}
\end{tabular}
\end{table}

\begin{table}
\caption{Non-zero regression coefficients as determined by LASSO for small ambulances ($\delta_{w,S}$).\label{FeatureResults2}}
\centering
\begin{tabular}{l c r }
\noalign{\smallskip} \hline \noalign{\smallskip}
Feature && Coefficient\\ \hline \noalign{\smallskip}
Intercept && $-13.600$ \\
Ratio of male to female population && $0.062$ \\
Male marriage rate (\%) && $-0.002$ \\
Population between 0-19 (\%) && $0.0077$ \\
Population over 60 (\%) && $-0.149$ \\
Disability rate (\%) && $-0.390$ \\
Pukka house (\%) && $-0.028$ \\
Access to a sanitary toilet with seal (\%) && $0.003$\\
Access to electricity (\%) && $0.126$\\
Ratio of male to female employment && $0.019$\\
\hline \noalign{\smallskip}
\end{tabular}
\end{table}

Although we focused on predicting the probability that a patient chooses an ambulance van or small ambulance, given that they require transportation to an ED, our approach can can be readily adapted to focus on other modes of transportation, such as private cars, rickshaws, or motorcycles.

\section{Travel time analysis}\label{Apd:TT}

In this section, we describe the travel time data collection methodology (Section~\ref{TT:Data}) and develop machine learning models to predict the baseline travel time between any two locations in both road networks (Section~\ref{TT:Preds}). 

\subsection{Travel time data}\label{TT:Data}
We gathered vehicle location data using custom GPS devices and an accompanying Android mobile application developed by our collaborators. These devices were used by five volunteer citizens over 16 days from March 14, 2014 to June 13, 2014 and over 14 days from February 28, 2015 to April 2, 2015. All drivers were instructed to drive normally, using typical routes and speed. A map matching algorithm was developed to map the GPS data to edges on the road network \citep{Ahmed2015}.

\begin{figure}[t]
\centering
\subfigure[\ Samples per edge \label{SamplesPerEdge}]{
\includegraphics[width=.4\textwidth]{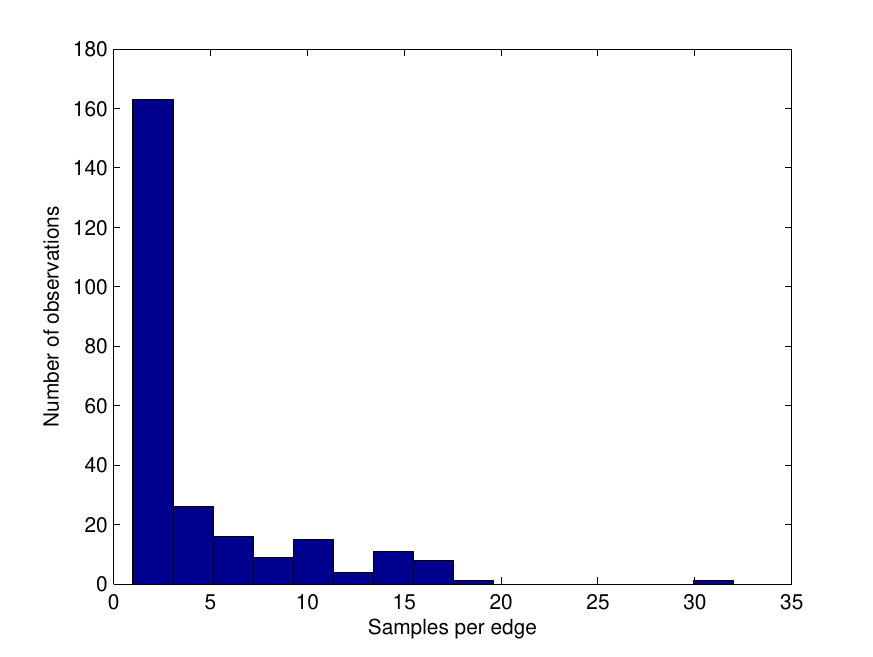}}
\subfigure[\ Edges per trip \label{EdgesPerTrip}]{
\includegraphics[width=.4\textwidth]{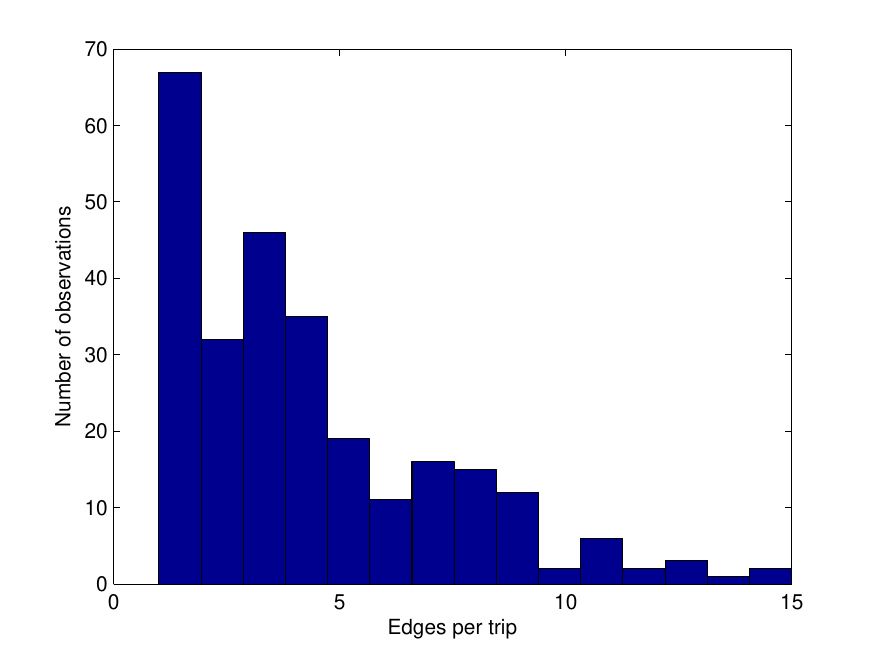}}
\caption{Histograms for collected GPS data. Only sampled edges are included in these figures.\label{SamplesE}}
\end{figure}

We obtained data for $269$ unique trips. A trip is defined as a path through the network from some origin node to some destination node. A destination node is defined as either one from which there is no subsequent GPS activity within 20 minutes on an edge emanating from that node or one with the last recorded GPS activity before the device was turned off by the driver. Trips ranged from 1 to 15 edges, with an average trip length of 4.1 edges. Edges in the network were present in a trip between 0 to 30 times; edges that were present in at least one trip had an average of 3.9 observations for a total of 1,103 edge observations (see Figure EC.\ref{SamplesE}). The median travel time of a trip was $592s$ (min: $10s$, max: $5543s$) while the median travel time on an edge was $105s$ (min: $5s$, max: $5062s$). Figures EC.\ref{SamplesPerEdge} and EC.\ref{EdgesPerTrip} display the number of data samples per sampled edge and the number of edges per trip in the ambulance road network, respectively.

\begin{figure}[t]
\centering
\includegraphics[width=0.5\textwidth]{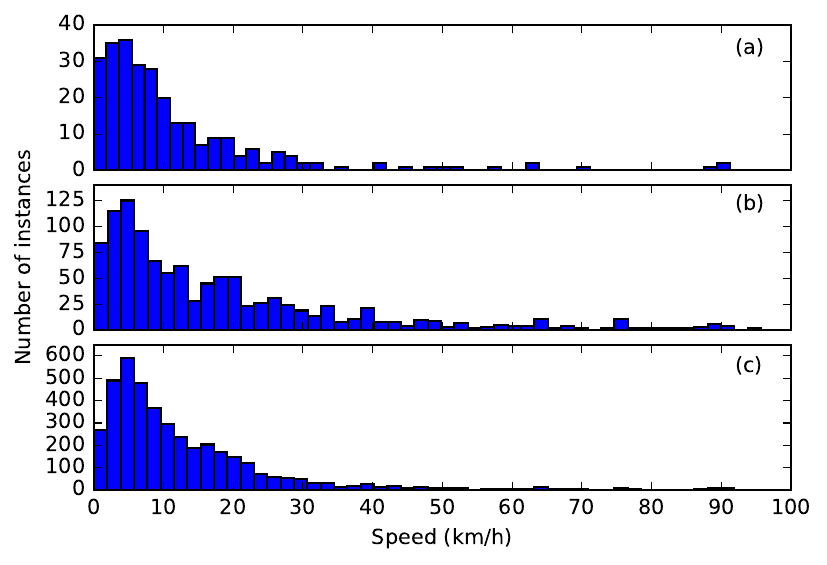}
\caption{Histograms of speed for (a) trip data, (b) edge data, and (c) sub-trip data.}\label{TVD}
\end{figure}

To predict the travel time between two nodes, one challenge we face is limited data. In particular, if we use trip data to train our models, we are limited to only 269 observations. On the other hand, if we use edge data, which is more plentiful and includes 1,103 observations, then we are unable to capture the delays caused at nodes between edges (i.e., intersections) or the impact of traveling through a ward because most edges lie wholly within one ward. To deal with this trade-off, we develop a modified bootstrapping method that simultaneously solves the limited data issue and the issues with using edge data. This bootstrapping method expands our dataset by partitioning each trip into all continguous \emph{sub-trips}. For example, a trip that begins at node $1$, visits nodes 2, then 3, and terminates at node $4$ (denoted $1-2-3-4$) would result in six sub-trips: $1-2$, $2-3$, $3-4$, $1-2-3$, $2-3-4$, and $1-2-3-4$. In other words, we include the original trip ($1-2-3-4$), all individual edges ($1-2$, $2-3$, $3-4$) and all sup-trips ($1-2-3$, $2-3-4$). This bootstrapping process results in a total of $4,086$ sub-trips, a $15$ times increase in the size of the training set. The new sub-trip data is not a direct replication of trip data because each sub-trip has unique features according to the origin/destination of that sub-trip. Figure~\ref{TVD} displays a histogram of speeds for trips, edges, and sub-trips. Note that the sub-trip data includes both trip and edge data. The average speed (standard deviation) for the trip, edge, and sub-trip data is 2.05 km/h (3.90), 3.30 km/h (5.31), and 2.45 km/h (3.64), respectively.

\begin{figure}[t]
\subfigure[\ Sub-trip\label{Ex}]{
\includegraphics[width=.3\textwidth]{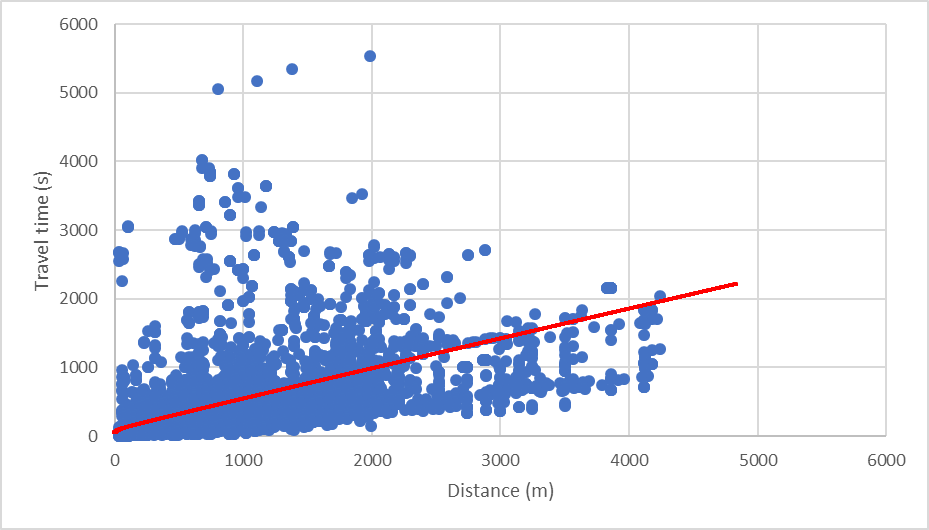}}
\subfigure[\ Trip \label{T}]{
\includegraphics[width=.3\textwidth]{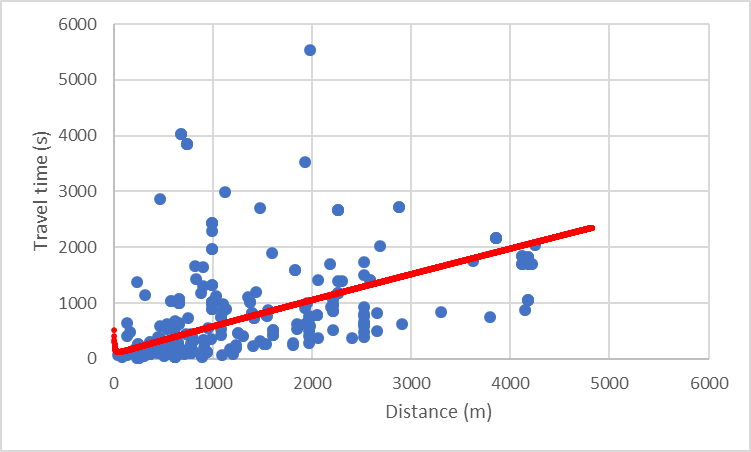}}
\subfigure[\  Edge \label{E}]{
\includegraphics[width=.3\textwidth]{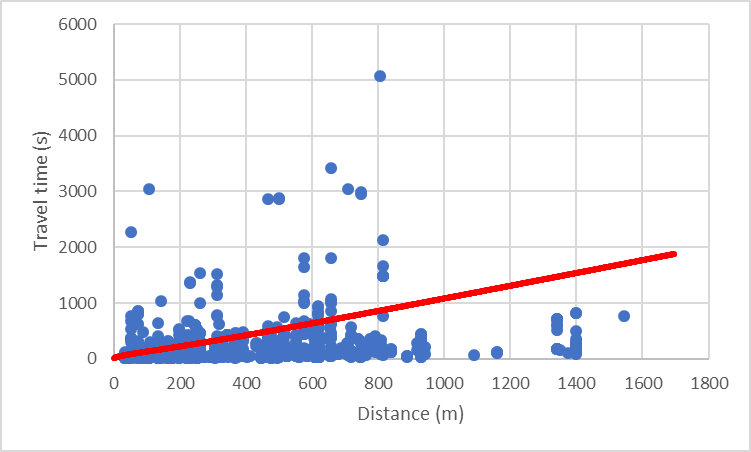}}
\caption{Scatter plots of distance vs. time with the fitted Kolesar model\label{TimeVDist}.}
\end{figure}

\begin{figure}[t]
\centering
\includegraphics[width=0.5\textwidth]{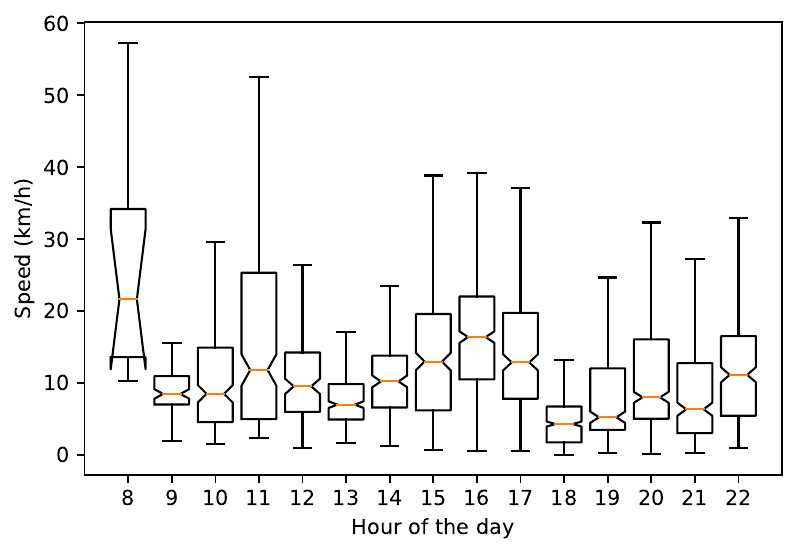}
\caption{Boxplot of speed for each hour of the day. The box represents the interquartile range (i.e., $IQR = Q3-Q1$), the red line indicates the median, the notches indicate the 95\% confidence interval around the median, and the upper (lower) whiskers correspond to $Q3 + 1.5*IQR$ ($Q1-1.5*IQR$).\label{SpeedBox}}
\end{figure}

Figure~\ref{TimeVDist} displays a scatter plot of distance vs. travel time for the sub-trip, trip, and edge data. The red curve is the \cite{Kolesar1975} model trained using all available data. Figure~\ref{SpeedBox} displays a boxplot of speed for each hour of the day where data was available. We find that speeds are slowest during the evening rush hour (i.e., 6pm and 7pm) and fastest in the mid-afternoon and nighttime/early morning. These results are consistent with our experience in Dhaka and with Google Maps traffic data (note that Google only started providing this service to Dhaka in late 2017/early 2018). \cite{McCormack2015} also found very similar results using data from the London Ambulance Service and this reference had been added to the paper.

\subsection{Travel time prediction}\label{TT:Preds}

Using the trip, edge, and sub-trip data, we compare four machine learning approaches for predicting the travel time in seconds between any two nodes (not necessarily adjacent) in the network. We use 73 features including the distance on the road network between the given nodes, the day of week, and time of day. In addition, for both the origin and destination wards, we include building-type information (e.g., the number of commercial or industrial buildings) and the $27$ demographic features used to predict emergency demand (see~\ref{Apd:DemandEst3}). The target is a real number that denotes the travel time in seconds between the two nodes. 

We compared the accuracy of four popular machine learning models: AdaBoost, Random Forest, linear regression with L1-regularization (LASSO), and K-nearest neighbors (KNN). For AdaBoost, we optimize the learning rate over $\{0.0001,0.001,0.01,0.1,1\}$ and number of weak learners over $\{100,250,500,750,1000\}$; for Random forest, we optimize the number of trees over $\{100,250,500,750,1000\}$; for LASSO, we optimize the regularization parameter over $\{0.0001,0.001,0.01,0.1,1\}$; for KNN we optimize the number of neighbors over $\{100,250,500,750,1000\}$. We train our models using repeated $10$-fold cross validation and we repeat the cross validation process $100$ times. We used a nested 3-fold cross validation loop on the training set for hyper-parameter tuning. In particular, as part of 10-fold cross validation, we partition the data into two sets: set1 and the test set. Set1 comprises 90\% of the data, while the test set comprises 10\%. We then split set1 into two sets: a training set and a validation set. We use 3-fold cross validation on set1 to conduct three training-validation instances and we use the average validation set error to choose our hyperparameters. Once the final hyperparameters have been selected, we apply the final model to the test set (that was not used as part of the model selection or fitting process in any way) to estimate generalization. We repeat the entire process 100 times to obtain more robust estimates.

We measure prediction accuracy using root mean squared error (RMSE). Once the final hyperparameters are determined, we train the model using all available data to obtain the final predictions of travel time, which are used in our optimization model. A previously developed model and two naive approaches serve as a baseline.  The first naive approach, \emph{Naive S}, predicts a constant equal to the average travel time from the empirical data, and the second, \emph{Naive D}, is a simple linear regression model fit to distance only. We also compared our machine learning approaches to the model developed by \cite{Kolesar1975} that we trained using the maximum likelihood methodology proposed by \cite{Budge2010}. All experiments were implemented using \texttt{Python 3.5}. 

\begin{figure}[t]
\centering
\includegraphics[width=.6\textwidth]{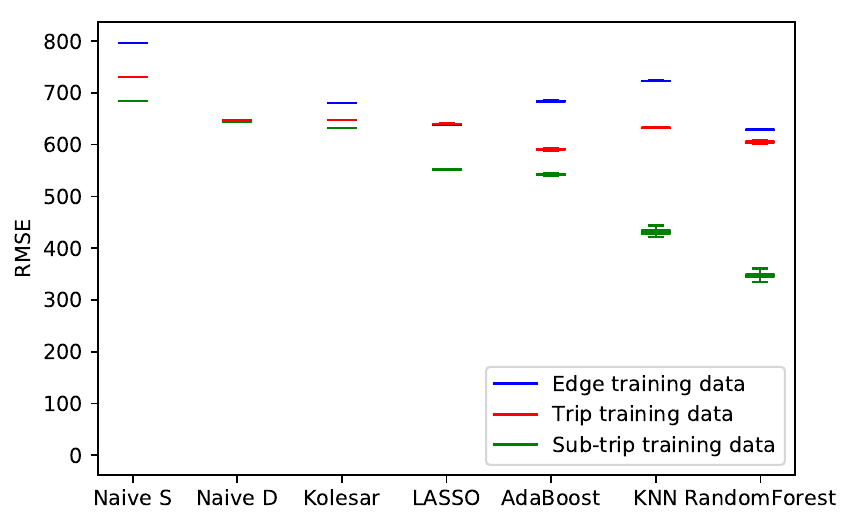}
\caption{Mean-squared error results for the models tested on sub-trip data.}\label{ExpandedTest}
\end{figure}

\begin{figure}[t]
\centering
\includegraphics[width=.5\textwidth]{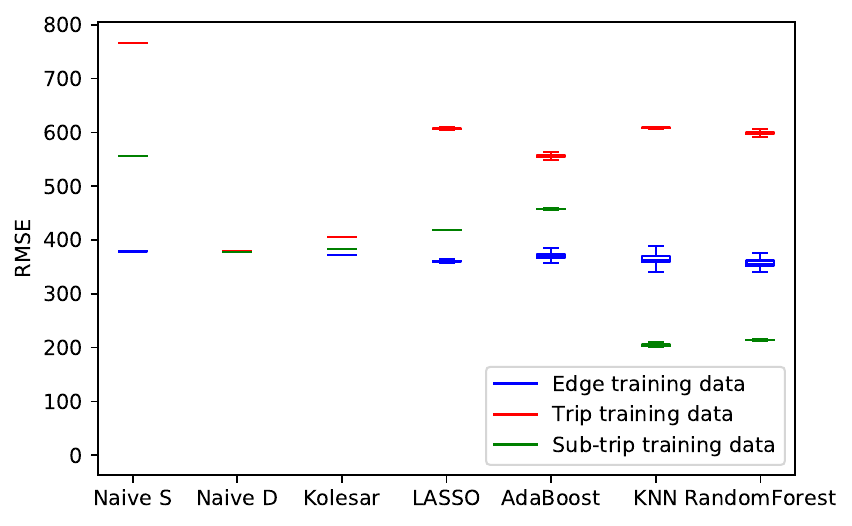}
\caption{Mean-squared error results for the models tested on edge data.}\label{EdgeTest}
\end{figure}

\begin{figure}[t]
\centering
\includegraphics[width=.5\textwidth]{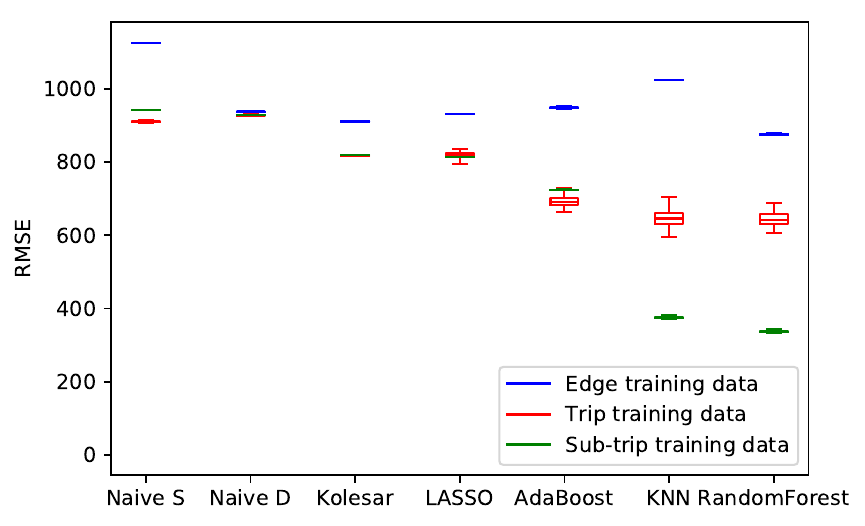}
\caption{Mean-squared error results for the models tested on trip data.}\label{TripTest}
\end{figure}

Figure~\ref{ExpandedTest} displays a boxplot of the root mean squared error (RMSE) distribution across $100$ repetitions for each of the prediction models tested on sub-trip data. The median RMSE for the Naive D (Naive S) approach was $648$s ($797$s) when trained on edge data, $647$s ($731$s) when trained on trip data, and $645$s ($684$s) when trained on sub-trip data. The random forest model performed the best with a median RMSE of $629$s, $605$s, and $348$s corresponding to improvements of 3\% (21\%), 6\% (17\%), and 46\% (49\%) over the Naive D (Naive S) approach when trained on edge, trip, and sub-trip data, respectively. All improvements were found to be statistically significant at $\alpha=0.01$ using the Wilcoxon signed-rank test. Figures~\ref{EdgeTest} and \ref{TripTest} display boxplots of the RMSE for each of the prediction models tested on edge and trip data, respectively. These results depict a similar finding: a random forest model trained with sub-trip data is the most accurate.

\begin{figure}[t]
\centering
\includegraphics[width=.5\textwidth]{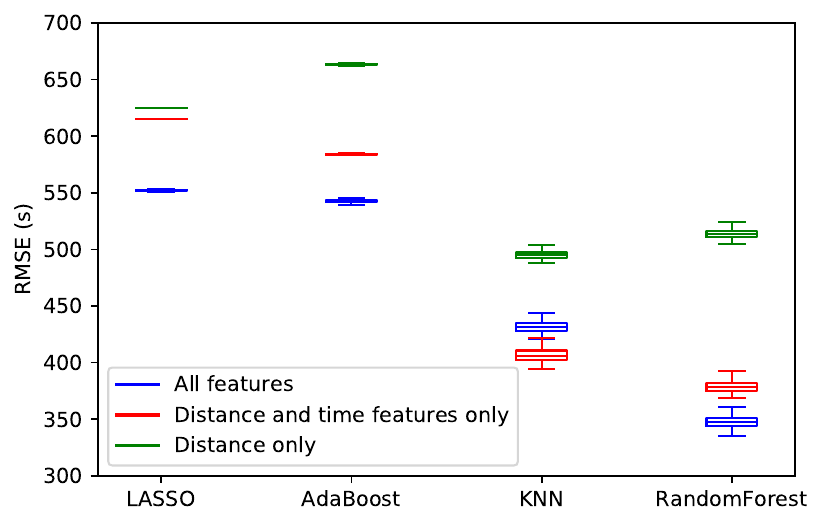}
\caption{Mean-squared error results for sub-trip models with only distance features, distance and time features, and all features.}\label{MEAN2d}
\end{figure}

To quantify the impact of time-based and geographical census features on prediction accuracy, we trained our sub-trip models with only distance features, distance and time features, and all features. Figure EC.\ref{MEAN2d} displays a boxplot comparing the RMSE of our models for these experiments. The Lasso, AdaBoost, and RandomForest models improved when both time and geographic features were included, and these improvements were found to be statistically significant using the Wilcoxon signed-rank test. In particular, including time-based features for the random forest model provides a RMSE improvement of 26\% (corresponding to a $134.5$s reduction in RMSE), over a model with access to only distance. The KNN model improved when adding time features, but it did not improve when geographical features were added because, unlike the other three models, KNN does not have an internal feature weighting process. In other words, the KNN model values all 73 census features equally. By using only the distance and time features, we are implicitly selecting the most important features for the model. These results reinforce the importance of considering time of day and day of the week for travel time estimation in urban areas in LMICs.

A random forest model comprising $1,000$ decision trees was selected as the final model and trained using all $4,086$ sub-trips. Each feature was available for inclusion to all 1,000 trees and relative feature importance was determined using the number of trees in the forest to which that feature contributes. Table~\ref{FeatureResultsTT} lists the features that had a relative importance greater than $0.01$. Our results suggest that travel distance, hour of day, and day of week are the three most important features. As expected, travel distance is the most dominant feature with a relative importance of $0.4128$. The hour of the day, which can be used as a proxy for peak traffic times, is the only other feature with an importance over $0.1$. Our findings are consistent with the results of previous traffic studies, which also found travel distance and the time of day to be the main factors \citep{Zhang2015, Vlahogianni2014}. As mentioned in Section~\ref{Lit:TTPred}, our approach extends previous work by incorporating demographic features for the origin and destination nodes. We found nine geographical census features with a relative importance of at least $0.01$. These additional features contribute to an 8\% reduction in RMSE relative to a random forest model that only has access to distance and time features.
\begin{table}
\caption{Relative feature importance as determined by the random forest model.\label{FeatureResultsTT}}
\centering
\begin{tabular}{l c r }
\noalign{\smallskip} \hline \noalign{\smallskip}
Feature && Relative importance \\ \hline \noalign{\smallskip}
Travel distance (m) && 0.413 \\
Hour of day && 0.147\\
Day of week && 0.084\\
Destination node medical facilities (no.) && 0.028 \\
Destination node ratio of male to female employment (\%) && 0.020\\
Destination node ratio of male to female industrial employment (\%) && 0.020\\
Destination node home owners (\%) && 0.020 \\
Destination node Jupri homes (\%) && 0.019 \\
Destination node population over 60 (\%) && 0.018\\
Origin node literacy rate (\%) && 0.015\\
Origin node home owners (\%) && 0.015 \\
Origin node non-sanitary toilets (\%) && 0.011 \\
\hline \noalign{\smallskip}
\end{tabular}
\end{table}

\section{Proof of Theorem 1.}\label{Apd:ProofThm1}

\proof{Proof.}

This proof establishes the Theorem 1 result by construction. Note that without vector notation, we can re-write NFF as:
\allowdisplaybreaks
\begin{equation} \label{NFF2}
\begin{aligned}
\mathrm{\underset{\by,\bff}{minimize}} \hspace*{1em} & \sum_{(i,j)\in \mathcal{E}} c_{ij}f_{ij} \\
\mbox{subject to}\hspace*{1em} & \sum_{i\in N} y_i= P, \\
& \sum_{j\in O(i)} f_{ij} - \sum_{j\in I(i)} f_{ji} \leq \alpha_i y_i - d_i, \forall i \in N, \\
& f_{ij}\geq 0, \forall (i,j)\in \mathcal{E}, \\
& y_i\in\{0,1\}, \forall i\in N,
\end{aligned}
\end{equation}
where $I(i) =\{ j\in N | (j,i)\in \mathcal{E}\}$ and $O(i)=\{ j\in N | (i,j)\in \mathcal{E}\}$. Recall the classic $p$-median formulation. The facility location variable is defined as $x_{ii}=1$ if a facility is located at node $i\in N$ and the assignment (routing) decision variables is denoted $x_{ij}=1$ if demand node $j$ has been assigned to facility node $i$. Using this notation, the $p$-median problem \citep{ReVelle1970} can be formulated as:
\allowdisplaybreaks
\begin{equation} \label{MED}
\begin{aligned}
\mathrm{\underset{\bx}{minimize}} \hspace*{1em} & \sum_{i\in N}\sum_{j\in N}d_jt_{ij}x_{ij} \\
\mbox{subject to}\hspace*{1em} &\sum_{i\in N} x_{ii} = P, \\
& \sum_{i\in N} x_{ij} = 1, \;\forall j \in N, \\
& x_{ii}\geq x_{ij}, \;\forall i, j, i\neq j, \in N, \\
& x_{ij}\geq 0, \forall i,j\in N, i\neq j, \\
& x_{ii}\in \{0,1\}, \forall i\in N,
\end{aligned}
\end{equation}
where  $t_{ij}$ denotes the shortest travel time between nodes $i,j\in N$ ($i$ and $j$ need not be adjacent) and $t_{ii}=0,\;\forall i\in N$. The optimal solution to \eqref{MED} is denoted by $\hat{\bx}$.

We first show that the $p$-median is polynomially reducible to NFF. That is, we show that the optimal solution from the $p$-median can be transformed into the optimal solution of NFF in polynomial time and both solutions have the same optimal cost. First, set $\tilde{y}_i=\hat{x}_{ii}, \forall i \in N$. By definition, the demand weighted shortest path length from $w\in N$ to $r\in N$ is given by $d_rt_{wr}$. To find the path from $w$ to $r$ (i.e., the sequence of nodes $i_w, ..., i_r$) along which flow must be directed, we solve:
\allowdisplaybreaks
\begin{equation} \label{SPP}
\begin{aligned}
\mathrm{\underset{\bff}{minimize}} \hspace*{1em} & \sum_{(i,j)\in E} c_{ij}f_{ij}^{rw} \\
\mbox{subject to}\hspace*{1em} & \sum_{i\in O(r)} f_{ri}^{rw} - \sum_{i\in I(r)} f_{ir}^{rw} = d_w,\\
& \sum_{j\in O(i)} f_{ij}^{rw} - \sum_{j\in I(i)} f_{ji}^{rw} = 0, \forall i \in N\setminus\{r,w\}, \\
& \sum_{i\in O(w)} f_{wi}^{rw} - \sum_{i\in I(w)} f_{iw}^{rw} = - d_w, \\
& f_{ij}\geq 0, \forall (i,j)\in E.
\end{aligned}
\end{equation}
We denote the optimal solution to (\ref{SPP}) as $\hat{\bff}^{rw}$. For the special case when $r=w$, Formulation (\ref{SPP}) is not well-defined and we assume that the shortest path has length zero (i.e., no flow is produced and $f_{ij}^{ww} = f_{ij}^{rr} = 0, \forall i,j \in N$). Set $\tilde{f}_{ij}=\sum\limits_{r\in N}\sum\limits_{w\in N} \hat{f}_{ij}^{rw} \hat{x}_{rw}$ to obtain a solution $(\tilde{\by}, \tilde{\bff})$ to Formulation~\eqref{NFF2}. We now show that the obtained solution $(\tilde{\by}, \tilde{\bff})$ is feasible with respect to \eqref{NFF2}. 

For the first constraint, we have: $$\sum\limits_{i\in N}\tilde{y}_i = \sum\limits_{i\in N} \hat{x}_{ii} = P.$$ 

For the second constraint, define $J=\{j\in N\:|\:y_j = 0\}$ and $I=\{i\in N\:|\: y_i=1\}$. Note that $I\cup J = N$. Consider some $k\in J$ (i.e., $y_k=0$),
\begin{align*}
\begin{split}
\sum_{j\in O(k)} \tilde{f}_{kj} - \sum_{j\in I(k)} \tilde{f}_{jk}
=& \sum_{j\in O(k)}\sum\limits_{r\in N}\sum\limits_{w\in N} \hat{f}_{kj}^{rw} \hat{x}_{rw} - \sum_{j\in I(k)} \sum\limits_{r\in N}\sum\limits_{w\in N} \hat{f}_{jk}^{rw} \hat{x}_{rw}, \\
=&\sum\limits_{r\in N}\sum\limits_{w\in N}\hat{x}_{rw}\left(\sum_{j\in O(k)} \hat{f}_{kj}^{rw} - \sum_{j\in I(k)} \hat{f}_{jk}^{rw}\right), \\
=&\sum\limits_{r\in N\setminus\{k\}}\hat{x}_{rk}\left(\sum_{j\in O(k)} \hat{f}_{kj}^{rk} - \sum_{j\in I(k)} \hat{f}_{jk}^{rk}\right) + \sum\limits_{w\in N\setminus\{k\}}\hat{x}_{kw}\left(\sum_{j\in O(k)} \hat{f}_{kj}^{kw} - \sum_{j\in I(k)} \hat{f}_{jk}^{kw}\right) \\ &+ \sum\limits_{r\in N\setminus\{k\}}\sum\limits_{w\in N\setminus\{k\}}\hat{x}_{rw}\left(\sum_{j\in O(k)} \hat{f}_{kj}^{rw} - \sum_{j\in I(k)} \hat{f}_{jk}^{rw}\right)+\hat{x}_{kk}\left(\sum_{j\in O(k)} \hat{f}_{kj}^{kk} - \sum_{j\in I(k)} \hat{f}_{jk}^{kk}\right), \\
=& \sum\limits_{r\in N\setminus\{k\}}\hat{x}_{rk}(-d_k) + \sum\limits_{w\in N\setminus\{k\}}\hat{x}_{kw} (d_w) + \sum\limits_{r\in N\setminus\{k\}}\sum\limits_{w\in N\setminus\{k\}}\hat{x}_{rw} (0) + \hat{x}_{kk}(0), \\
=&  -d_k.
\end{split}
\end{align*}
Consider some $k\in I$ (i.e., $y_k=1$),
\begin{align*}
\begin{split}
\sum_{j\in O(k)} \tilde{f}_{kj} - \sum_{j\in I(k)} \tilde{f}_{jk}
=& \sum_{j\in O(k)}\sum\limits_{r\in N}\sum\limits_{w\in N} \hat{f}_{kj}^{rw} \hat{x}_{rw} - \sum_{j\in I(k)} \sum\limits_{r\in N}\sum\limits_{w\in N} \hat{f}_{jk}^{rw} \hat{x}_{rw}, \\
=&\sum\limits_{r\in N}\sum\limits_{w\in N}\hat{x}_{rw}\left(\sum_{j\in O(k)} \hat{f}_{kj}^{rw} - \sum_{j\in I(k)} \hat{f}_{jk}^{rw}\right), \\
=&\sum\limits_{r\in N\setminus\{k\}}\hat{x}_{rk}\left(\sum_{j\in O(k)} \hat{f}_{kj}^{rk} - \sum_{j\in I(k)} \hat{f}_{jk}^{rk}\right) + \sum\limits_{w\in N\setminus\{k\}}\hat{x}_{kw}\left(\sum_{j\in O(k)} \hat{f}_{kj}^{kw} - \sum_{j\in I(k)} \hat{f}_{jk}^{kw}\right) \\ &+ \sum\limits_{r\in N\setminus\{k\}}\sum\limits_{w\in N\setminus\{k\}}\hat{x}_{rw}\left(\sum_{j\in O(k)} \hat{f}_{kj}^{rw} - \sum_{j\in I(k)} \hat{f}_{jk}^{rw}\right)+\hat{x}_{kk}\left(\sum_{j\in O(k)} \hat{f}_{kj}^{kk} - \sum_{j\in I(k)} \hat{f}_{jk}^{kk}\right), \\
=& \sum\limits_{r\in N\setminus\{k\}}\hat{x}_{rk}(-d_k) + \sum\limits_{w\in N\setminus\{k\}}\hat{x}_{kw} (d_w) + \sum\limits_{r\in N\setminus\{k\}}\sum\limits_{w\in N\setminus\{k\}}\hat{x}_{rw} (0) + \hat{x}_{kk}(0), \\
\leq & \sum\limits_{w\in N\setminus\{k\}} d_w = \sum_{w\in N} d_w - d_k = \alpha - d_k.
\end{split}
\end{align*}

Lastly, we show that the objective function values of both solutions are equal,
\begin{align*}
\begin{split}
\sum_{(i,j)\in E} c_{ij}\tilde{f}_{ij}
&= \sum_{(i,j)\in E} c_{ij}\sum_{r\in N}\sum_{w\in N}\hat{f}_{ij}^{rw} \hat{x}_{rw}, \\
&= \sum_{r\in N}\sum_{w\in N} \hat{x}_{rw}\left(\sum_{(i,j)\in E} c_{ij}\hat{f}_{ij}^{rw}\right), \\
& =\sum_{r\in N}\sum_{w\in N} \hat{x}_{rw} d_w t_{rw}.
\end{split}
\end{align*}

We now prove the reverse direction. That is, a solution from NFF can be transformed into a solution for the $p$-median with the same optimal cost. We denote the optimal solution to NFF as $(\hat{\by},\hat{\bff})$.

First, we set $\tilde{x}_{kk} = \hat{y}_k, \forall k\in N$. Define $J=\{r\in N\:|\:\hat{y}_r = 0\}$ and $I=\{w\in N\:|\: \hat{y}_w=1\}$. Note that $I\cup J = N$. Compute $t_{rw},$ $\forall r\in I$ and $\forall w\in J$ (i.e., the shortest path between nodes $r$ and $w$). This can be done by using Dijkstra's algorithm or by extracting the path lengths directly from the given optimal solution to NFF. Both methods are polynomial time.

Now, consider some $k\in J$, and solve $\mathrm{argmin_{i\in I}} \;t_{ik}$. Denote the optimal index as $i^k$ and the optimal value as $t_{i^kk}$. Set $\tilde{x}_{i^k k} =1$, $\tilde{x}_{k j} =0, \forall j\in N$, and $\tilde{x}_{ik}=0, \forall i\in N\setminus\{i^k\}$. Consider some $k\in I$, which implies that $\tilde{x}_{kk}=1$. Set $\tilde{x}_{ik}=0,\forall i\in N\setminus\{k\}$ to obtain the solution, $\tilde{\bx}$. We now show that the obtained solution $\tilde{\bx}$ is feasible for the $p$-median. 

For the first constraint, we have: $$\sum\limits_{i\in N}\tilde{x}_{ii} = \sum\limits_{i\in N} \hat{y}_{i} = P.$$

For the second constraint, consider $r\in J$. By our construction, $\tilde{x}_{i^r r} =1$ and $\tilde{x}_{ir}=0, \forall i\in N\setminus\{i^r\}$. Therefore, $\sum\limits_{i\in N} \tilde{x}_{ir} = \tilde{x}_{i^r r} +\sum\limits_{i\in N\setminus\{i^r\}}\tilde{x}_{ir}  = 1$. Consider, $w\in I$. By our construction, $\tilde{x}_{ww}=1$ and $\tilde{x}_{iw}=0 , \forall i\in N\setminus\{w\}$. Therefore, $\sum\limits_{i\in N} \tilde{x}_{iw} = \tilde{x}_{ww} +\sum\limits_{i\in N\setminus\{w\}}\tilde{x}_{iw}  = 1$. Combining these implies $\sum_{i\in N} \tilde{x}_{ij} = 1, \forall j\in N$.

For the third constraint, consider $r\in J$. By our construction $\tilde{x}_{rr} =0$ and $\tilde{x}_{rk}=0, \forall k\in N\setminus\{r\}$. Therefore, $\tilde{x}_{rr}\geq\tilde{x}_{rk}, \forall r\in J, k\in N\setminus\{r\}$. Consider $w\in I$. By our construction, $\tilde{x}_{ww}=1$ and $\tilde{\bx}\in\{0,1\}$ (i.e., $\tilde{\bx}\leq 1$), therefore we have $\tilde{x}_{ww}\geq\tilde{x}_{wk}, \forall w\in I, k\in N\setminus\{w\}$. Combining these implies $\tilde{x}_{ii}\geq\tilde{x}_{ij}, \forall i\in N, j\in N\setminus\{i\}$.

Lastly, we show that the objective function values are equal. First we must derive some intermediate information. Consider the following optimization problem with $\hat{\by}$ fixed,
\begin{equation} \label{SPP1}
\begin{aligned}
\mathrm{\underset{\bff}{minimize}} \hspace*{1em} & \sum\limits_{(i,j)\in E} f_{ij}c_{ij} \\
\mbox{subject to}\hspace*{1em} & \sum_{j\in O(i)}f_{ij} - \sum_{j\in I(i)}f_{ji} \leq \alpha \hat{y}_i - d_i, \forall i\in N,\\
& f_{ij}\geq 0, \forall (i,j)\in E. \\
\end{aligned}
\end{equation}
Denote the optimal solution of (\ref{SPP1}) by $\hat{\bff}$. The dual of (\ref{SPP1}) is given by,
\begin{equation} \label{SPP2}
\begin{aligned}
\mathrm{\underset{\bp}{maximize}} \hspace*{1em} & \sum_{i\in N}p_i(\alpha\hat{y}_i - d_i) \\
\mbox{subject to}\hspace*{1em} & p_i - p_j \leq c_{ij}, \forall (i, j)\in E, \\
& p_i\leq 0, \forall i \in N. 
\end{aligned}
\end{equation}
Denote the optimal solution to (\ref{SPP2}) by $\hat{\bp}$. The dual variable $\hat{p}_k$ represents the change in optimal cost due to increasing $d_k$ by one unit. If we increase $d_k$ by one unit, the optimal solution will increase by the length of the shortest path from $i^k\in I$ to $k$. Therefore, at optimality, the value of -$p_k$ (because $p_k$ is negative in \ref{SPP2}) is equal to the length of the shortest path from $i^k\in I$ to $k\in J$. Mathematically, $-p_k = t_{i^k, k}$. Note that this implies that $p_k=0, \forall k\in I$ because the shortest path from a facility to itself, has length zero.

Now we show that optimal costs are equal:
\begin{align*}
\begin{split}
\sum_{w\in N}\sum_{r\in N} \tilde{x}_{rw} d_w t_{rw} &=\sum_{r\in J}\sum_{w\in N} \tilde{x}_{rw} d_w t_{rw} +\sum_{r\in I}\sum_{w\in J} \tilde{x}_{rw} d_w t_{rw} +  \sum_{r\in I}\sum_{w\in I\setminus\{r\}} \tilde{x}_{rw} d_w t_{rw} +\tilde{x}_{ww}d_wt_{ww} \\
&= \sum_{w\in J}\sum_{r\in I} \tilde{x}_{rw} d_w t_{rw} \hspace*{5.4em} (t_{ww}=0, \;\tilde{x}_{rw} = 0, \forall r\in J, w\in N, \\ &\quad\quad\quad\quad\quad\quad\quad\quad\quad\quad\quad\quad\quad\quad\mbox{and} \,x_{wr}=0, \forall w\in I, r\in I\setminus\{w\}) \\
&= \sum_{w\in J} d_w\sum_{r\in I} \tilde{x}_{rw} t_{rw} \\
& = \sum_{w\in J} d_w t_{r^ww} \hspace*{8em}  \mbox{(By our construction)} \\
& = - \sum_{w\in J} d_w \hat{p}_w \hspace*{8em}  \mbox{(From duality)}\\
& = - \sum_{w\in N} d_w \hat{p}_w \hspace*{8em}  (\hat{p}_w=0, \forall w\in I)\\
& = \alpha\sum_{w\in N} y_w \hat{p}_w- \sum_{w\in N} d_w \hat{p}_w \hspace*{3em}  (\hat{p}_w=0, \forall w\in I \mbox{  and  } y_w=0, \forall w\in J)\\
& = \sum_{w\in N} \hat{p}_w(\alpha y_w - d_w) \\
& = \sum_{(r,w)\in E} \hat{f}_{rw}c_{rw} \hspace*{7.3em}  \mbox{(By strong duality)}
\end{split}
\end{align*}
Combining both directions, we have that $$\sum_{(i,j)\in E}\bar{f}_{ij}c_{ij}=\sum_{i\in N}\sum_{j\in N} \hat{x}_{ij} d_j t_{ij}\leq\sum_{i\in N}\sum_{j\in N} \tilde{x}_{ij} d_j t_{ij},\forall \tilde{\bx},$$ and that $$\sum_{i\in N}\sum_{j\in N} \tilde{x}_{ij} d_j t_{ij}=\sum_{(i,j)\in E}\hat{f}_{ij}c_{ij}\leq\sum_{(i,j)\in E}\bar{f}_{ij}c_{ij},\forall \bar{\bff}.$$Therefore,$$\sum_{i\in N}\sum_{j\in N} \hat{x}_{ij} d_j t_{ij}=\sum_{(i,j)\in E}\hat{f}_{ij}c_{ij}.\Halmos$$
\endproof

\section{Proof of Theorem 2.}\label{Apd:ProofThm2}
\proof{Proof.}
Let $\mathbb{Y}=\{\by\; | \;\beee'\by = P, \by\geq \bzero\}$ and $\mathbb{F}(\by,\bd) = \{\bff \;|\;\bA\bff\leq\balpha\bI\by-\bd,\bff\geq\bzero\}$. Then, \textbf{R-NFF} can be written as
\allowdisplaybreaks
\begin{equation*}
\begin{aligned}
\mathrm{\min_{\by\in\mathbb{Y}}\;\max_{\bc\in\mathcal{C}, \bd\in\mathcal{D}}\;\min_{\bff\in\mathbb{F}(\by,\bd)}} \hspace*{1em} & \bc'\bff, \\
\end{aligned}
\end{equation*}
or in epigraph form
\begin{equation*}
\begin{aligned}
\mathrm{\underset{\by\in\mathbb{Y},t}{minimize}}\quad&t \\
\mbox{subject to}\quad & t\geq\max_{\bc\in\mathcal{C}, \bd\in\mathcal{D}}\;\min_{\bff\in\mathbb{F}(\by,\bd)} \hspace*{.5em}\bc'\bff. \\
\end{aligned}
\end{equation*}
Enumerating the elements of $\mathcal{D}$, the model becomes
\begin{equation*}
\begin{aligned}
\mathrm{\underset{\by\in\mathbb{Y},t}{minmize}}\quad&t \\
\mbox{subject to}\quad & t\geq\max_{\bc\in\mathcal{C}}\;\min_{\bff\in\mathbb{F}(\by,\bd^k)} \hspace*{.5em}\bc'\bff,\quad k = 1,...,N. \\
\end{aligned}
\end{equation*}
Since $\mathcal{C}$ and $\mathbb{F}(\by,\bd^k)$ are disjoint, we can swap the min and max using the min-max theorem \citep{VonNeumann1928}:
\begin{equation*}
\begin{aligned}
\mathrm{\underset{\by\in\mathbb{Y},t}{minimize}}\quad&t \\
\mbox{subject to}\quad & t\geq\min_{\bff\in\mathbb{F}(\by,\bd^k)}\;\max_{\bc\in\mathcal{C}} \hspace*{.5em}\bc'\bff,\quad k = 1,...,N.\\
\end{aligned}
\end{equation*}
We then replace $\bff$ by $\bff^k$ for each scenario $k$, which yields: 
\begin{equation*}
\begin{aligned}
\mathrm{\underset{\by\in\mathbb{Y},t}{minimize}}\quad&t \\
\mbox{subject to}\quad & t\geq\min_{\bff^k\in\mathbb{F}(\by,\bd^k)}\;\max_{\bc\in\mathcal{C}} \hspace*{.5em}\bc'\bff^k,\quad k = 1,...,N.\\
\end{aligned}
\end{equation*}
We can now move $\bff^k$ to the outer minimization problem:
\begin{equation*}
\begin{aligned}
\mathrm{\underset{\by\in\mathbb{Y},t, \bff^k}{minimize}}\quad&t \\
\mbox{subject to}\quad & t\geq\mathrm{\underset{\bc\in\mathcal{C}}{maximize}} \hspace*{.5em}\bc'\bff^k,\quad k = 1,...,N,\\
& \bff^k\in\mathbb{F}(\by,\bd^k), \quad k = 1,...,N.\\
\end{aligned}
\end{equation*}
Finally, for each $k$, we take the dual of the inner maximization problem to obtain the required result. \Halmos\endproof


\section{Comparison of solution approaches}\label{CompExp}

In this section, we present results from a set of computational experiments that compare the effectiveness of our exact and heuristic scenario generation algorithms. To do so, we use smaller randomly-generated problem instances that can be solved to optimality.

\subsection{Experimental setup.}\label{CompExp:Setup}

We use three random network instances to conduct our experiments. The first network has 30 nodes and 90 edges, the second has 50 nodes and 150 edges, and the third network has 75 nodes and 226 edges. For each graph, we vary the number of scenarios ($|\mathcal{S}|$) in $\mathcal{D}$, the interdiction budget ($B$) in $\mathcal{C}$, and the number of vehicle outposts ($P$). Specifically, we consider: $|\mathcal{S}|\in\{1, 10, 100, 1000, 10000\}$, $P\in\{1,2,5,10,25\}$, and $B\in\{0,10,50,100,250,500,1000\}$. Hence, we solve 175 problem instances for each random network, for a total of 525 problem instances. We solve each instance using 1) a commercial solver (Gurobi), 2) our exact scenario generation algorithm (SGen), and 3) our heuristic scenario generation algorithm (HSGen) with $10$ random starts and $10$ interchanges. We chose this number of random starts and interchanges after testing our heuristic with different values (see Figure~\ref{RandomStarts}). We set a maximum time limit of 36,000 seconds for each instance. All experiments were programmed using \texttt{MATLAB2016a} and run on a desktop computer with an Intel Core i7-4790K 4.0 GHz processor and 32 GB of RAM.%

To estimate realistic edge-lengths for these instances, we randomly sample from the edge data distribution introduced in Figure~\ref{TVD}. To generate node-weights and demand scenarios, we used a modified version of the methodology outlined in Section~\ref{Opt:DUnc}. We estimate the population and $\xi$ as in Section~\ref{Opt:DUnc}. Without underlying ward features, our logistic regression model is not applicable so we use the naive approach from~\ref{Apd:DemandEst3} instead.

\subsection{Scenario generation algorithm performance.}\label{CompExp1}
The scenario generation algorithm was able to solve all 525 problem instances to optimality, while Gurobi struggled with larger instances. Table~\ref{HSGen1} compares solution times as a function of uncertainty set size. Gurobi was not able to solve any of the instances that had 100 or more demand scenarios, except the one with no travel time uncertainty. Table~\ref{HSGen2} compares the solution times for instances that vary in the size of the underlying network and the number of outposts located. The scenario generation approach enjoys the largest speed up for intermediate values of $P$.

\subsection{Heuristic algorithm performance.}\label{CompExp2}
Table~\ref{HSGen1} also compares the optimal cost and solution time between SGen and HSGen as a function of the number of scenarios and the interdiction budget. The objective function value is displayed as mean response time, in seconds. To determine this value, we divide the actual objective function value by the total number of trips. The performance of the heuristic algorithm remains relatively stable as the number of scenarios increases with solutions times that are an order or magnitude less than SGen.  Table~\ref{HSGen2} compares the optimal cost and solution time between SGen and HSGen as a function of the size of the graph and the number of outposts, while holding both the interdiction budget and the number of scenarios constant at 100. Across all instances, the heuristic algorithm was able to obtain the optimal solution when the number of outposts was small. The performance also remains relatively stable as the size of the graph grows. However, the performance degrades as $P$ grows. This degradation in performance is balanced by up to an order of magnitude speed-up in certain cases. While HSGen does not close the optimality gap as the size of the problem increases, for the large-scale, real-world instances of the robust problems that we solve in Section~\ref{DhakaExp}, it is the only method capable of generating a solution in a reasonable time limit. 

\begin{figure}[t]
\centering
\subfigure[\ Number of scenarios \label{RS1}]{
\includegraphics[width=.5\textwidth]{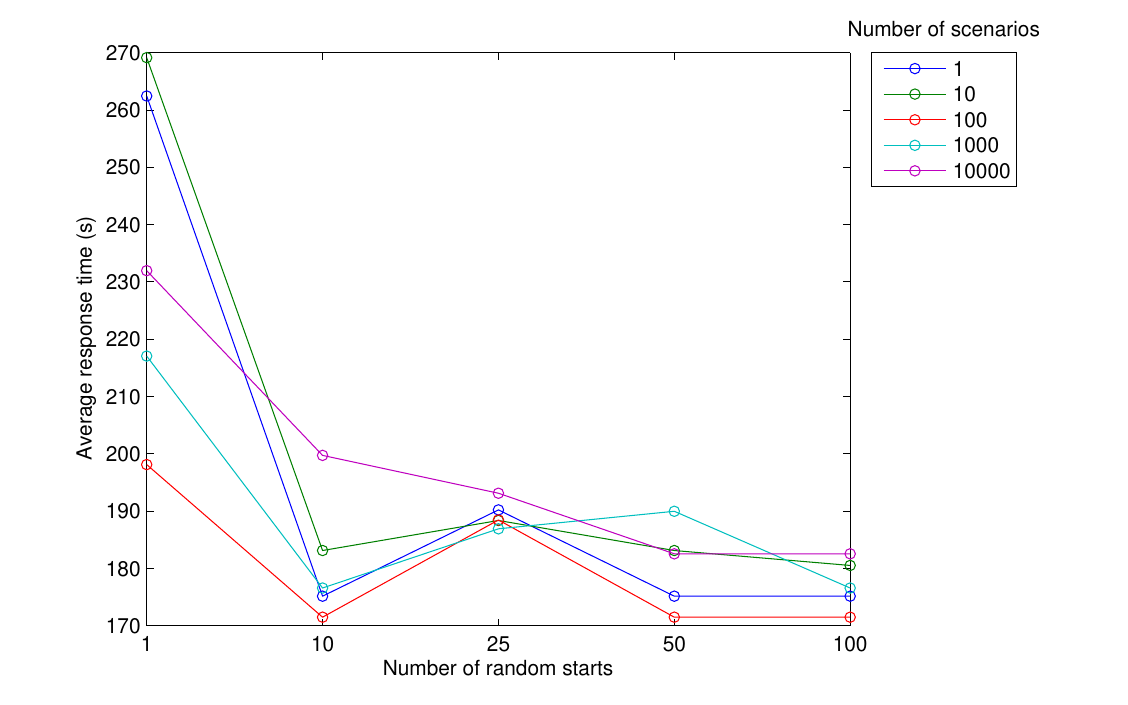}}
\subfigure[\ Number of outposts \label{RS2}]{
\includegraphics[width=.46\textwidth]{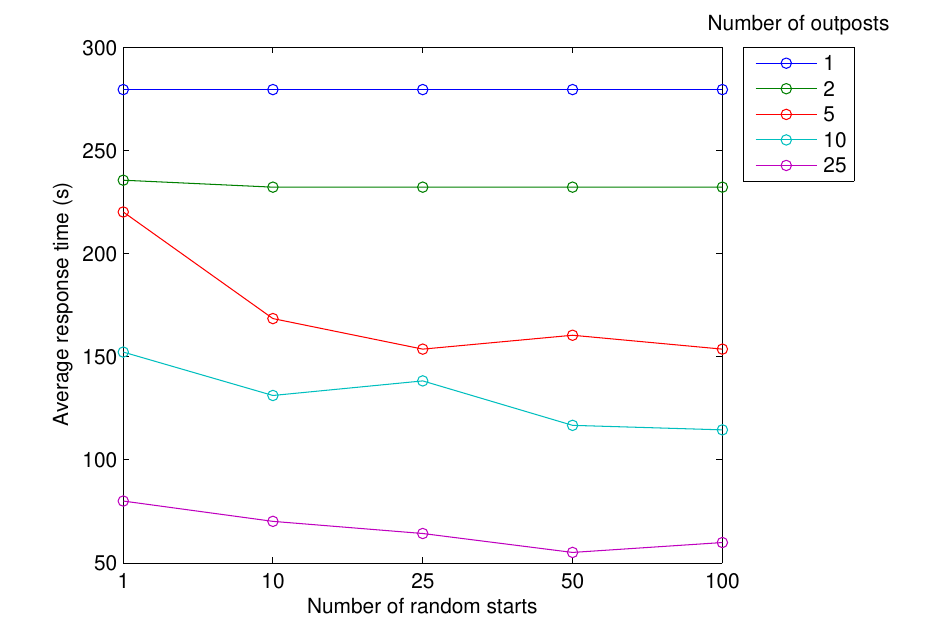}}
\subfigure[\ Number of scenarios \label{RS3}]{
\includegraphics[width=.5\textwidth]{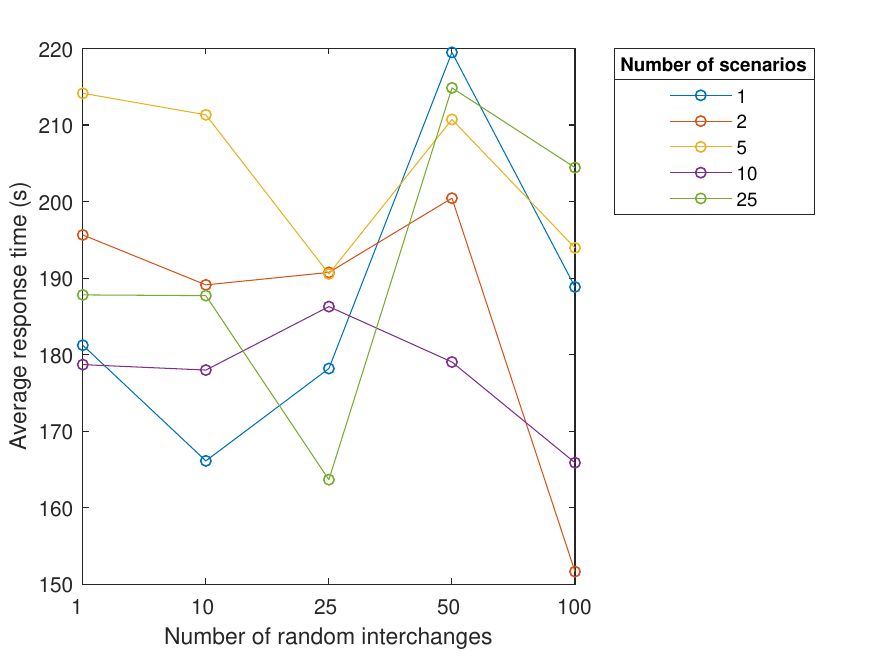}}
\subfigure[\ Number of outposts \label{RS4}]{
\includegraphics[width=.46\textwidth]{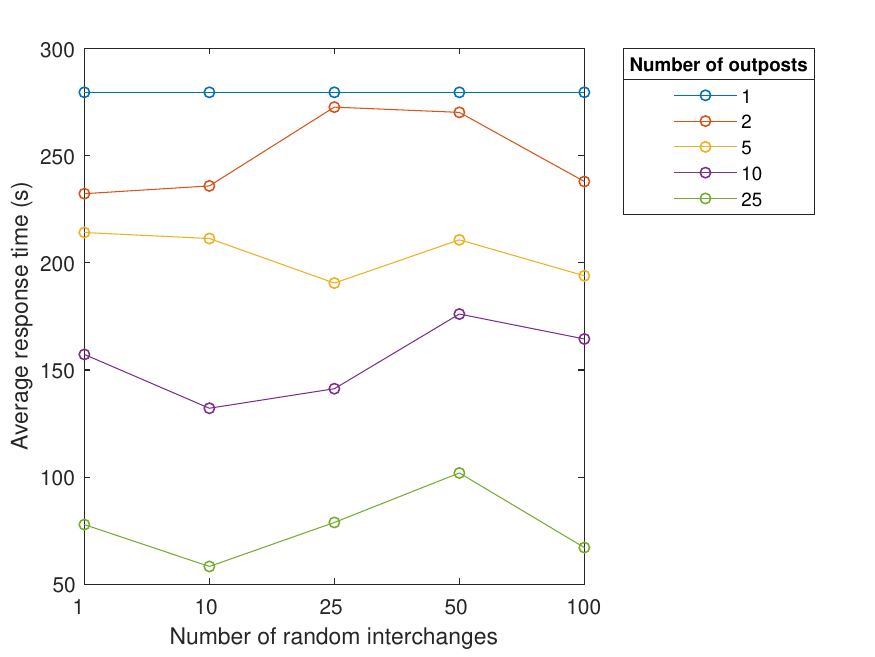}}
\caption{Objective function value as a function of the number of random starts and interchanges used for HSGen.\label{RandomStarts}}
\end{figure}

Figure~\ref{RandomStarts} displays the performance of HSGen as a function of the number of random starts and random interchanges for different numbers of scenarios and different numbers of outposts. Figure~\ref{RS1} shows that HSGen improves significantly from one to ten random starts, but does not appear to improve much beyond ten. By contrast, Figure~\ref{RS2} displays a small improvement from one to ten random starts and only marginal improvements thereafter. Thus, we use ten random starts to conduct our real experiments on the Dhaka road network. Figures~\ref{RS3} and~\ref{RS4} show that there does not appear to be a correlation between increasing the number of random interchanges and the overall solution quality. However, in all cases, there is a small improvement from one to 10 random interchanges. Thus, we use we use ten random interchanges to conduct our real experiments on the Dhaka road network.

\begin{table}
\caption{Comparison of objective function values and solution times between Gurobi, SGen, and HSGen as a function of the number of scenarios and uncertainty budget. The number of vehicles and the size of the graph are held constant at 5 and, 75 nodes and 226 edges, respectively. \label{HSGen1}}
\centering
\begin{tabular}{l c r r r r r r}
&& \multicolumn{3}{c}{\underline{     Objective function value     }} & \multicolumn{3}{c}{\underline{      Solution time      }} \\
$|S|$ & Budget & SGen & HSGen & Optimality gap (\%) & Gurobi (s) & SGen (s) & HSGen (s) \\ \hline
\multirow{4}{*}{1} &  1 &  133.4    &    142.8    &     6.6    & 3.4&    3.4       & 3.7 \\
& 10 & 134.8    &      141.0     &    4.4    &   4.3 &  4.3    &    3.6  \\
& 100 &  145.9    &   171.7     &      15.0    &  5.6 &   5.7   &     3.7   \\
& 1000 & 205.1    &   214.7   &      4.4    &   10.7 & 11.1   &     4.3   \\ \hline
\multirow{4}{*}{10} &  1 &  135.1    &   156.7    &    13.7    & 179.5 &  37.3   &     9.4\\
& 10 &  136.6    &   147.8      &   7.5    & 192.5 &  31.2   &     9.2  \\
& 100 &   148.4    &   169.0     &    12.2    &  625.3 & 109.7   &    10.7    \\
& 1000 & 208.3    &    209.8     &    0.7    &  4581.7 & 72.0   &    12.1   \\ \hline

\multirow{4}{*}{100} &  1  & 143.2     &  164.9    &    13.1    &  26857.0 &  137.9   &    19.1    \\
& 10 &  144.5     &  163.1    &    11.4    &   - & 170.1   &    24.1  \\
& 100 & 153.7     &  167.1     &    8.0    &  - &  197.5   &    16.9    \\
& 1000 & 218.4     &  227.2     &    3.9    &  - &  2974.8   &    24.4  \\  \hline
\multirow{4}{*}{1000} &  1 &  140.1    &   153.8      &   8.9    &  - & 585.6  &      84.2      \\
& 10 & 141.4 & 161.8 & 12.6 & - &945.3 & 63.0 \\
& 100 & 152.5 & 158.4 & 3.8 & - &1160.0 & 93.7 \\
& 1000 & 212.1  &     219.0   &     3.1 &   - &       3721.8 &       34.8    \\ \hline
\multirow{4}{*}{10000} &  1 & 142.5   &    162.4 &      12.2   &   - &  3856.3 &       581.2  \\
& 10 &  143.8   &    161.6     &  11.0  &  - &      3264.3 &       462.7  \\
& 100 & 154.7   &    164.0   &    5.7 &    - &     9306.5 &       581.5 \\
& 1000 & 217.8    &   227.0   &    4.0 &     - &    5475.0&       350.7 \\ \hline
\end{tabular}
\end{table}

\setlength{\extrarowheight}{2pt}
\begin{table}
\caption{Comparison of objective function values and solution times between SGen and HSGen as a function of graph size and $P$. The number of scenarios and the interdiction budget are held constant at 100 and 100, respectively. \label{HSGen2}}
\centering
\begin{tabular}{c c c r r r r r r}
&&& \multicolumn{3}{c|}{\underline{     Objective function value     }} & \multicolumn{3}{c}{\underline{      Solution time      }} \\
Nodes & Edges & P & SGen & HSGen & Optimality gap (\%) & Gurobi (s) & SGen (s) & HSGen (s) \\ \hline
\multirow{5}{*}{30} & \multirow{5}{*}{90} & 1 & 213.9   &    213.9   &         0   &    87.8 & 12.7    &    2.0    \\
& & 2 &  155.8  &     155.8  &          0  &   82.2 &   25.5   &     9.2   \\
& & 5 & 93.0  &     102.8  &       9.6   &  473.9 &  25.1   &     7.7 \\
& & 10 &  46.1  &      60.3  &      23.5   &  210.4 &  94.8   &    42.9   \\
&& 25 &     7.3  &      11.3  &      35.4   &  10.6 &   142.1  &    40.1   \\ \hline
\multirow{5}{*}{50} & \multirow{5}{*}{150} & 1 &       287.0  &     287.0  &          0   &  525.1 &  38.5   &     9.6    \\
& & 2 &        191.0  &     191.0 &           0   &  1213.5 &  25.9   &     5.9  \\
& & 5 &       109.9  &     130.0 &       15.5   &  3534.0 &  89.6   &    25.4    \\
& & 10 & 62.8  &      88.7 &       29.2   &  3410.6 &  100.0   &    86.7    \\
&& 25 &   17.2  &      31.8 &       45.7   &   211.6 & 197.4   &    99.4    \\ \hline
\multirow{5}{*}{75} & \multirow{5}{*}{226} & 1 &  279.5  &     279.5 &           0   &    2671.8 & 31.7   &     6.9    \\
& & 2 & 232.2  &     232.2 &           0   &   - &  88.8  &     13.1    \\
& & 5 &   153.7  &     167.1 &        8.0   &  - &  197.5  &     16.9  \\
& & 10 &  99.4  &     137.3 &       27.6   &   - &  400.9   &    32.0    \\
&& 25 &     40.4  &      65.9 &       38.7   &  2591.1 &  176.9   &    113.9   \\ \hline
\end{tabular}
\end{table}

\section{Tactical Simulation Model}\label{Apd:Sim}

Algorithm~\ref{SimAlgPsuedo} displays high-level pseudo-code for our simulation framework. Note that $t_C$ denotes the current time, $t_W$ denotes the waiting time, $t_D$ denotes the drive time to the emergency location, and $t_S$ denotes the scene time sampled from an exponential distribution with mean of 15 min. The DISPATCH function uses a greedy dispatching policy that assigns the closest ambulance and we determine the closest ambulance by solving the robust shortest path problem with $B=1000$ for each available ambulance. After reaching the scene and picking up the patient, ambulances are routed to the closest hospital determined by solving the robust shortest path problem with $B=1000$ and $t_H$ denotes the drive time to the hospital. The ROUTEHOME function  determines the time until an ambulance has returned to its home base location ($t_B$) by solving the robust shortest path problem with $B=100$.

\begin{algorithm}
\caption{Tactical ambulance simulator}\label{SimAlgPsuedo}
\begin{algorithmic}[1]
\Function{SIMULATE(C,R,Y)}{}
\State $C \gets$ Load the simulated call times and locations 
\State $R \gets$ Load the road network with hospital locations and travel times for each day/hour 
\State $Y \gets$ Load the optimized outpost locations and the number of ambulances per outpost
\State $Q \gets \emptyset$ \Comment{initialize empty call queue}
\State $\mathcal{E} \gets C $ \Comment{initialize event queue with calls} 
\State $A \gets Y$ \Comment{initialize ambulance availability list} 
\While{$|\mathcal{E}|>0$}
\State Remove next event $e$ from $\mathcal{E}$
\State Update current time $t_C$
\If{$e$ is a new call}
\If{$|A| = 0$} \Comment{No available ambulances}
\State $Q \gets Q + e$ \Comment{Queue the call}
\Else 
\State $V(c) \gets$ DISPATCH(c,A,R) \Comment{Dispatch closest ambulance}
\State $A = A - V(c)$ \Comment{Remove dispatched vehicle from available list}
\State $e_{new}$ \Comment{Create new event for when ambulance is free at hospital}
\State $\mathcal{E} \gets e_{new}$ \Comment{Insert new event at time $t = t_C + t_W + t_D + t_S + t_H$}
\EndIf
\ElsIf{$e$ is ambulance becomes available}
\If{$|Q| > 0$}
\State $V(c) \gets$ DISPATCH(c,A,R) \Comment{Dispatch newly available ambulance}
\State $e_{new}$ \Comment{Create new event for when ambulance is free at hospital}
\State $\mathcal{E} \gets e_{new}$ \Comment{Insert new event at time $t = t_C + t_W + t_D + t_S + t_H$}
\Else
\State ROUTEHOME(c,Y,R) \Comment{Route ambulance home}
\State $e_{new}$ \Comment{Create new event for when ambulance is free at its base}
\State $\mathcal{E} \gets e_{new}$ \Comment{Insert new event at time $t = t_C + t_B$}
\EndIf
\ElsIf{$e$ is ambulance returned to base}
\If{$|Q| > 0$}
\State $V(c) \gets$ DISPATCH(c,A,R) \Comment{Dispatch newly available ambulance}
\State $e_{new}$ \Comment{Create new event for when ambulance is free at hospital}
\State $\mathcal{E} \gets e_{new}$ \Comment{Insert new event at time $t = t_C + t_W + t_D + t_S + t_H$}
\Else 
\State $A = A + V(c)$ \Comment{add ambulance to available list}
\EndIf
\EndIf
\EndWhile
\EndFunction
\end{algorithmic}
\end{algorithm}


\section{Dhaka Policy Experiments}\label{Apd:Results}

\begin{figure}[t]
\centering
\includegraphics[width=.25\textwidth]{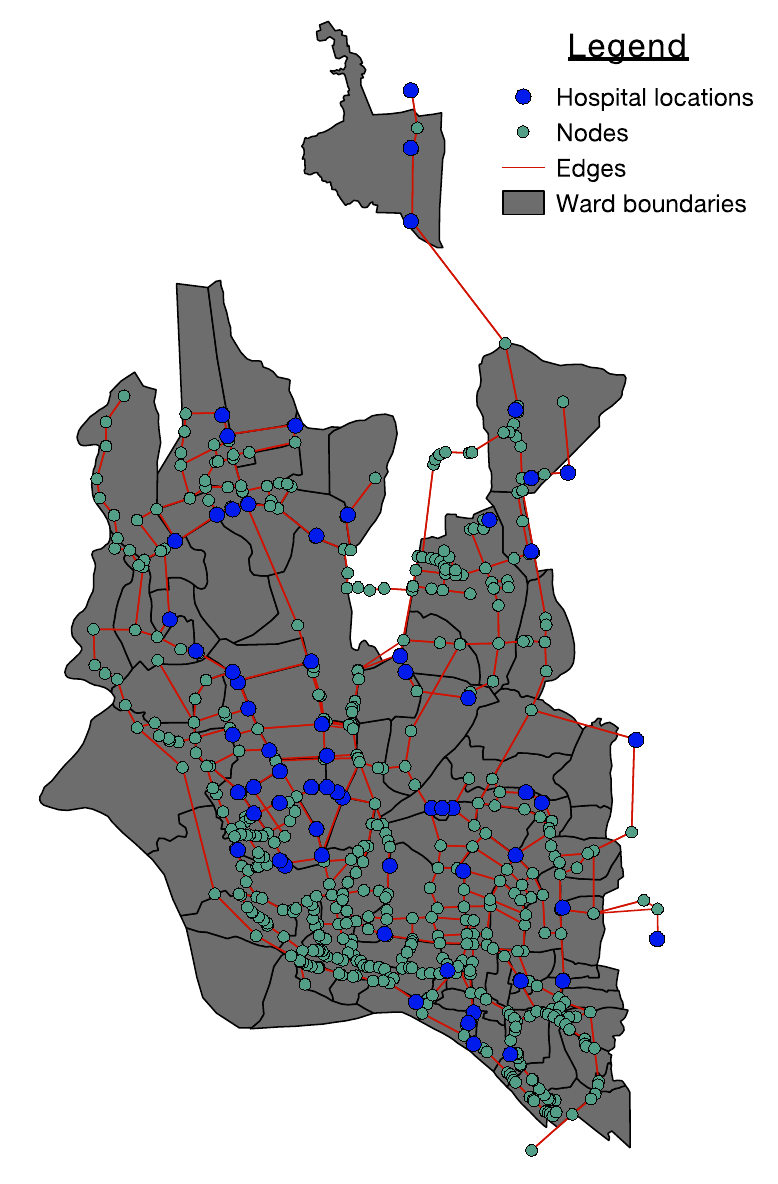}
\caption{Map of all 67 current outpost locations. \label{CurrentLocations}}
\end{figure}

Figure~\ref{CurrentLocations} depicts the locations of all 67 current outposts.

\subsection{What is the impact of the number of ambulances per outpost?}\label{Apd:NumAper}

\begin{figure}[t]
\centering
\includegraphics[width=.5\textwidth]{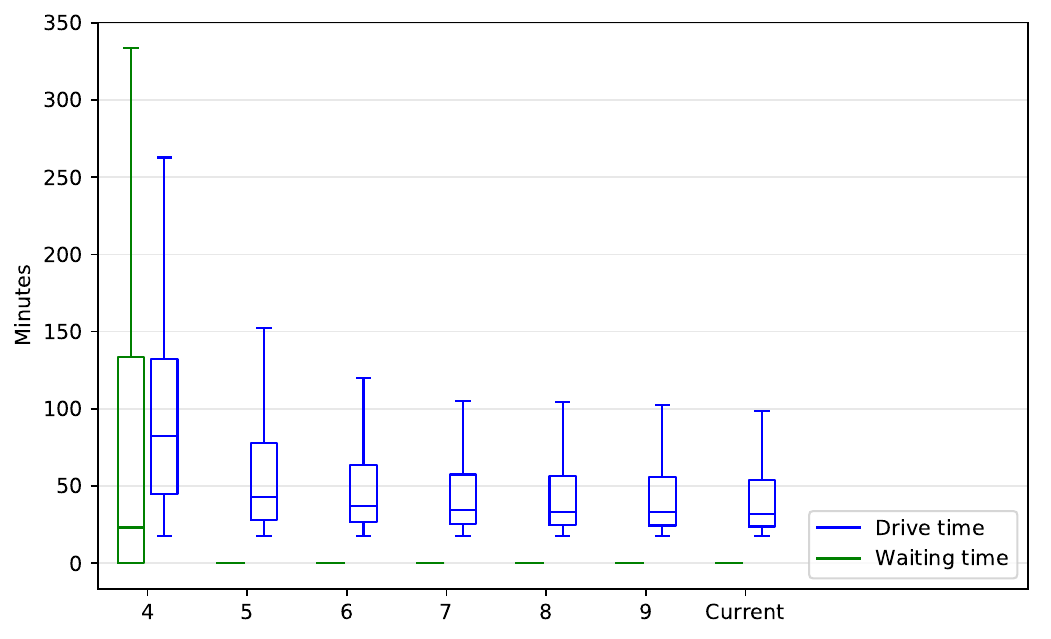}
\caption{Drive time and waiting time as a function of the number of ambulances per outpost for 20 outpost locations. ``Current" is the baseline scenario with 67 hospital-based outposts. \label{WTRTnumA}}
\end{figure}

Figure~\ref{WTRTnumA} displays the two major components of response time (drive time and waiting time) as a function of the number of ambulances per outpost for 20 outpost locations and the current baseline scenario (67 hospital-based outposts). The current baseline scenario includes a total of 269 ambulances spread across 67 outposts, while the 20 outpost solution with nine ambulances per outpost includes only 180 total ambulances. In other words, similar response time performance can be achieved with far fewer resources, if the resources are utilized more effectively. We find that diminishing returns are reached with seven ambulances per outpost (a total of 140) and we use seven per outpost for our policy experiments. Note that three or fewer ambulances per outpost result in a system with waiting times over 24 hours.

\subsection{What is the impact of the travel time budget?}\label{Apd:BudgetSens}

In this section, we conduct a sensitivity analysis to determine the impact of the travel time budget. We use the HSGen algorithm with 10 random starts and 10 random interchanges to solve \eqref{NFFsolve} with $P=20$ and $B=\{0,100,1000,2500,5000,7500,10000\}$. We then apply the outpost locations resulting from a budget of 1000 seconds to all seven budget instances. We compare these results with the response time of outpost locations optimized and evaluated on the same budget. We conduct separate experiments for each of the three temporal combinations.

\begin{figure}[t]
\centering
\subfigure[\ Rush hour \label{BS0}]{
\includegraphics[width=.3\textwidth]{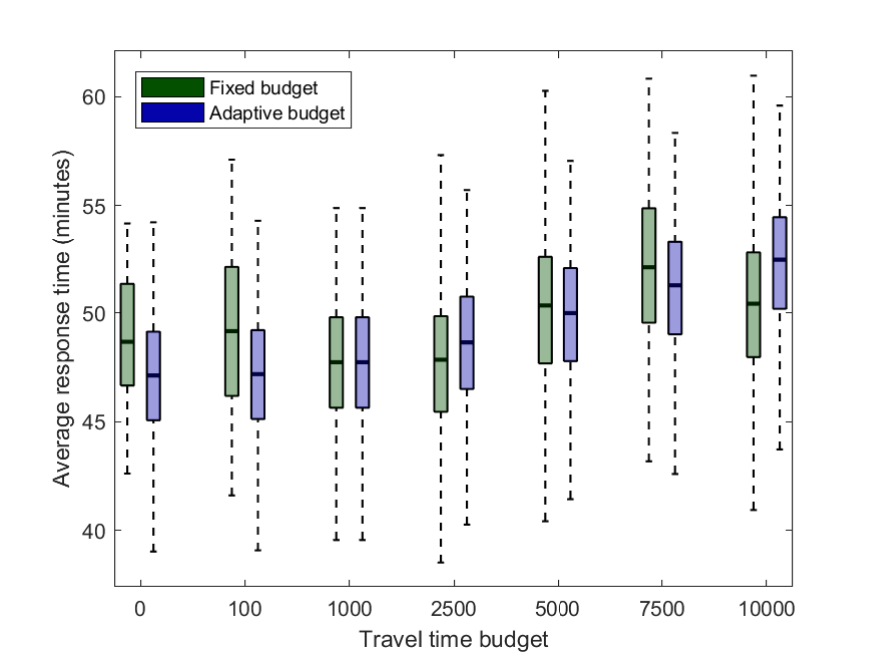}}
\subfigure[\ Overnight \label{BS1}]{
\includegraphics[width=.3\textwidth]{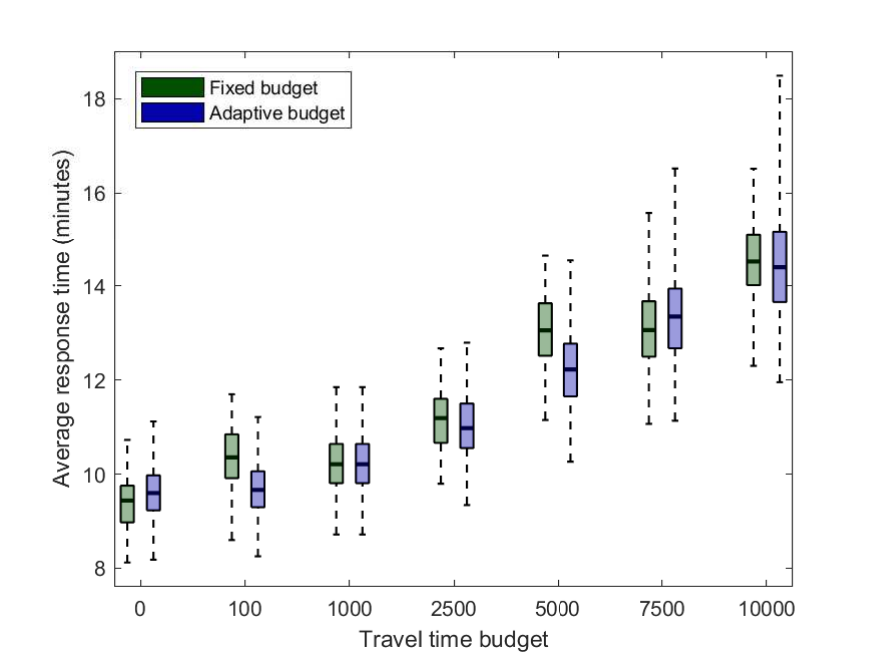}}
\subfigure[\ Weekend \label{BS2}]{
\includegraphics[width=.3\textwidth]{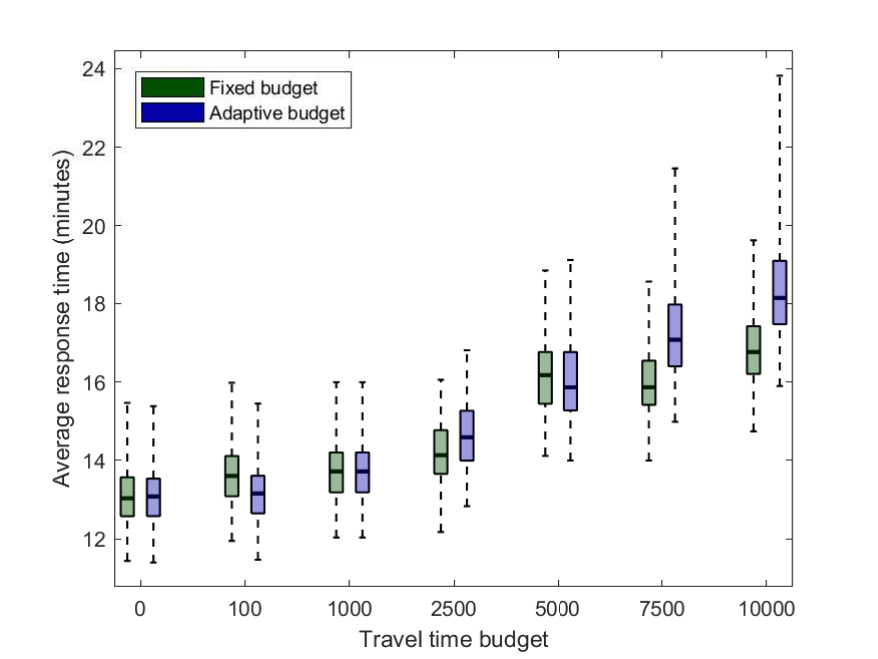}}
\caption{Response time performance of outpost locations determined using a fixed budget of 1000 seconds versus outpost locations determined using an adaptive travel time budget.\label{BSens}}
\end{figure}

Figure~\ref{BSens} compares the response time performance between a fixed travel time budget and a problem-specific travel time budget for each of the three temporal combinations. For rush hour, the fixed budget outpost locations perform better when evaluated on networks with budgets of 2500 and 10000 with an average improvement of 1.45 minutes (4.5\%). For all other instances, the fixed budget performed worse with an average degradation of 1.29 minutes (4\%). Similar results are observed for overnight and weekend baseline traffic scenarios.

\subsubsection{Discussion and policy implications.}
Our results suggest that the outpost locations determined using a travel time budget of 1000 seconds are relatively insensitive to changes the travel time budget. This is an important result that implies that ambulance providers in Dhaka can use the optimal outpost locations from a budget of 1000 seconds without concern that these locations will perform significantly worse for other travel time budgets.

\subsection{How do the optimization-estimated response times compare to the simulation-estimated response times?}\label{Apd:OptVSim}

\begin{figure}[t]
\centering
\includegraphics[width=.5\textwidth]{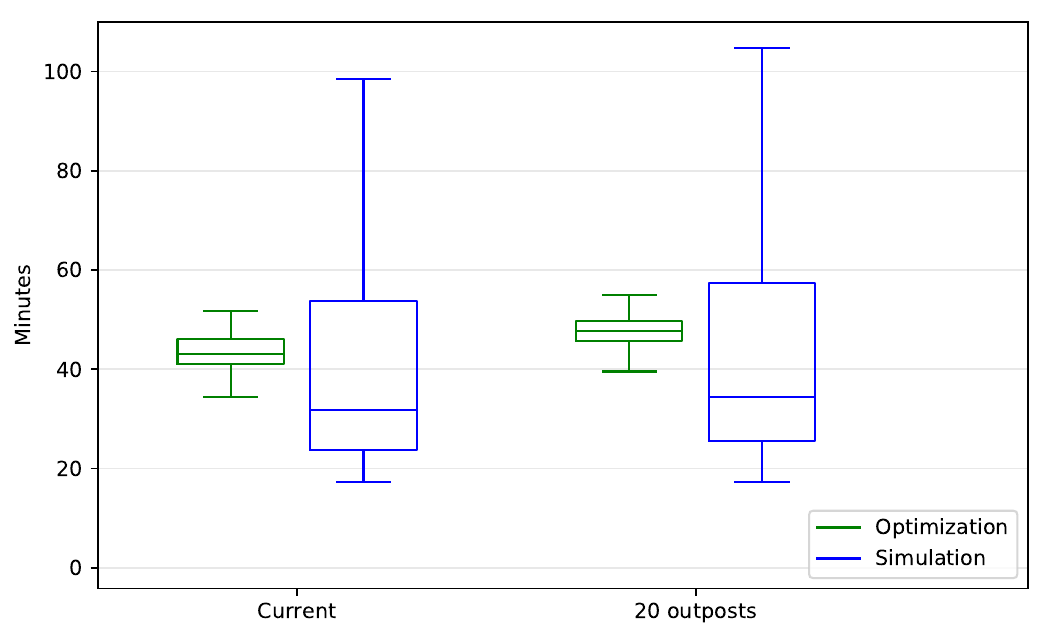}
\caption{A comparison of optimization and simulation estimated response times for the current outpost locations and a network with 20 outpost locations. \label{OptvSim}}
\end{figure}

Figure~\ref{OptvSim} compares the response times estimated by the optimization and simulation models for the current and 20 outpost solutions. The median response time estimated via optimization overestimate the median response time estimated via simulation by 11.3 min (26.2\%) and 13.3 min (27.7\%) for the current and 20 outpost solutions, respectively. Although the optimization results fall within the interquartile range of the simulation results, the optimization model underestimates the total range of response times by 208.5 min and 311.6 min for the current and 20 outpost solutions, respectively. In summary, we find that the response times estimated by the optimization model provide a conservative estimate on the median response time, but significantly underestimate the total range of response times.

\end{document}